\documentclass[12pt]{article}
\usepackage{theorem}
\usepackage{amssymb}
\usepackage{graphicx}
\usepackage[mathscr]{eucal}
\newtheorem{prop}{}[section]
\newtheorem{defi}[prop]{}
\newtheorem{lemma}[prop]{}
{\theorembodyfont{\upshape} \newtheorem{rema}[prop]{}}
\newcommand{\boma}[1]{{\mbox{\boldmath $#1$} }}

\begin{document}
\newcommand{\uper}[1]{\stackrel{\barray{c} {~} \\ \mbox{\footnotesize{#1}}\farray}{\longrightarrow} }
\newcommand{\nop}[1]{ \|#1\|_{\piu} }
\newcommand{\no}[1]{ \|#1\| }
\newcommand{\nom}[1]{ \|#1\|_{\meno} }
\newcommand{\UU}[1]{e^{#1 \AA}}
\newcommand{\UD}[1]{e^{#1 \Delta}}
\newcommand{\bb}[1]{\mathbb{{#1}}}
\newcommand{\HO}[1]{\bb{H}^{{#1}}}
\newcommand{\Hz}[1]{\bb{H}^{{#1}}_{\zz}}
\newcommand{\Hs}[1]{\bb{H}^{{#1}}_{\ss}}
\newcommand{\Hg}[1]{\bb{H}^{{#1}}_{\gg}}
\newcommand{\HM}[1]{\bb{H}^{{#1}}_{\so}}
\def\tvainf{\vspace{-0.4cm} \barray{ccc} \vspace{-0,1cm}{~}
\\ \vspace{-0.2cm} \longrightarrow \\ \vspace{-0.2cm} \scriptstyle{T \vain + \infty} \farray}
\def\Mm{\mathscr M}
\def\XS{\boma{x}}
\def\TS{\boma{t}}
\def\Lam{\boma{\eta}}
\def\DS{\boma{\rho}}
\def\KS{\boma{k}}
\def\LS{\boma{\lambda}}
\def\PR{\boma{p}}
\def\VS{\boma{v}}
\def\ski{\! \! \! \! \! \! \! \! \! \! \! \! \! \!}
\def\h{L}
\def\EM{M}
\def\EMP{M'}
\def\R{R}
\def\Rr{\Upsilon}
\def\E{E}
\def\FFf{\mathscr{F}}
\def\A{F}
\def\Xim{\Xi_{\meno}}
\def\Ximn{\Xi_{n-1}}
\def\lan{\lambda}
\def\om{\omega}
\def\Om{\Omega}
\def\Sim{\Sigm}
\def\Sip{\Delta \Sigm}
\def\Sigm{{\mathscr{S}}}
\def\Ki{{\mathscr{K}}}
\def\zz{{\scriptscriptstyle{0}}}
\def\ss{{\scriptscriptstyle{\Sigma}}}
\def\gg{{\scriptscriptstyle{\Gamma}}}
\def\so{\ss \zz}
\def\Dz{\bb{\DD}'_{\zz}}
\def\Ds{\bb{\DD}'_{\ss}}
\def\Dg{\bb{\DD}'_{\gg}}
\def\Ls{\bb{L}^2_{\ss}}
\def\Lg{\bb{L}^2_{\gg}}
\def\bF{{\bb{V}}}
\def\Fz{\bF_{\zz}}
\def\Fs{\bF_\ss}
\def\Fg{\bF_\gg}
\def\Pre{P}
\def\UUU{{\mathcal U}}
\def\fiapp{\phi}
\def\PU{P1}
\def\PD{P2}
\def\PT{P3}
\def\PQ{P4}
\def\PC{P5}
\def\PS{P6}
\def\Q{Q1}
\def\X{Q2}
\def\Xp{Q3}
\def\Vi{V}
\def\bVi{\bb{V}}
\def\K{V}
\def\Ks{\bb{\K}_\ss}
\def\Kz{\bb{\K}_0}
\def\KM{\bb{\K}_{\, \so}}
\def\HGG{\bb{H}^\G}
\def\HG{\bb{H}^\G_{\so}}
\def\EG{{\mathfrak{P}}^{\G}}
\def\G{G}
\def\de{\| f_0 \|}
\def\esp{\sigma}
\def\dd{\displaystyle}
\def\LP{\mathfrak{L}}
\def\dive{\mbox{div}}
\def\la{\langle}
\def\ra{\rangle}
\def\um{u_{\meno}}
\def\uv{\mu_{\meno}}
\def\Fp{ {\textbf F_{\piu}} }
\def\Ff{ {\textbf F} }
\def\Fm{ {\textbf F_{\meno}} }
\def\piu{\scriptscriptstyle{+}}
\def\meno{\scriptscriptstyle{-}}
\def\omeno{\scriptscriptstyle{\ominus}}
\def\Tt{ {\mathscr T} }
\def\Xx{ {\textbf X} }
\def\Yy{ {\textbf Y} }
\def\Ee{ {\textbf E} }
\def\VP{{\mbox{\tt VP}}}
\def\CP{{\mbox{\tt CP}}}
\def\cp{$\CP(f_0, t_0)\,$}
\def\cop{$\CP(f_0)\,$}
\def\copn{$\CP_n(f_0)\,$}
\def\vp{$\VP(f_0, t_0)\,$}
\def\vop{$\VP(f_0)\,$}
\def\vopn{$\VP_n(f_0)\,$}
\def\vopdue{$\VP_2(f_0)\,$}
\def\leqs{\leqslant}
\def\geqs{\geqslant}
\def\mat{{\frak g}}
\def\tG{t_{\scriptscriptstyle{G}}}
\def\tN{t_{\scriptscriptstyle{N}}}
\def\TK{t_{\scriptscriptstyle{K}}}
\def\CK{C_{\scriptscriptstyle{K}}}
\def\CN{C_{\scriptscriptstyle{N}}}
\def\CG{C_{\scriptscriptstyle{G}}}
\def\CCG{{\mathscr{C}}_{\scriptscriptstyle{G}}}
\def\tf{{\tt f}}
\def\ti{{\tt t}}
\def\ta{{\tt a}}
\def\tc{{\tt c}}
\def\tF{{\tt R}}
\def\C{{\mathscr C}}
\def\P{{\mathscr P}}
\def\V{{\mathscr V}}
\def\TI{\tilde{I}}
\def\TJ{\tilde{J}}
\def\Lin{\mbox{Lin}}
\def\Hinfc{ H^{\infty}(\reali^d, \complessi) }
\def\Hnc{ H^{n}(\reali^d, \complessi) }
\def\Hmc{ H^{m}(\reali^d, \complessi) }
\def\Hac{ H^{a}(\reali^d, \complessi) }
\def\Dc{\DD(\reali^d, \complessi)}
\def\Dpc{\DD'(\reali^d, \complessi)}
\def\Sc{\SS(\reali^d, \complessi)}
\def\Spc{\SS'(\reali^d, \complessi)}
\def\Ldc{L^{2}(\reali^d, \complessi)}
\def\Lpc{L^{p}(\reali^d, \complessi)}
\def\Lqc{L^{q}(\reali^d, \complessi)}
\def\Lrc{L^{r}(\reali^d, \complessi)}
\def\Hinfr{ H^{\infty}(\reali^d, \reali) }
\def\Hnr{ H^{n}(\reali^d, \reali) }
\def\Hmr{ H^{m}(\reali^d, \reali) }
\def\Har{ H^{a}(\reali^d, \reali) }
\def\Dr{\DD(\reali^d, \reali)}
\def\Dpr{\DD'(\reali^d, \reali)}
\def\Sr{\SS(\reali^d, \reali)}
\def\Spr{\SS'(\reali^d, \reali)}
\def\Ldr{L^{2}(\reali^d, \reali)}
\def\Hinfk{ H^{\infty}(\reali^d, \KKK) }
\def\Hnk{ H^{n}(\reali^d, \KKK) }
\def\Hmk{ H^{m}(\reali^d, \KKK) }
\def\Hak{ H^{a}(\reali^d, \KKK) }
\def\Dk{\DD(\reali^d, \KKK)}
\def\Dpk{\DD'(\reali^d, \KKK)}
\def\Sk{\SS(\reali^d, \KKK)}
\def\Spk{\SS'(\reali^d, \KKK)}
\def\Ldk{L^{2}(\reali^d, \KKK)}
\def\Knb{K^{best}_n}
\def\sc{\cdot}
\def\k{\mbox{{\tt k}}}
\def\x{\mbox{{\tt x}}}
\def\g{ {\textbf g} }
\def\QQQ{ {\textbf Q} }
\def\AAA{ {\textbf A} }
\def\gr{\mbox{gr}}
\def\sgr{\mbox{sgr}}
\def\loc{\mbox{loc}}
\def\PZ{{\Lambda}}
\def\PZAL{\mbox{P}^{0}_\alpha}
\def\epsilona{\epsilon^{\scriptscriptstyle{<}}}
\def\epsilonb{\epsilon^{\scriptscriptstyle{>}}}
\def\lgraffa{ \mbox{\Large $\{$ } \hskip -0.2cm}
\def\rgraffa{ \mbox{\Large $\}$ } }
\def\restriction{\upharpoonright}
\def\M{{\scriptscriptstyle{M}}}
\def\m{m}
\def\Fre{Fr\'echet~}
\def\I{{\mathcal N}}
\def\ap{{\scriptscriptstyle{ap}}}
\def\fiap{\varphi_{\ap}}
\def\dfiap{{\dot \varphi}_{\ap}}
\def\DDD{ {\mathfrak D} }
\def\BBB{ {\textbf B} }
\def\EEE{ {\textbf E} }
\def\GGG{ {\textbf G} }
\def\TTT{ {\textbf T} }
\def\KKK{ {\textbf K} }
\def\HHH{ {\textbf K} }
\def\FFi{ {\bf \Phi} }
\def\GGam{ {\bf \Gamma} }
\def\sc{ {\scriptstyle{\bullet} }}
\def\a{a}
\def\ep{\epsilon}
\def\c{\kappa}
\def\parn{\par\noindent}
\def\teta{M}
\def\elle{L}
\def\ro{\rho}
\def\al{\alpha}
\def\si{\sigma}
\def\be{\beta}
\def\ga{\gamma}
\def\te{\vartheta}
\def\ch{\chi}
\def\et{\eta}
\def\complessi{{\bf C}}
\def\len{{\bf L}}
\def\reali{{\bf R}}
\def\interi{{\bf Z}}
\def\naturali{{\bf N}}
\def\To{ {\bf T} }
\def\Td{ {\To}^d }
\def\Zd{ \interi^d }
\def\Zet{{\mathscr{Z}}}
\def\Ze{\Zet^d}
\def\T1{{\textbf To}^{1}}
\def\es{s}
\def\ee{{E}}
\def\FF{\mathcal F}
\def\FFu{ {\textbf F_{1}} }
\def\FFd{ {\textbf F_{2}} }
\def\GG{{\mathcal G} }
\def\EE{{\mathcal E}}
\def\KK{{\mathcal K}}
\def\PP{{\mathcal P}}
\def\PPP{{\mathscr P}}
\def\PN{{\mathcal P}}
\def\PPN{{\mathscr P}}
\def\QQ{{\mathcal Q}}
\def\J{J}
\def\Np{{\hat{N}}}
\def\Lp{{\hat{L}}}
\def\Jp{{\hat{J}}}
\def\Pp{{\hat{P}}}
\def\Pip{{\hat{\Pi}}}
\def\Vp{{\hat{V}}}
\def\Ep{{\hat{E}}}
\def\Gp{{\hat{G}}}
\def\Kp{{\hat{K}}}
\def\Ip{{\hat{I}}}
\def\Tp{{\hat{T}}}
\def\Mp{{\hat{M}}}
\def\La{\Lambda}
\def\Ga{\Gamma}
\def\Si{\Sigma}
\def\Upsi{\Upsilon}
\def\Gam{\Gamma}
\def\Gag{{\check{\Gamma}}}
\def\Lap{{\hat{\Lambda}}}
\def\Upsig{{\check{\Upsilon}}}
\def\Kg{{\check{K}}}
\def\ellp{{\hat{\ell}}}
\def\j{j}
\def\jp{{\hat{j}}}
\def\BB{{\mathcal B}}
\def\LL{{\mathcal L}}
\def\MM{{\mathcal U}}
\def\SS{{\mathcal S}}
\def\DD{D}
\def\Dd{{\mathcal D}}
\def\VV{{\mathcal V}}
\def\WW{{\mathcal W}}
\def\OO{{\mathcal O}}
\def\RR{{\mathcal R}}
\def\TT{{\mathcal T}}
\def\AA{{\mathcal A}}
\def\CC{{\mathcal C}}
\def\JJ{{\mathcal J}}
\def\NN{{\mathcal N}}
\def\HH{{\mathcal H}}
\def\XX{{\mathcal X}}
\def\XXX{{\mathscr X}}
\def\YY{{\mathcal Y}}
\def\ZZ{{\mathcal Z}}
\def\CC{{\mathcal C}}
\def\cir{{\scriptscriptstyle \circ}}
\def\circa{\thickapprox}
\def\vain{\rightarrow}
\def\parn{\par \noindent}
\def\salto{\vskip 0.2truecm \noindent}
\def\spazio{\vskip 0.5truecm \noindent}
\def\vs1{\vskip 1cm \noindent}
\def\fine{\hfill $\square$ \vskip 0.2cm \noindent}
\def\ffine{\hfill $\lozenge$ \vskip 0.2cm \noindent}
\newcommand{\rref}[1]{(\ref{#1})}
\def\beq{\begin{equation}}
\def\feq{\end{equation}}
\def\beqq{\begin{eqnarray}}
\def\feqq{\end{eqnarray}}
\def\barray{\begin{array}}
\def\farray{\end{array}}
\makeatletter \@addtoreset{equation}{section}
\renewcommand{\theequation}{\thesection.\arabic{equation}}
\makeatother
\begin{titlepage}
{~}
\vspace{-2cm}
\begin{center}
{\huge On approximate solutions of semilinear evolution equations
II. Generalizations, and applications to Navier-Stokes equations.}
\end{center}
\vspace{0.5truecm}
\begin{center}
{\large
Carlo Morosi${}^1$, Livio Pizzocchero${}^2$} \\
\vspace{0.5truecm} ${}^1$ Dipartimento di Matematica, Politecnico
di
Milano, \\ P.za L. da Vinci 32, I-20133 Milano, Italy \\
e--mail: carlo.morosi@polimi.it \\
${}^2$ Dipartimento di Matematica, Universit\`a di Milano\\
Via C. Saldini 50, I-20133 Milano, Italy\\
and Istituto Nazionale di Fisica Nucleare, Sezione di Milano, Italy \\
e--mail: livio.pizzocchero@mat.unimi.it
\end{center}
\begin{abstract}
In our previous paper \cite{uno}, a general framework was outlined to treat the approximate solutions
of semilinear evolution equations; more precisely, a scheme was presented to infer from an
approximate solution the existence (local or global
in time) of an exact solution, and to estimate their distance. In the first half
of the present work the abstract framework of \cite{uno} is extended, so as to be applicable
to evolutionary PDEs whose nonlinearities contain derivatives in the space variables.
In the second half of the paper this
extended framework is applied to the incompressible Navier-Stokes equations, on a torus $\Td$ of
any dimension. In this way a number of results are obtained in the setting of the
Sobolev spaces $\HO{n}(\Td)$, choosing the approximate solutions in a number of different ways.
With the simplest choices we recover local existence of the exact solution for
arbitrary data and external forces, as well as global existence
for small data and forces. With the supplementary assumption
of exponential decay in time for the forces, the same decay law is derived for
the exact solution with small (zero mean) data and forces.
The interval of existence for arbitrary data, the upper bounds on data and forces
for global existence, and all estimates on the exponential decay of the exact solution are derived
in a fully quantitative way
(i.e., giving the values of all the necessary constants; this makes a difference with
most of the previous literature). Nextly, the Galerkin
approximate solutions are considered and precise, still quantitative estimates are
derived for their $\HO{n}$ distance from
the exact solution; these are global in time for small data and forces (with exponential
time decay of the above distance, if the forces decay similarly).
\end{abstract}
\vspace{0.2cm} \noindent
\textbf{Keywords:} Differential equations, theoretical approximation, Navier-Stokes \hfill \parn equations,
Galerkin method.
\par \vspace{0.05truecm} \noindent \textbf{AMS 2000 Subject classifications:} 35A35, 35Q30, 65M60.
\end{titlepage}
\tableofcontents
\setcounter{footnote}{0}
\section{Introduction.}
\label{intr}
This is a continuation of our previous paper \cite{uno} on the approximate solutions of semilinear
Cauchy problems in a Banach space $\Ff$, and on
their use to get fully quantitative estimates on: (i) the interval of existence of the exact solution;
(ii) the distance at any time between the exact and the approximate solution. \parn
In \cite{uno}, we mentioned the potential interest of (i) (ii) in relation to the equations
of fluid dynamics. Here we treat specifically the incompressible Navier-Stokes (NS) equations
on a torus $\Td$ of any dimension $d \geqs 2$, taking for $\Ff$ a Sobolev space of vector fields over
$\Td$. To be precise, we consider the Sobolev space $\HO{n}(\Td) \equiv \HO{n}$
of the "velocity fields" $v : \Td \vain \reali^d$ whose derivatives of order $\leqs n$ are square integrable;
then we choose $\Ff := \HM{n}$, where the subscripts ${~}_{\so}$ indicate the subspace of $\HO{n}$ formed by the
divergence free, zero mean velocity fields $f : \Td \vain \reali^d$
(of course, the condition of zero divergence represents incompressibility;
the mean velocity can always be supposed to vanish, passing to a convenient
moving frame). We always
take $n > d/2$. \parn
The choice of $\Td$ as a space domain allows a rather simple treatment, based
on Fourier analysis; we presume that the results of this paper could be extended to bounded
domains of $\reali^d$, with suitable boundary conditions.
In our notations, the NS Cauchy problem is written
\beq {\dot \varphi}(t) = \Delta \varphi(t) - \LP(\varphi(t) \sc \, \partial \varphi(t)) + \xi(t)~,
\qquad \varphi(0) = f_0~,
\label{ns} \feq
where $f = \varphi(t)$ is the velocity field at time $t$, $\LP$ the
Leray projection on the divergence free vector fields and $\xi(t)$ is
the external forcing at time $t$ (more precisely, what remains of the external force field
after applying $\LP$ and subtracting the mean value).
We can regard \rref{ns} as a realization of the abstract
semilinear Cauchy problem
\beq {\dot \varphi}(t) = \AA \varphi(t) + \PP(\varphi(t), t)~, \qquad \varphi(t_0) = f_0 \label{ans} \feq
where $\AA : f \mapsto \AA f$ is a linear operator and $\PP : (f, t) \mapsto \PP(f,t)$ is a
nonlinear map. Of course, in the NS case $\AA$ is the Laplacian $\Delta$ and
$\PP(f,t) := - \LP(f \sc \, \partial f) + \xi(t)$.
By a standard method, both \rref{ns} and its abstract version \rref{ans} can be
reformulated as a Volterra integral equation, involving the semigroup $(e^{t \AA})_{t \geqs 0}$. \parn
In \cite{uno}, a general setting was proposed for semilinear Volterra problems,
when the nonlinearity $\PP(\cdot, t)$ is a sufficiently smooth map of a Banach space $\Ff$ into itself.
This setting cannot be applied
to Cauchy problems like \rref{ns}. In fact, due to the presence of the derivatives
$\partial f$, the map $f \mapsto \LP(f \sc \, \partial f)$ cannot be
seen as a smooth map of a Sobolev space, say $\HM{n}$, into itself; on the contrary, the above map
is smooth from $\HM{n}$ to $\HM{n-1}$. The external forcing $\xi: t \mapsto \xi(t)$ fits well
to this situation if we require it to be a sufficiently smooth map from $[0,+\infty)$ to
$\HM{n-1}$. \parn
In view of the applications to \rref{ns}, in the first half of the present paper (Sections \ref{prelim}-\ref{quadr})
we extend the abstract
framework of \cite{uno} to the case where, at each time $t$,
$\PP(\cdot, t)$ is a smooth map between $\Ff$ and a larger Banach space $\Fm$. A general scheme to treat approximate
solutions is developed along these lines; this could be applied not only to
\rref{ns}, but also to other evolutionary PDEs (essentially, of parabolic type)
with space derivatives in the nonlinear part. \parn
In the second half of the paper (Sections \ref{inns}-\ref{nume})
we fix the attention on the NS equations, in the framework
of the above mentioned $\HM{n}$ spaces (incidentally, we wish to point out that
other function spaces could be used
to analyse the same equations within our general scheme). \parn
Some technicalities related to either the first or the second
half are presented in Appendices \ref{appez}-\ref{appekdue}. \parn
Of course, there is an enormous literature on NS equations,
their approximation methods and the intervals of existence of the exact solutions:
references \cite{Che} \cite{Cos} \cite{Foj} \cite{Fojt} \cite{Gal} \cite{Fuj} \cite{Katf} \cite{Kat}
\cite{Lem} \cite{Sin} \cite{Tem} are examples including seminal works, classical treatises and
recent contributions. Some differences between the present analysis and
most of the published literature are the following: \parn
(i) Our discussion of the NS approximate solutions is
part of a more general framework, in the spirit of the first half of the paper. \parn
(ii) Our analysis is \textsl{fully quantitative}: any function, numerical constant, etc., appearing in our
estimates on the solutions is given explicitly. In the end, our approach
gives bounds on the interval of existence of the exact NS solution and on its distance
from the approximate solution in terms of fully computable numbers; such computations are
exemplified in a number of cases. \parn
(iii) If compared with other contributions, our approach seems to be more suitable
to derive the existence of global exact solutions from suitable approximate solutions,
under specific conditions (typically, of small initial data); a
comment on this point appears in Remark \ref{remache} (iii). \salto
Hereafter we give more details about the contents of the paper. \salto
\textbf{First half: a general setting for the approximate solutions of \rref{ans}.} We have just mentioned
the assumption $\PP(\cdot, t) : \Ff \mapsto \Fm$. We furtherly suppose
$\AA : \Fp \vain \Fm$ where $\Fp$ is a dense subspace of $\Ff$, to be equipped with the graph norm of $\AA$;
in the end, this gives a triple of spaces $\Fp \subset \Ff \subset \Fm$.
\parn
To go on, we require $\AA$ to generate a semigroup on $\Fm$, with the fundamental regularizing
property $\UU{t} (\Fm) \subset \Ff$ for all $t > 0$. A more precise description of all
these assumptions is given in Section \ref{prelim}: here we suppose, amongst else,
the availability of an upper bound  $\um(t) \in (0,+\infty)$ for
the operator norm of $e^{t \AA}$, regarding the latter as a map from $\Fm$ to $\Ff$. The bound $\um$
is allowed to diverge (mildly) for $t \vain 0^{+}$, an indication that $\UU{t} \Fm \not\subset \Ff$ for $t = 0$:
the precise assumption is $\um(t) = O(1/t^{1-\esp})$ with $0 < \esp \leqs 1$.
In applications to the NS system, $\Ff = \HM{n}$ and $\Ff_{\mp} = \HM{n \mp 1}$; the semigroup
$(e^{t \Delta})$ of the Laplacian has the prescribed regularizing features, with $\esp = 1/2$. \parn
In Section \ref{teoria} we present a general theory of the approximate solutions,
for an abstract Cauchy (or Volterra) problem of the type sketched above. The basic idea
is to associate to any approximate solution $\fiap: [t_0, T) \vain \Ff$ of the problem an integral
\textsl{control inequality} for an unknown function $\RR : [t_0, T) \vain [0,+\infty)$;
this has the form
\beq \EE(t) + \int_{t_0}^t~ d s~ \um(t - s)~ \ell(\RR(s), s) \leqs \RR(t) \label{hsf} \feq
where $\EE : [t_0, T) \vain [0,+\infty)$ is an estimator for the (integral) error of
$\fiap$, and $\ell$ is a function describing the growth of $\PP$ from $\fiap$.
The main result in this framework is the following: if the control inequality is fulfilled
by some function $\RR$ on $[t_0, T)$, then the semilinear Volterra problem has an exact solution
$\varphi : [t_0,T) \vain \Ff$, and
$$ \| \varphi(t) - \fiap(t) \| \leqs \RR(t) $$
for all $t$ in this interval ($\|~\|$ is the norm of $\Ff$). \parn
When $\Fm = \Ff$, we recover from here the framework of \cite{uno}. Similarly to the result of
\cite{uno}, the present theorem about
$\RR$, $\fiap$ and $\varphi$ can be considered as the abstract and unifying
form of many statements, appearing in the literature about specific systems. \parn
The available literature would suggest to prove the above theorem along this path: (i) derive
an existence theorem for $\varphi$ on small intervals; (ii) use
some nonlinear Gronwall lemma to prove that $\| \varphi(t) - \fiap(t) \| \leqs \RR(t)$ on any
interval $[t_0, T') \subset [t_0, T)$
where $\varphi$ is defined; (iii) show the existence of $\varphi$ on
the full domain $[t_0, T)$ of $\RR$ by the following \textsl{reductio ad absurdum}:
if not so, $\| \varphi(t) - \fiap(t) \|$ would
diverge before $T$ and its upper bound via $\RR(t)$ would be violated. \parn
Our proof of the theorem on $\RR$, $\fiap$ and $\varphi$, presented in
Section \ref{prova}, replaces the above strategy with a more constructive
approach. The main idea is to interpret
the control inequality \rref{hsf} as individuating a tube of radius $\RR = \RR(t)$ around
$\fiap$, invariant under the action of the semilinear Volterra operator $\JJ$ for
our problem. This makes possible to construct the solution by an iteration of
Peano-Picard type, starting from $\fiap$; the result is a Cauchy sequence
of functions $\varphi_k = \JJ^k (\fiap)$ on $[t_0, T)$, ($k=0,1,2,...$), whose $k \vain +\infty$ limit
is an exact solution of the given Volterra problem. \parn
From this viewpoint,
existence of the solution on a short time interval, with any datum $\varphi(t_0) = f_0$,
is a very simple corollary of the previous theorem based
the choice $\fiap(t) :=$ constant $= f_0$. \parn
Even though there is a basic analogy with \cite{uno}, proving the main theorem on approximate
solutions is technically more difficult in the present case, mainly due to the divergence of $\um(t)$ for $t \vain 0^{+}$.
Such a divergence is also relevant in applications: in fact, differently from
\cite{uno}, Eq. \rref{hsf} with $\leqs$ replaced by $=$ cannot be
reduced to an ordinary differential equation.
Our assumption $\um(t) = O(1/t^{1 - \esp})$ relates \rref{hsf} to the framework
of singular integral equations of fractional type (which could be interpreted in
terms of the so-called "fractional differential calculus"). \parn
In spite of these pathologies, solving
\rref{hsf} is rather simple when the semigroup $(e^{t \AA})$ and $\fiap$ have suitable features,
and the nonlinear function $\PP$ has the (affine) quadratic structure
\beq \PP(f,t) = \PPP(f,f) + \xi(t)~, \label{quas} \feq
with $\PPP : \Fm \times \Fm \vain \Ff$ a continuous bilinear form
and $\xi : [0,+\infty) \vain \Fm$ a (locally Lipschitz) map; this is the
subject of Section \ref{quadr} (where the datum $f_0$ of \rref{ans} is
always specified at $t_0 = 0$).
\parn
The section starts from a fairly general statement
on the control inequaility \rref{hsf}, which is subsequently applied with specific choices of the
approximate solution.
First of all, we consider the choice $\fiap(t) := 0$. In this case, for any datum $f_0$ and
external forcing $\xi$, we construct for the control inequality a solution $\RR$
with domain a suitable interval $[0,T)$; this implies the existence on $[0,T)$
of the solution $\varphi$ of \rref{ans}, and gives an estimate $\| \varphi(t) \|
\leqs \RR(t)$ on the same interval.
If $f_0$ and $\xi$ are
sufficiently small,  $T = +\infty$ and so $\varphi$ is global.
With the stronger assumption that $\xi(t)$ decays exponentially for $t \vain +\infty$,
we derive for the control inequality a solution $t \mapsto \RR(t)$ which is also exponentially
decaying; so, the same can be said for $\| \varphi(t) \|$.
Next we consider, for a small $f_0$ and a small, exponentially decaying
$\xi$, the approximate solution $\fiap$ obtained solving the linear Cauchy problem
${\dot \fiap}(t) = \AA \fiap(t) + \xi(t)$, $\fiap(0) = f_0$. In this case
the control equation still possesses a global, exponentially decaying solution $\RR$, giving
a precise estimate on the distance $\| \varphi(t) - \fiap(t) \|$.
\salto
\textbf{Second half: applications to the NS equations.}
In Section \ref{inns} we review the Sobolev
spaces of vector fields on $\Td$, and the Leray formulation of the incompressible NS equations
within this framework; furthermore, we show that the Cauchy problem with
mean initial velocity $m_0$ can be reduced to an equivalent Cauchy problem where the initial velocity
has zero mean, by a change of space-time coordinates $(x, t) \mapsto (x - h(t), t)$, where
the function $t \mapsto h(t)$ is suitably determined. \parn
In the same section we give explicitly the constants $K_{n d} \equiv K_n$ such that $\| f \sc \partial g \|_{n-1}
\leqs K_n \| f \|_n \| g \|_n$ for all velocity fields $f, g$ on $\Td$, $\|~\|_n$ and $\|~\|_{n-1}$ denoting
the Sobolev norms of orders $n$ and $n-1$. The study of these constants, inspired by our previous work
\cite{MP}, prepares the fully quantitative application of the methods presented in the first half of the paper. \parn
Section \ref{basic} starts from
the formulation \rref{ns} of the Cauchy problem, in the already mentioned
Sobolev spaces $\Ff = \HM{n}$, $\Ff_{\mp} = \HM{n \mp 1}$. We check that \rref{ns} fulfills
all requirements of the general theory for quadratic nonlinearities, and construct
the estimator $\um$ for the semigroup $(e^{t \Delta})$. \parn
In Section \ref{resns} we rephrase for the NS equations all the results of Section \ref{quadr} on the
abstract quadratic case \rref{quas}. The estimates on the time of existence $T$, for arbitrary data
and forcing, have a fully explicit form; the same happens for the bounds on the norms $\| f_0 \|_n$,
$\| \xi(t) \|_{n-1}$ which ensure global existence and, possibly, exponential decay of
$\varphi(t)$ for $t \vain + \infty$.
\parn
In Section \ref{gale} we discuss the approximate NS solutions provided by the Galerkin method.
More precisely, for each finite set $\G (\not\ni 0)$ of wave vectors we consider the subspace $\HG$
spanned by the exponentials $e^{i k x}$ ($k \in \G$), and the projection on $\HG$ of the NS Cauchy problem;
this has a solution $t \mapsto \varphi^\G(t)$ (in general, on a sufficiently
small interval; with special assumptions, also involving the forcing,
$\varphi^\G$ is global and decays exponentially
for $t \vain +\infty$). Applying the framework of Section \ref{basic} with $\fiap = \varphi^\G$ we
derive the following results (with $p > n$ and $|G | :=
\inf_{k \in \Zd_0 \setminus \G} \sqrt{1 + |k|^2}$). \parn
(i) For any initial datum $f_0 \in \HM{p}$ of the NS Cauchy problem \rref{ns}, and
each external forcing $\xi$ with values in $\HM{p-1}$, both
$\varphi^G$ and $\varphi$ exist on a suitable interval $[0,T)$, and there is an estimate
$$ \| \varphi(t) - \varphi^\G(t) \|_n \leqs {\WW_{n p \, |G|}(t) \over |G|^{p-n}}
\qquad \mbox{for $t \in [0,T)$}~; $$
$T$ can be $+\infty$, if the datum and the forcing
are sufficiently small. Both  $T$ and the
function $t \mapsto \WW_{n p \, |G|}(t)$ are given explicitly. \parn
(ii) If $f_0 \in \HM{p}$ is sufficiently small and there is a small,
exponentially decaying forcing $t \mapsto \xi(t) \in \HM{p-1}$, then
$$ \| \varphi(t) - \varphi^\G(t) \|_n \leqs {W_{n p \, |G|} \over |G|^{p-n}} \, e^{-t} \qquad \mbox{for
$t \in [0,+\infty)$}~; $$
the upper bounds for $f_0$, $\xi$ and the coefficient $W_{n p \, |G|}$ are also
given explicitly. \parn
The results (i) (ii) imply convergence of $\varphi^\G$ to $\varphi$ as $|G| \vain +\infty$,
on the time interval where the previous estimates hold (which can be $[0,+\infty)$,
as pointed out). \parn
In Section \ref{nume} we exemplify our estimates giving the numerical values of $T$ and
of the error estimators in (i) (ii) for certain data and forcing, with $d=3$
and $n=2$, $p=4$~.
\section{Introducing the abstract setting.}
\label{prelim}
\textbf{Notations.} (i) All Banach spaces considered in this paper are over the same
field, which can be $\reali$ or $\complessi$. \parn
(ii) If $\Xx$ and $\Yy$ are Banach spaces, we write
\beq \Xx \hookrightarrow \Yy \feq
to indicate that $\Xx$ is a dense vector subspace of $\Yy$ and that its natural inclusion into
$\Yy$ is continuous (i.e., $\| x \|_{{}_{\footnotesize{\Yy}}} \leqs$
constant $\| x \|_{{}_{\footnotesize{\Xx}}}$ for all  $x \in \Xx$). \parn
(iii) Consider two sets $\Theta, \Xx$ and a function $\chi : \Theta \vain \Xx$, $t \mapsto \chi(t)$. The graph of $\chi$ is
\beq \gr \chi := \{ (\chi(t), t)~|~t \in \Theta \} \subset \Xx \times \Theta ~. \feq
If $\Xx = [0,+\infty]$, we define the subgraph of $\chi$ as
\beq \sgr \chi := \{ (r, t)~|~t \in \Theta, r \in [0, \chi(t)) \} \subset  [0,+\infty) \times \Theta~. \label{subgr} \feq
(iv) Consider a function $\chi : \Theta \vain \Xx$, where $\Theta$ is a real interval and
$\Xx$ a Banach space.
This function is locally Lipschitz if, for each compact subset $I$
of $\Theta$, there is a constant $M = M(I) \in [0,+\infty)$ such that
\beq \| \chi(t) - \chi(t') \|_{{}_{\footnotesize{\Xx}}} \leqs M |t - t'| \qquad \mbox{for all $t, t' \in I$}~. \feq
As usually, we denote with $C^{0,1}(\Theta, \Xx)$ the set of these functions.
\salto
\textbf{General assumptions.} Throughout the section, we will consider a
set
\beq (\Fp, \Ff, \Fm, \AA, u, \um, \PP) \feq
with the following properties.
\vskip 0.1cm\noindent
(\PU) $\Fp$, $\Ff$ and $\Fm$ are Banach spaces with norms $\nop{~}$, $\no{~}$ and $\nom{~}$, such that
\beq \Fp \hookrightarrow \Ff \hookrightarrow \Fm~. \feq
Here and in the sequel, $\BBB(f_0, r)$ will denote the open ball
$\{ f \in \Ff~|~\no{ f - f_0} < r \}$ (the radius $r$ can be $+\infty$, and
in this case $\BBB(f_0, r) = \Ff$). \parn
(\PD) $\AA$ is a linear operator such that
\beq \AA : \Fp \vain \Fm~, \qquad f \mapsto \AA f~.  \feq
Viewing $\Fp$ as a subspace of $\Fm$, the norm $\nop{~}$ is equivalent to the
graph norm $f \in \Fp \mapsto \nom{f} + \nom{A f}$. \vskip 0.1cm\noindent
(\PT) Viewing $\AA$ as a densely defined linear operator in $\Fm$,
it is assumed that $\AA$ generates a strongly continuous semigroup $(e^{t \AA})_{t \in [0,+\infty)}$
on $\Fm$ (of course, from the standard theory of linear semigroups, we have
$e^{t \AA}(\Fp) \subset \Fp$ for all $t \geqs 0$). \vskip 0.1cm\noindent
(\PQ) One has
\beq e^{t \AA}(\Ff) \subset \Ff \qquad \mbox{for $t \in [0, +\infty)$}~; \label{reguf} \feq
the function $(f, t) \mapsto e^{t \AA} f$ gives a strongly continuous semigroup on
$\Ff$ (i.e., it is continuous from $\Ff \times [0,+\infty)$ to $\Ff$). \parn
Furthermore,
$u \in C([0,+\infty), (0,+\infty))$ is a function such that
\beq \no{e^{t \AA} f} \leqs u(t) \no{f} \qquad \mbox{for $t \geqs 0$, $f \in \Ff$}~; \label{equ} \feq
this function will be referred to as an \textsl{estimator} for the semigroup
$(e^{t \AA})$ with respect to the norm of $\Ff$. \vskip 0.1cm\noindent
(\PC) One has
\beq e^{t \AA}(\Fm) \subset \Ff \qquad \mbox{for $t \in(0, +\infty)$}~; \label{regum} \feq
the function $(f, t) \mapsto e^{t \AA} f$
is continuous from $\Fm \times (0,+\infty)$ to $\Ff$
(in a few words: for all $t > 0$, $e^{t \AA}$ regularizes the vectors of $\Fm$, sending them
into $\Ff$ continuously). Furthermore,
$\um \in C((0,+\infty), (0,+\infty))$ is a function such that
\beq \no{e^{t \AA} f} \leqs \um(t) \nom{f} \qquad \mbox{for $t > 0$, $f \in \Fm$}~; \label{eqv} \feq
\beq \um(t) = O({1 \over t^{1- \esp}}) \qquad \mbox{for $t \vain  0^{+}$}~,~~\esp \in (0,1]~.
\label{integrab} \feq
The function $\um$ will be referred to as an estimator for the semigroup $e^{t \AA}$ with respect
to the norms of $\Ff$ and $\Fm$; Eq. \rref{integrab} ensures its integrability in any right neighbourhood
of $t=0$. \vskip 0.1cm\noindent
(\PS) One has
\beq \PP : Dom \PP \subset \Ff \times \reali \vain \Fm~, \qquad (f,t) \mapsto \PP(f,t)~, \feq
and the domain of $\PP$ is semi-open in $\Ff \times \reali$: by this we mean that, for any $(f_0, t_0)
\in Dom \PP$, there are $\delta, r \in (0,+\infty]$ such that $\BBB(f_0, r)
\times [t_0, t_0 + \delta) \subset Dom \PP$. Furthermore,
$\PP$ is Lipschitz on each closed, bounded subset $\CC$ of $\Ff \times \reali$
such that $\CC \subset Dom\PP $; by this, we mean that there are constants $L =
L(\CC)$ and $M = M(\CC) \in [0,+\infty)$ such that
\beq \nom{\PP(f, t) - \PP(f', t')} \leqs L \no{f - f'} + M | t - t' |
\qquad \mbox{for all $(f, t), (f', t') \in \CC$}~. \label{lip} \feq
\begin{rema}
\textbf{Remark.}
As anticipated, our aim is to discuss the Cauchy problem $\dot \varphi(t) = \AA
\varphi(t) + \PP(\varphi(t),t)$, $\varphi(t_0) = f_0$ (and its equivalent
formulation as a Volterra problem) for a system $(\Fp,\Ff, \Fm, \AA, u, \um, \PP)$
with the previously mentioned properties (\PU),...,(\PS). In
comparison with the present work, the analysis of \cite{uno} corresponds to the
special case
\beq \Fm = \Ff~, \qquad \um = u \feq
in which, by the continuity of $u$ at $t=0$, Eq. \rref{integrab} is fulfilled with
$\esp=1$ ({\footnote{In \cite{uno} $e^{t \AA}$ was written ${\mathcal U}(t)$, and
$\Fp$ was simply indicated with $Dom \AA$; furthermore, we assumed $Dom \PP$ to be
open in $\reali \times \Ff$.}}).
\end{rema}
\salto
\textbf{Preliminaries to the analysis of the Cauchy and Volterra problems.} \parn
(i) In the sequel, whenever we consider an interval $[t_0, T)$, we intend
$- \infty < t_0 < T \leqs +\infty$ ({\footnote{In \cite{uno}, we also
considered solutions of the Cauchy or Volterra problems with domain
a closed, bounded interval $[t_0, T]$; the symbol
$[t_0, T|$ was employed to denote an interval of either type.
Here we only consider the first case (semiopen, possibly unbounded), simply to avoid
tedious distinctions.}}). \parn
(ii) Let us consider a function $\om \in C([t_0, T), \Fm)$ and the function
\beq \Om : t \in [t_0, T) \mapsto \Om(t) := \int_{t_0}^t d s \, \UU{(t - s)} \om(s)~. \feq
Using the regularising properties
(\PC) of the semigroup with respect to the spaces $\Ff$ and $\Fm$, one easily proves that
for each fixed $t$ the function $s \mapsto \UU{(t - s)} \om(s)$ belongs to $L^1((t_0, t), dt, \Ff)$ and
$\Om \in C([t_0, T), \Ff)$. \parn
(iii) All the results in (ii) apply in particular to the case $\om(s) := \PP(\psi(s), s)$
where $\psi \in C([t_0, T), \Ff)$ and $\gr \psi \subset Dom \PP$. Functions of this form
will often appear in the forthcoming analysis of Volterra problems. \parn
(iv) Many facts stated in the sequel depend on the basic identity, here recalled for future citation,
\beq \psi(t) = \UU{(t - t_0)} \psi(t_0) + \int_{t_0}^t d s~
\UU{(t - s)} \left[ {\dot \psi}(s) - \AA \psi(s) \right] \label{iden}~, \feq holding
for any function  $\psi \in C([t_0, T), \Fp) \, \cap \, C^1([t_0, T), \Fm)$ and
each $t$ in this interval. \parn
\salto
\textbf{Formal definitions of the Cauchy and Volterra problems.} These definitions are similar to the ones adopted in
\cite{uno}, with slight changes due to the present use of two different spaces $\Ff$, $\Fm$.
\begin{prop}
\label{decau} \textbf{Definition.} Consider a pair $(f_0, t_0) \in Dom \PP$, with
$f_0 \in \Fp$. The \textsl{Cauchy problem} \cp with datum
$f_0$ at time $t_0$ is the following one:
$$ \mbox{\textsl{Find}}~
\varphi \in C([t_0, T), \Fp) \cap C^1([t_0, T), \Fm) ~\mbox{\textsl{such that $\gr
\varphi \subset Dom \PP$ and}} $$
\beq {\dot \varphi(t)} = \AA \varphi(t) + \PP(\varphi(t),t)
\quad \mbox{\textsl{for all} $t \in [t_0, T)$}~, \qquad
\varphi(t_0) = f_0 \label{cau} ~.\feq
\end{prop}
We note that $C([t_0, T), \Fp) \subset C([t_0, T), \Ff)$. This fact,
with the properties of $\AA$ and $\PP$, implies the following: if $\varphi
\in C([t_0, T), \Fp)$ and $\gr \varphi \subset Dom \PP$, the right hand side of
the differential equation in \rref{cau} defines a function in $C([t_0, T), \Fm)$.
\begin{prop}
\label{volter}
\textbf{Definition.} Consider a pair $(f_0, t_0) \in Dom \PP$.
The \textsl{Volterra problem} \vp \, with \textsl{datum $f_0$ at time $t_0$} is the following one:
$$ \mbox{\textsl{Find}}~
\varphi \in C([t_0, T), \Ff) \quad \mbox{\textsl{such that $\gr \varphi \subset Dom \PP$ and }}$$
\beq \varphi(t) = \UU{(t - t_0)}
f_0 + \int_{t_0}^t~ d s~ \UU{(t - s)} \PP(\varphi(s),s) \quad
\mbox{\textsl{for all} $t \in [ t_0, T)$}~. \label{int} \feq
\end{prop}
\begin{prop}
\label{integral}
\textbf{Proposition.} For $(f_0, t_0) \in Dom \PP$ and $f_0 \in \Fp$, we have the following.
\parn
(i) a solution $\varphi$ of \cp is also solution of \vp; \parn
(ii) a solution $\varphi$ of \vp is also a solution of \cp, if $\Fm$ is reflexive.
\end{prop}
\textbf{Proof} It is based on \rref{iden}: see \cite{Caz}.
The derivation of (ii), which is the most technical part, uses
the Lipschitz property (\PS) of $\PP$ and the reflexivity of $\Fm$ to show that
a solution of \vp has the necessary regularity to fulfill \cp.  \fine
\parn
\begin{prop}
\label{unique}
\textbf{Proposition.} (Uniqueness theorem for the Volterra problem). Consider a pair
$(f_0, t_0) \in Dom \PP$, and assume that \vp has
two solutions $\varphi \in C([t_0, T), \Ff)$, $\varphi' \in C([t_0, T'),
\Ff)$. Then
\beq \varphi(t) = \varphi'(t) \qquad \mbox{for $t \in [t_0, \min(T, T'))$}~. \feq
\end{prop}
\textbf{Proof.} We consider any $\tau \in [t_0, \min(T, T'))$, and show that
$\varphi = \varphi'$ in $[t_0, \tau]$. To this purpose,
we subtract Eq.\rref{int} for $\varphi$ from the analogous
equation for $\varphi'$; taking the norm $\no{~}$ and using Eqs.
\rref{eqv} \rref{integrab} \rref{lip}, for each $t \in [t_0, \tau]$ we obtain:
\beq
\no{ \varphi(t) - \varphi'(t) } \leqs \int_{t_0}^t \! \! \! d s \, \um(t - s)
\nom{\PP(\varphi(s), s) - \PP(\varphi'(s), s)}
\label{dfnt}
\feq
$$
 \leqs
U L \int_{t_0}^t \! \! \! d s \, {\no{ \varphi(s) - \varphi'(s) } \over (t - s)^{1 - \esp}}
~.
$$In the above: $L \geqs 0$ is a constant fulfilling the Lipschitz condition
\rref{lip} for $\PP$ on the set $\CC := \gr (\varphi \restriction [t_0, \tau])
\cup \gr(\varphi' \restriction [t_0, \tau])$; $U \geqs 0$ is a constant such that
$\um(t') \leqs U/{t'}^{1-\esp}$ for all $t' \in (0, \tau]$ (which exists due to
\rref{integrab}). \parn
Eq. \rref{dfnt} implies $\no{ \varphi(t) - \varphi'(t) }
= 0$ for all $t \in [t_0, \tau]$; in fact, this result follows applying to the function
$z(t) := \no{ \varphi(t) - \varphi'(t) }$ the forthcoming Lemma. \fine
\begin{prop}
\label{lemmaz}
\textbf{Lemma.} Consider a function $z \in C([t_0, \tau], [0,+\infty))$
(with $-\infty < t_0 < \tau < + \infty$), and assume there are
$\Lambda \in [0,+\infty)$, $\esp \in (0,1]$ such that
\beq z(t) \leqs \Lambda \int_{t_0}^t d s \, {z(s)\over (t - s)^{1 - \esp}} \qquad \mbox{for $t \in [t_0, \tau]$}~.
\label{diseqz} \feq
Then, $z(t) = 0$ for all $t \in [t_0, \tau]$. \parn
\end{prop}
\textbf{Proof.} See Appendix \ref{appez}. \fine
\begin{rema}
\label{locex} \textbf{Remark.} For \vp we will grant existence as well, on
sufficiently small time intervals (see the forthcoming Proposition
\ref{esloc}, where local existence is obtained as a simple application
of the general theory of approximate solutions).
\end{rema}
\salto
\textbf{The Volterra integral operator.} This is the (nonlinear) integral
operator appearing in problem \vp. More precisely, let us state the following.
\begin{prop}
\label{defj} \textbf{Definition.} Let $(f_0, t_0) \in Dom \PP$. The \textsl{Volterra
integral operator}
$\JJ_{(f_0, t_0)} \equiv \JJ$ associated to this pair is the following map: \parn
(i) $Dom \JJ$ is made of the functions $\psi \in C([t_0, T), \Ff)$
(with  arbitrary $T \in (t_0, + \infty]$)
such that $\gr \psi \subset Dom \PP$; \parn
(ii) for each $\psi$ in this domain, $\JJ(\psi) \in C([t_0, T), \Ff)$ is the function
\beq t \in [t_0, T) \mapsto \JJ(\psi)(t) := \UU{(t - t_0)}
f_0 + \int_{t_0}^t~ d s~ \UU{(t - s)} \PP(\psi(s),s)~. \label{dej} \feq
\end{prop}
\section{Approximate solutions of the Volterra and Cau\-chy problems: the main
result.}
\label{teoria}
Throughout the section, we
consider again a set $(\Fp, \Ff, \Fm, \AA, u, \um, \PP)$, with the properties (\PU)-(\PS)
of the previous section.
The definitions that follow generalize similar notions, introduced in \cite{uno}.
\salto
\textbf{Approximate solutions, and their errors.} We introduce them in the following way.
\begin{prop}
\label{allap}
\textbf{Definition.} Let $(f_0, t_0) \in Dom \PP$. \parn
(i) An \textsl{approximate solution} of \vp is any function $\fiap
\in C([t_0, T), \Ff)$ such that $\gr \fiap \subset
Dom \PP$. \parn
(ii) The \textsl{integral error} of $\fiap$ is the function
\beq E(\fiap) := \fiap - \JJ(\fiap) \in C([t_0, T), \Ff)~, \label{ier} \feq
i.e., $E(\fiap)(t) = \fiap(t) - \UU{(t - t_0)} f_0 -
\int_{t_0}^t d s~\UU{(t - s)} \PP(\fiap(s), s)$. \parn
An \textsl{integral error estimator} for $\fiap$ is a
function $\EE \in  C([t_0, T), [0, +\infty))$ such that, for all $t$ in this interval,
\beq \no{E(\fiap)(t)} \leqs \EE(t)~. \label{defee} \feq
\end{prop}
\begin{prop}
\label{allapp}
\textbf{Definition.} Let $(f_0, t_0) \in Dom \PP$, and $f_0 \in \Fp$. \parn
(i) An \textsl{approximate solution} of \cp is any function $\fiap
\in C([t_0, T), \Fp)$ $\cap$ $C^1([t_0, T), \Fm)$ such that $\gr \fiap \subset
Dom \PP$. \parn
(ii) The \textsl{datum error} for $\fiap$ is the difference
\beq d(\fiap) := \fiap(t_0) - f_0 \in \Fp \subset \Ff~; \feq
a \textsl{datum error estimator} for $\fiap$ is a nonnegative real
number $\delta$ such that
\beq \no{ d(\fiap) } \leqs \delta~. \label{defde} \feq
(iii) The \textsl{differential error} of $\fiap$ is the function
\beq e(\fiap) \in C([t_0, T), \Fm),~~ t \mapsto e(\fiap)(t) := \dfiap(t) -
\AA \fiap(t) - \PP(\fiap(t), t)~; \feq
a \textsl{differential error estimator} for $\fiap$ is a function $\ep \in C([t_0, T), [0,+\infty))$
such that, for $t$ in this interval,
\beq \nom{e(\fiap)(t)} \leqs \ep(t)~. \label{defep}  \feq
\end{prop}
\begin{rema}
\textbf{Remarks.} {\bf (i)} A function $\fiap$ as in Definition \ref{allap} (resp., Definition \ref{allapp})
is a solution of \vp (resp., of \cp)  if and only if $E(\fiap) = 0$
(resp., $d(\fiap) =0$ and $e(\fiap)=0$). \parn
{\bf (ii) } Of course, the previous definitions of the error estimators can be fulfilled setting
$\EE(t) := \no{E(\fiap)(t)}$, $\delta := \no{d(\fiap)}$,
$\ep(t) := \nom{e(\fiap)(t)}$.
\end{rema}
\begin{prop}
\label{lap}
\textbf{Lemma.} Let $(f_0, t_0) \in Dom \PP$, $f_0 \in \Fp$, and
$\fiap$ be an approximate solution of \cp with datum
and differential errors $d(\fiap)$, $e(\fiap)$. Then:  \parn
(i) $\fiap$ is also an approximate solution of \vp, with integral error
\beq E(\fiap)(t) = \UU{(t - t_0)} \, d(\fiap) + \int_{t_0}^t d s~
\UU{(t - s)} e(\fiap)(s)~. \label{Ede}
\feq
(ii) If $\delta, \ep$ are datum and differential error estimators for $\fiap$, an integral error estimator
for $\fiap$ is
\beq \EE(t) := u(t - t_0) \,\delta + \int_{t_0}^t d s~\um(t - s) \ep(s) \qquad \mbox{for all
$t \in [t_0, T)$.} \label{hasthe} \feq
\end{prop}
\textbf{Proof.} (i) To derive Eq.\rref{Ede}, use the definitions of
$E(\fiap)$, $d(\fiap)$, $e(\fiap)$ and the identity \rref{iden} with $\psi := \fiap$. \parn
(ii) To derive the estimator \rref{hasthe}, apply $\no{~}$ to both sides of \rref{Ede},
using Eqs. \rref{equ}, \rref{eqv} for $u, \um$, and Eqs. \rref{defde}, \rref{defep}
for $\delta, \ep$. \fine
\salto
\textbf{Growth of $\boma{\PP}$ from a curve.} To introduce this notion, we need some notations. \parn
Let us
consider a function $\rho \in C([t_0, T), (0,+\infty])$; we recall
that, according to \rref{subgr}, the subgraph of $\rho$ is the set
$\sgr \rho := \{ (r, t)~|~t \in [t_0, T)~,~r \in [0, \rho(t))\}$.
Furthermore, let $\fiapp \in C([t_0, T), \Ff)$.
We define the \textsl{$\rho$-tube around $\fiapp$} as the set
\beq \Tt(\fiapp, \rho) := \{ (f, t)~|~f \in \Ff, t \in [t_0, T),  \no{f - \fiapp(t)} < \rho(t) \}~; \feq
of course, the above tube is the whole space $\Ff$ if $\rho(t) = +\infty$ for all $t$
({\footnote{In \cite{uno}, this notion was presented in the case $\rho =$ constant; the
present generalization is harmless, and could have been employed in our previous work as well.}}). \parn
\begin{prop}
\label{grow} \textbf{Definition.} Let $\fiapp \in C([t_0, T),
\Ff)$, with $\gr \fiapp \subset Dom \PP$. A \textsl{growth
estimator for $\PP$ from $\fiapp$} is a function $\ell$ with these
features. \parn
(i) The domain of $\ell$ is the subgraph of some
function $\rho \in C([t_0, T),(0,+\infty])$, and \beq \ell \in
C(\sgr \rho, [0, +\infty))~, \qquad (r, t) \mapsto \ell(r, t)~;
\label{cont} \feq $\ell$ is nondecreasing in the first variable,
i.e., $\ell(r, t) \leqs \ell(r', t)$ for $r \leqs r'$ and any $t$.
\parn
(ii) The function $\rho$ in (i) is such that  $\Tt(\fiapp, \rho) \subset Dom \PP$. For all
$(f, t) \in \Tt(\fiapp, \rho)$, it is \beq \hspace{3.5cm}
\nom{\PP(f, t) - \PP(\fiapp(t), t) } \leqs \ell(\no{f - \fiapp(t)},
t)~.  \label{esti} \feq
\end{prop}
\vskip 0.1cm\noindent
\begin{rema}
\textbf{Remark.} Consider any tube $\Tt(\fiapp, \rho') \subset Dom \PP$. Using the Lipschitz
property \rref{lip} of $\PP$, one can easily construct a growth estimator $\ell$ of domain
$\sgr(\rho'/2)$, depending linearly on $r$: $\ell(r, t) = \lambda(t) r$.  \end{rema}
\salto
\textbf{The main result on approximate solutions.} This is contained in the following
\begin{prop}
\label{main}
\textbf{Proposition.} Let $(f_0, t_0) \in Dom \PP$, and consider the problem \vp. Suppose that: \parn
(i) $\fiap \in C([t_0, T), \Ff)$ is an approximate solution of \vp,
$\EE \in C([t_0, T),$ $[0,+\infty))$ is an estimator for the integral error $E(\fiap)$;
\parn
ii) $\ell \in C(\sgr \rho, [0, +\infty))$ is a growth estimator for
$\PP$ from $\fiap$ (for a suitable $\rho \in$ $ C([t_0, T), (0,+\infty])$. \parn
Consider the following problem:
$$ \mbox{\textsl{Find}}~
\RR \in C([t_0, T), [0, +\infty)) \quad \mbox{\textsl{such that $\gr \RR \subset \sgr \rho$,~ and}} $$
\beq \EE(t) + \int_{t_0}^t~ d s~ \um(t - s)~ \ell(\RR(s), s) \leqs \RR(t)
\quad \mbox{\textsl{for} $t \in [ t_0, T)$}~. \label{monod} \feq If \rref{monod}
has a solution $\RR$ on $[t_0, T)$, then \vp has a solution $\varphi$ with the same
domain, and
\beq \no{\varphi(t) - \fiap(t)} \leqs \RR(t) \qquad \mbox{for $t \in [t_0, T)$}~. \feq
The solution $\varphi$ is constructed
by a Peano-Picard iteration of $\JJ$, starting from $\fiap$.
\end{prop}
\textbf{Proof.} See the next section. \fine
\begin{prop}
\textbf{Definition.} Eq.\rref{monod} will be referred to as the \textsl{control inequality}.
\end{prop}
\begin{rema}
\label{runk}
\textbf{Remarks. (i)} It is worthwhile stressing the following: the estimators $\um, \ell, \EE$
in the control inequality \rref{monod} depend on $\AA, \PP, \fiap$, and
should be regarded as known when the Volterra problem and the approximate
solution are specified. So, \rref{monod} is a problem in one unknown $\RR$,
that one sets up using only informations about $\fiap$. After $\RR$ has been
found, it is possible to draw conclusions about the (exact) solution $\varphi$ of \vp. In the usual
language: the control inequality allows predictions on $\varphi$ through an \textsl{a posteriori}
analysis of $\fiap$.
\parn
\textbf{(ii)} (extending to the present framework
a comment in \cite{uno}). Typically, one meets this situation: $\fiap$,
$\EE$, $\ell$ are defined for $t$ in some interval $[t_0, T')$, and
the control inequality \rref{monod} has a solution $\RR$ on an interval $[t_0, T)
\subset [t_0, T')$; in this case one renames $\fiap$, $\EE$, etc. the restrictions of the
previous functions to $[t_0, T)$, and applies
Proposition \ref{main} to them (as an example, this occurs essentially in the proof
of the forthcoming result).
\end{rema}
\salto
\textbf{A first implication of Proposition \ref{main}: local existence.}
The most general and simple consequence of Proposition
\ref{main} is the fact anticipated in Remark \ref{locex}, i.e.,
the local existence for the Volterra problem. Here we formulate this statement precisely.
\begin{prop}
\label{esloc}
\textbf{Proposition.} Let $(f_0, t_0) \in Dom \PP$. Then, there are $R', T', \EE, \ell$ such
that (i)-(iii) hold: \parn
(i) $R' \in (0,+\infty]$, $T' \in (t_0, + \infty]$ and $\BBB(f_0, R') \times [t_0, T') \subset Dom \PP$; \parn
(ii) $\EE \in C([t_0, T'), [0,+\infty))$ and, for all $t \in [t_0, T')$,
\beq \| f_0 - \UU{(t - t_0)} f_0 - \int_{t_0}^t d s~\UU{(t - s)} \PP(f_0, s)~\| \leqs \EE(t)~, \qquad
\EE(t_0) = 0~;  \label{rhs} \feq
(iii) $\ell \in C([0, R') \times [t_0, T'), [0,+\infty))$, $(r, t) \mapsto \ell(r,t)$;
this function is non decreasing in the first variable and, for $(f, t) \in \BBB(f_0, R') \times [t_0, T')$,
\beq \nom{ \PP(f, t) - \PP(f_0, t)} \leqs \ell(\| f - f_0 \|, t)~. \label{eqell} \feq
Given $R', T', \EE, \ell$ with properties (i)-(iii), we have (a)(b): \parn
(a) there are $R \in (0, R')$ and $T \in (t_0, T']$ such that, for all $t \in [t_0, T)$,
\beq \EE(t) + \int_{t_0}^t d s \, \um(t-s) \ell(R, s) \leqs R~; \label{tesa} \feq
(b) if $T$ and $R$ are as in item (a), \vp has a solution $\varphi$ of domain
$[t_0, T)$ and, for all $t$ in this interval,
\beq \| \varphi(t) - f_0 \| \leqs R~. \label{tesb} \feq
\end{prop}
\textbf{Proof.} \textsl{Step 1. Existence of $R',T', \EE, \ell$}.
A pair $(R', T')$ as in (i) exists because $Dom \PP$ is semi-open (see the explanations in
(\PS)). A function $\EE$ as in (ii) can be constructed setting $\EE(t) := $ the left hand side of
the inequality in \rref{rhs}.
A function $\ell$ as in (iii) is constructed putting
$\ell(r,t) := \sup_{f \in \overline{\bf{B}}(f_0, r)} \nom{ \PP(f,t) - \PP(f_0, t)}$; this sup
is proved to be finite using the Lipschitz type inequality \rref{lip}
for $\PP$ on each set $\CC := \overline{\BBB}(f_0, r) \times \{t\}$. One checks
that $\ell(r,t)$ is continuous in $(r,t)$ and non decreasing in $r$. \parn
\textsl{Step 2. Proof of (a)}. Let us pick up any $R \in (0,R')$, and define
\beq \GG : [t_0, T') \vain [0,+\infty)~, ~~t \mapsto \GG(t) := \EE(t) + \int_{t_0}^t d s \, \um(t-s) \ell(R, s)~.
\feq
Then $\GG$ is continuous and $\GG(t_0) = 0$; this fact, with the positivity of $R$,
implies the existence of $T \in [t_0, T')$ such that $\GG(t) \leqs R$ for all $t \in [t_0, T)$;
the last inequality is just the thesis \rref{tesa}. \parn
\textsl{Step 3. Proof of (b)}. We apply Proposition \ref{main} to the approximate solution
\beq \fiap(t) := f_0 \qquad \mbox{for all $t \in [t_0, T)$}~. \feq
Due to (i)-(iii), the function $\EE \restriction [t_0, T)$ is an integral error estimator for $\fiap$, and
the function $\ell \restriction [0, R') \times [t_0, T)$ is a growth estimator for $\PP$ from
$\fiap$ (the function $\rho$ in the general Definition \ref{grow} of growth estimator
is given in this case by $\rho(t)=$ \, constant $= R'$). Eq. \rref{tesa} tells us that
the general control inequality \rref{monod} is fulfilled by the function $\RR(t):=$\,constant $= R$
for all $t \in [t_0, T)$. So, Proposition \ref{main} implies the existence of a solution $\varphi$ of
\vp on $[t_0, T)$, and also gives the inequality \rref{tesb}. \fine
\salto
\section{Proof of Proposition \ref{main}.}
\label{prova}
Let us make all the assumptions in the statement of the Proposition.
We begin the proof introducing an appropriate topology for the space of
continous functions $[t_0, T) \vain \Ff$.
\begin{prop}
\label{ctf}
\textbf{Definition.} From now on, $C([t_0, T), \Ff)$ will be viewed as
a (Hausdorff, complete) locally convex space with the topology of uniform convergence on all compact subintervals
$[t_0, \tau] \subset [t_0, T)$. By this, we mean the
topology induced by the seminorms $(\|~\|_{\tau})_{\tau \in [t_0, T)}$, where
\beq \|~\|_{\tau} : C([t_0, T), \Ff) \vain [0,+\infty), \qquad \psi \mapsto \| \psi \|_{\tau} :=
\sup_{t \in [t_0, \tau]} \| \psi(t) \| \feq
($\no{~}$ is the usual norm of $\Ff$).
\end{prop}
\begin{rema}
\textbf{Remark.} The uncountable family of seminorms $(\|~\|_{\tau})$ is
topologically equivalent to the countable subfamiliy $(\|~\|_{\tau_n})$, where
$(\tau_n)$ is any sequence of points of $[t_0, T)$ such that $\lim_{n \vain +
\infty} \tau_n = T$. Therefore, $C([t_0, T), \Ff)$ is a \Fre space. \ffine
To go on we introduce a basic set, whose definition depends on $\fiap$ and on the function $\RR$
in the control inequality \rref{monod}. \end{rema}
\begin{prop}
\textbf{Definition.} We put
\beq \DDD := \{ \psi \in C([t_0, T), \Ff)~|~\| \psi(t) - \fiap(t) \| \leqs \RR(t) ~~\mbox{for $t
\in [t_0, T)$}~\}~. \feq
\end{prop}
\begin{prop}
\textbf{Lemma.} (i) $\DDD$ is
a closed subset of $C([t_0, T), \Ff)$, in the topology of Definition \ref{ctf}. \parn
(ii) For all $\psi \in \DDD$, one has $\gr \, \psi \subset Dom \PP$ (and so, $\JJ(\psi)$
is well defined).
\end{prop}
\textbf{Proof.} (i) Suppose $\psi \in C([t_0, T), \Ff)$ and
$\psi = \lim_{n \vain \infty} \psi_n$, where $(\psi_{n})$ is
a sequence of elements of $\DDD$. Then, for all $t \in [t_0, T)$ we have
$\| \psi(t) - \fiap(t) \| = \lim_{n \vain \infty} \| \psi_n(t) -
\fiap(t) \| \leqs \RR(t)$. \parn (ii) Let us consider the function
$\rho \in C([t_0, T),(0,+\infty])$ mentioned in the statement of
Proposition \ref{main}. Then, $\psi \in \DDD$ $\Rightarrow$ $\|
\psi(t) - \fiap(t) \| \leqs \RR(t) < \rho(t)$ for all $t \in [t_0,
T)$ $\Rightarrow$ $\gr \psi \subset \Tt(\fiap, \rho) \subset Dom
\PP$. \fine
From now on, our attention will be focused on the map
\beq \DDD \vain C([t_0, T], \Ff)~, \qquad \psi \mapsto \JJ(\psi)~. \feq
Of course, for $\varphi \in \DDD$, we have the equivalence
\beq \varphi~ \mbox{solves \vp} \quad \Longleftrightarrow \quad \JJ(\varphi) = \varphi~.
\label{ofco} \feq
To clarify the sequel, let us recall that $\esp \in (0,1]$ is the constant appearing in Eq.
\rref{integrab}. \parn
\begin{prop}
\label{islip}
\textbf{Lemma.} (i) For each $\tau \in [t_0, T)$ there is a constant $\Lambda_{\tau} \in [0,+\infty)$
such that, for all $\psi, \psi' \in \DDD$,
\beq \| \JJ(\psi)(t) - \JJ(\psi')(t) \| \leqs \Lambda_\tau \int_{t_0}^t d s~ {\| \psi(s) - \psi'(s) \| \over
(t-s)^{1 - \esp}}
\qquad \mbox{for $t \in [t_0, \tau]$}~. \label{verific} \feq
(ii) For all $\tau \in [t_0, T)$ and $\psi, \psi' \in \DDD$, the above equation implies a Lipschitz type
inequality
\beq \| \JJ(\psi) - \JJ(\psi') \|_{\tau} \leqs {\Lambda_\tau (\tau - t_0)^{\esp} \over \esp} \| \psi  - \psi' \|_{\tau}
\label{lypt} \feq
(which ensures, amongst else, the continuity of $\JJ$ on $\DDD$).
\end{prop}
\textbf{Proof.} (i) Let $\psi, \psi' \in \DDD$.
We consider Eq. \rref{dej} for $\JJ(\psi)(t)$, and subtract from it the analogous one
for $\JJ(\psi')(t)$. After applying the norm of $\Ff$ and using \rref{eqv}, we infer
\beq \| \JJ(\psi)(t) - \JJ(\psi')(t) \| \leqs
\int_{t_0}^t ~d s~ \um(t-s) \| \PP(\psi(s),s) - \PP(\psi'(s), s) \|_{-} \label{ddaa} \feq
for all $t \in [t_0, T)$. To go on, we fix $\tau \in [t_0, T)$ and define
\beq C_{\tau} := \{ (f, s)~\in \Ff \times [t_0, \tau]~|~\| f - \fiap(s) \| \leqs \RR(s) \}~. \feq
This is a closed, bounded subset of $\Ff \times \reali$, and $C_{\tau} \subset \Tt(\fiap, \rho) \subset Dom \PP$.
Therefore, by the Lipschitz property (\PS) of $\PP$, there is a nonnegative constant $L(C_\tau) \equiv L_{\tau}$ such
that
\beq \| \PP(f,s) - \PP(f', s) \|_{-}  \leqs L_{\tau} \| f - f' \| \qquad \mbox{for $(f, s)$, $(f',s) \in C_{\tau}$}~.
\label{ppp1} \feq
Furthermore, recalling Eq. \rref{integrab} for $\um$, we see that there is another constant $U_{\tau}$ such that
\beq \um(t') \leqs {U_{\tau} \over {t'}^{1 - \esp}} \qquad \mbox{for $t' \in (0, \tau - t_0]$}~. \label{ppp2} \feq
Inserting Eqs. \rref{ppp1} \rref{ppp2} into \rref{ddaa}, we obtain the thesis \rref{verific} with
$\Lambda_{\tau} := U_{\tau} L_{\tau}$. \parn
(ii) For each $t \in [t_0,\tau]$, Eq. \rref{verific} implies
$$\| \JJ(\psi)(t) - \JJ(\psi')(t) \| \leqs \Lambda_\tau \| \psi - \psi'\|_{\tau}
\int_{t_0}^t {d s \over (t-s)^{1 - \esp}}
 = \Lambda_\tau \| \psi - \psi'\|_{\tau} {(t - t_0)^\esp \over \esp}~.$$
Taking the sup over $t$, we obtain the thesis \rref{lypt}. \fine
\begin{prop}
\label{cent}
\textbf{Lemma} (Main consequence of the control inequality for $\RR$). One has
\beq \JJ(\DDD) \subset \DDD~. \feq
\end{prop}
\textbf{Proof.} Let $\psi \in \DDD$; then
\beq \JJ(\psi) - \fiap = [\JJ(\fiap) - \fiap] + [\JJ(\psi) - \JJ(\fiap)] = - E(\fiap) +
[\JJ(\psi) - \JJ(\fiap)]~. \feq
We write this equality at any $t \in [t_0, T)$, explicitating $\JJ(\psi)(t) - \JJ(\fiap)(t)$;
this gives
\beq \JJ(\psi)(t) - \fiap(t) =  - E(\fiap)(t) + \int_{t_0}^t d s~\UU{(t-s)}~
[~\PP(\psi(s), s) - \PP(\fiap(s), s)~]~. \label{senzanor} \feq Now, we apply the
norm of $\Ff$ to both sides and use Eqs. \rref{defee} for $E(\fiap)(t)$, \rref{eqv} for $\UU{(t-s)}$,
\rref{esti} for the growth of $\PP$ from $\fiap$: in this way we obtain
\beq \| \JJ(\psi)(t) - \fiap(t) \| \leqs \EE(t) + \int_{t_0}^t d s~\um(t - s) \,\ell(\| \psi(s) - \fiap(s) \|, s)~.
\label{eqacc} \feq On the other hand, $\| \psi(s) - \fiap(s) \| \leqs \RR(s)$ implies
$\ell(\| \psi(s) - \fiap(s) \|, s) \leqs \ell(\RR(s), s)$; inserting this into
\rref{eqacc}, and using the control inequality \rref{monod} for $\RR$, we conclude
\beq \hspace{2cm} \| \JJ(\psi)(t) - \fiap(t) \| \leqs \EE(t) + \int_{t_0}^t d s
~\um(t-s) \ell(\RR(s),s) \leqs \RR(t)~, \feq
i.e., $\JJ(\psi) \in \DDD$. \fine
The invariance of $\DDD$ under $\JJ$ is a central result; with the
previously shown properties of $\JJ$, it allows to set up the Peano-Picard iteration
and get ultimately a fixed point of this map.
\begin{defi}
\label{defik}
\textbf{Definition.} $(\varphi_k)$ $(k \in \naturali)$ is the sequence in $\DDD$ defined recursively by
\beq  \hspace{3cm} \varphi_0 := \fiap~, \qquad \varphi_k := \JJ(\varphi_{k-1}) \quad (k \geqs 1)~.  \feq
\end{defi}
\begin{prop}
\label{senzanom}
\textbf{Lemma.} Let $\tau \in [t_0, T)$. For all $k \in \naturali$, one has
\beq \| \varphi_{k+1}(t) - \varphi_{k}(t) \| \leqs \Sigma_{\tau}~{\Lambda_{\tau}^k \, \Gamma(\esp)^k \,
(t - t_0)^{k \esp} \over \Gamma(k \esp + 1)}
\qquad \mbox{for $t \in [t_0, \tau]$}~, \label{num1} \feq
where $\Lambda_\tau$ is the constant of Eq.\rref{verific} and $\Sigma_\tau := \max_{t \in [t_0, \tau]} \EE(t)$. So,
\beq \| \varphi_{k + 1} - \varphi_{k} \|_\tau \leqs \Sigma_{\tau}~{\Theta^k_{\tau \esp} \over \Gamma(k \esp + 1)}~,
\qquad \Theta_{\tau \esp} := \Lambda_{\tau} \, \Gamma(\esp) \, (\tau - t_0)^{\esp}~. \label{eqc} \feq
\end{prop}
\textbf{Proof.} Eq.\rref{eqc} is an obvious consequence of \rref{num1}. We prove
\rref{num1} by recursion, indicating with a subscript ${}_k$ the thesis at a
specified order. We have $\varphi_{1} - \varphi_0 =$ $\JJ(\fiap) - \fiap = -
E(\fiap)$, whence
$\| \varphi_{1}(t) - \varphi_{0}(t) \| \leqs \EE(t) \leqs \Sigma_{\tau}$;
this gives \rref{num1}$_0$. Now, we suppose \rref{num1}$_k$ to hold and infer its
analogue of order $k+1$. To this purpose, we keep in mind Eq. \rref{verific} and write
\beq \| \varphi_{k+2}(t) - \varphi_{k+1}(t) \| =
\| \JJ(\varphi_{k+1})(t) - \JJ(\varphi_{k})(t) \|
\label{sopra} \feq
$$ \leqs \Lambda_{\tau} \int_{t_0}^ t d s {\| \varphi_{k+1}(s) - \varphi_{k}(s) \| \over
(t- s)^{1 - \esp}} \leqs
\Lambda_{\tau} \Sigma_{\tau}~{\Lambda_{\tau}^k \, \Gamma(\esp)^k
\over \Gamma(k \esp + 1)} \int_{t_0}^t d s {(s - t_0)^{k \esp} \over (t-s)^{1 -\esp}}~. $$
On the other hand, we have
({\footnote{To check this, make in the integral the change of variable $s = t_0 + x (t-t_0)$, with
$x \in [0,1]$, and then use the general identity
$$ \int_{0}^1 d x \, x^{\alpha-1} (1- x)^{\beta - 1} = {\Gamma(\alpha) \Gamma(\beta) \over \Gamma(\alpha + \beta)}
\qquad \mbox{for $\alpha, \beta > 0$}~. $$ }})
\beq \int_{t_0}^t~d s~{(s - t_0)^{k \esp} \over (t - s)^{1 - \esp}} =
{\Gamma(k \esp + 1) \Gamma(\esp)\over \Gamma((k+1) \esp + 1)}~(t- t_0)^{(k+1) \esp}~;
\label{hand} \feq
inserting \rref{hand} into \rref{sopra} we obtain the thesis \rref{num1}$_{k+1}$. \fine
The next (and final) Lemma is a generalization of the inequality \rref{eqc}, based on the
Mittag-Leffler function $E_{\esp}$ (see, e.g., \cite{Wei}). For any $\sigma >0$, this is the entire function defined by
\beq E_{\esp} : \complessi \vain \complessi~, \qquad z \mapsto E_{\esp}(z) :=
\sum_{\ell=0}^{+\infty} {z^\ell \over \Gamma(\ell \esp + 1)}~.
\label{mittag} \feq
In particular, $E_{\esp}(z) \in [1,+\infty)$ for all $z \in [0,+\infty)$ and
$E_1(z) = e^z$ for all $z \in \complessi$~.
\begin{lemma}
\label{legen} \textbf{Lemma.} For all $\tau \in [t_0, T)$ and $k, k' \in
\naturali$,
\beq \| \varphi_{k'} - \varphi_k \|_{\tau} \leqs \Sigma_{\tau}~{\Theta^h_{\tau \esp} \over
\Gamma(h \esp + 1)} E_{\esp}(\Theta_{\tau \esp})~,
\qquad h := \min(k, k')
\label{ttesi} \feq
($\Theta_{\tau \esp}$ being defined by \rref{eqc}); this implies that $(\varphi_k)$ is a Cauchy sequence.
\end{lemma}
\textbf{Proof.} It suffices to consider the case $k' > k$ (so that
$h=k$). Writing $\varphi_{k'} - \varphi_k =$ $\sum_{j=k}^{k'- 1} (\varphi_{j+1} -
\varphi_{j}) $ and using Eq.\rref{eqc} we get
\beq \| \varphi_{k'} - \varphi_k \|_{\tau} \leqs \Sigma_{\tau} ~ \sum_{j=k}^{k' -1}
{\Theta^j_{\tau \esp} \over \Gamma(j \esp + 1)}
~. \label{fm} \feq
On the other hand, for each $z \geqs 0$,
\beq \sum_{j=k}^{k'-1} {z^j \over \Gamma(j \esp + 1)} \leqs \sum_{j=k}^{+\infty} {z^j \over \Gamma(j \esp + 1)}
= z^k \sum_{\ell=0}^{+\infty} {z^\ell \over \Gamma((k + \ell) \esp + 1)} \label{lst} \feq
$$ \leqs {z^k \over \Gamma(k \esp + 1)} \sum_{\ell=0}^{+\infty} {z^\ell \over \Gamma(\ell \esp + 1)} =
{z^k \over \Gamma(k \esp + 1)} E_{\esp}(z) $$
(the last inequality depends on the relation
$\Gamma(\alpha + \beta + 1) \geqs \Gamma(\alpha+1) \Gamma(\beta+1)$ for $\alpha, \beta \geqs 0$).
With $z = \Theta_{\tau \esp}$, from \rref{fm} \rref{lst} we obtain \rref{ttesi}. This
equation, with the obvious fact that
$z^h / \Gamma(h \esp + 1)$ $\vain 0$ for $h \vain \infty$ and fixed $z \in \complessi$,
implies
\beq \| \varphi_{k'} - \varphi_k \|_{\tau} \vain 0 \qquad \mbox{for $(k, k') \vain \infty$}~, \feq
for each fixed $\tau \in [t_0, T)$. In conclusion, $(\varphi_k)$ is a Cauchy sequence. \fine
\vskip 0.2cm \noindent
\textbf{Proof of Proposition \ref{main}.} $(\varphi_k)$ being a Cauchy sequence,
$ \lim_{k \mapsto \infty} \varphi_k := \varphi $ exists in $C([t_0, T), \Ff)$;
$\varphi$ belongs to $\DDD$, because this set is closed.
By the continuity of $\JJ$, we have
\beq \JJ(\varphi) = \lim_{k \mapsto \infty} \JJ(\varphi_k) =
\lim_{k \mapsto \infty} \varphi_{k+1} = \varphi~. \feq
Now, recalling \rref{ofco} we get the thesis. \fine
\begin{rema}
\textbf{Remark.}
The Mittag-Leffler function $E_{\esp}$ on $[0,+\infty)$ is strictly related
to a linear integral equation. More precisely, given
$\sigma >0$ let us consider the following problem: find $G \in C([0,+\infty), \reali)$ such that
\beq G(t) = 1 + {1 \over \Gamma(\esp)} \int_{0}^t \! \! d s  {G(s) \over (t - s)^{1- \esp}} \quad
\mbox{for all $t \in [0,+\infty)$.} \label{prob} \feq
This has a unique solution
\beq G(t) := E_{\esp}(t^\esp) \qquad \mbox{for $t \in [0,+\infty)$} ~.\feq
One checks directly that the above $G$ solves \rref{prob}
({\footnote{To this purpose
one inserts into \rref{prob} the series expansion coming from \rref{mittag}, and then uses
the identity in the previous footnote.}}); uniqueness of the solution follows from
the linearity of the problem and from Lemma \ref{lemmaz}. \parn
Integral equations like \rref{prob} are related to the so-called
"fractional differential equations" (see e.g. \cite{Del}, also mentioning $E_{\sigma}$).
\end{rema}
\section{Applications of Proposition \ref{main} to systems with quadratic nonlinearity:
local and global results.}
\label{quadr}
\textbf{The setting.} Throughout this section we consider a set
$(\Fp, \Ff, \Fm, \AA, u, \um, \PPP, \xi)$ with the following features. \parn
$\Fp,\Ff,\Fm$ are Banach spaces, $\AA$ is an operator and $u, \um$ are semigroup estimators fulfilling conditions
(\PU)-(\PC). Furthermore: \parn
(\Q) $\PPP$ is a bilinear map such that
\beq \PPP : \Ff \times \Ff \vain \Fm~, \qquad (f,g) \mapsto \PPP(f,g)~; \feq
we assume continuity of $\PPP$, which is equivalent to the
existence of a constant $K \in (0,+\infty)$ such that
\beq \nom{\PPP(f,g)} \leqs K \| f \| \| g \| \label{defk} \feq
for all $f, g \in \Ff$ ({\footnote{Of course, in the trivial case $\PPP = 0$
we could fulfill \rref{defk} with $K=0$ as well. In the sequel we always assume
$K > 0$, to avoid tedious specifications in many subsequent statements and formulas.
In any case, such statements and formulas could be
extended to $K=0$ by elementary limiting procedures.}}).
 \parn
(\X) We have
\beq \xi \in C^{0,1}([0,+\infty), \Fm)~, \feq
(recall that  $C^{0,1}$ stands for the locally Lipschitz maps).
\salto
Having made the above assumptions, let us fix some notations.
\begin{prop}
\label{defix}
\textbf{Definition.} From now on: \parn
(i) $\MM \in C([0,+\infty), [0,+\infty))$ is a nondecreasing function such that
\beq \int_{0}^t d s \, \um(s) \leqs \MM(t) \qquad \mbox{for $t \in [0,+\infty)$}~,
\qquad  \MM(0) = 0 \label{defmm} \feq
(e.g., $\MM(t) := \int_{0}^t d s \, \um(s)$; in
any case, \rref{defmm} and the positivity of $\um$ imply $\MM(t) > 0$
for all $t > 0$). We
put $\MM(+\infty) := \lim_{t \vain + \infty} \MM(t) \in (0,+\infty]$. \parn
(ii) $\Xim \in C([0,+\infty), [0,+\infty))$ is any nondecreasing function such that
\beq \nom{\xi(t)} \leqs \Xim(t) \qquad \mbox{for $t \in [0,+\infty)$}~ \label{dexi} \feq
(e.g., $\Xim(t) := \sup_{s \in [0,t]} \nom{\xi(s)}$). We put \,
$\Xim(+\infty) := \lim_{t \vain + \infty} \Xim(t) \in [0,+\infty]$. \parn
(iii) $\PP$ is the (affine) quadratic map induced by $\PPP$ and $\xi$ in the following way:
\beq \PP : \Ff \times [0,+\infty) \vain \Fm~, \qquad (f,t) \mapsto \PP(f,t) := \PPP(f,f) + \xi(t)~. \feq
\end{prop}
\salto
\textbf{Properties of $\boma{\PP}$.} Let us analyse this map, so as to match the schemes of
the previous sections.
First of all we note that $Dom \PP$ is semiopen in $\Ff \times \reali$, in the sense
defined in (\PS). Furtheremore, we have the following.
\begin{prop}
\label{samest}
\textbf{Proposition.} For all $(f,t)$ and $(f',t') \in \Ff \times [0,+\infty)$,
\beq \nom{ \PP(f, t) - \PP(f', t') } \leqs 2 K \| f' \| \, \| f - f ' \| + K \| f - f' \|^2 +
\nom{ \xi(t) - \xi(t') }~. \label{kd} \feq
\end{prop}
\textbf{Proof.} For the sake of brevity, we define $h := f - f'$. Then
$$ \PP(f, t) = \PPP(f'+h, f'+h) + \xi(t) $$
$$ = \PPP(f' ,f' ) + \PPP(f', h) + \PPP(h, f') + \PPP(h, h) + \xi(t) ~; $$
subtracting $\PP(f', t')$, we get
$$ \PP(f, t) - \PP(f', t') = \PPP(f', h) +
\PPP(h, f') + \PPP(h, h) + \xi(t) - \xi(t') ~. $$
We apply $\nom{~}$ to both sides, taking into account Eq. \rref{defk}; this gives
$\nom{\PP(f, t) - \PP(f', t')} \leqs 2 K \| f' \| \, \| h \|  + K \| h \|^2 + \nom{\xi(t) - \xi(t')}$,
yielding the thesis \rref{kd}. \fine
The previous proposition has two straightforward consequences.
\begin{prop}
\label{corlip}
\textbf{Corollary.} For each bounded set $\CC$ of $\Ff \times [0,+\infty)$,
there are two constants $L= L(\CC)$, $M= M(\CC)$ such that
\beq \nom{ \PP(f, t) - \PP(f', t') }
\leqs L \| f - f'\| + M | t - t'| \quad \mbox{for $(f,t), (f',t') \in \CC$}~; \feq
so, $\PP$ fulfills condition (\PS).
Let us denote with $\BB$ and $I$ the projections of $\CC$ on $\Ff$ and $[0,+\infty)$, respectively.
Then we can take $L := 4 K \| \BB \|$, where $\| \BB \| := \sup_{f \in \BB} \| f \|$; furthermore,
we can take for $M$ any constant such that $\nom{\xi(t) - \xi(t')} \leqs M | t - t'|$ for $t, t' \in I$.
\end{prop}
\textbf{Proof.} In Eq. \rref{kd}, we substitute the relations $\| f ' \| \leqs \| \BB \|$,
$\|  f - f' \|^2 \leqs (\| f \| + \| f'\|) \| f - f'\| \leqs 2 \| \BB \| \| f - f'\|$,
and the inequality defining $M$. \fine
\begin{prop}
\label{corell}
\textbf{Corollary.} Let us consider any function $\fiapp \in C([0, T), \Ff)$. Then
\beq \nom{ \PP(f, t) - \PP(\fiapp(t), t) } \leqs \ell(\| f - \fiapp(t) \|, t) \qquad \mbox{for all $f \in \Ff$,
$t \in [0,T)$}~, \feq
\beq \ell : [0,+\infty) \times [0, T) \vain [0,+\infty)~,
\quad (r,t) \mapsto \ell(r,t) := 2 K \| \fiapp(t) \| \, r + K r^2~. \label{moregen} \feq
The function $\ell$ is a growth estimator for $\PP$ from $\fiapp$, in the sense of Definition
\ref{grow} (with a radius of the tube $\rho(t) := + \infty$ for all $t$).
\end{prop}
\textbf{Proof.} Use Eq. \rref{kd} with $(f', t') := (\fiapp(t), t)$. \fine
\salto
\textbf{Cauchy and Volterra problems; approximate
solutions.} These problems will always be considered taking $t_0 :=0$ as
the initial time; we will write
\beq \mbox{\vop} := \VP(f_0, 0)~, \quad \mbox{\cop} := \CP(f_0, 0) \qquad \mbox{for each $f_0$}~. \feq
For the above problems, we have the following results. \parn (a) If $f_0 \in \Fp$ and
$\Fm$ is reflexive, \vop is equivalent to \cop (see Proposition \ref{integral}).
\parn (b) For any $f_0 \in \Ff$, uniqueness and local existence are granted for
\vop: see Propositions
\ref{unique} and \ref{esloc}. \parn
(c) We can apply to \vop Proposition \ref{main} on approximate solutions,
choosing arbitrarily the approximate solution $\fiap$; as an error estimator, we can use
at will the function $\ell$ in Corollary \ref{corell} (or any upper bound for it). This
yields the following statement.
\begin{prop}
\label{geprop}
\textbf{Proposition.}
Let $f_0 \in \Ff$ and $T \in (0,+\infty]$. Let us consider for \vop an
approximate solution $\fiap \in C([0,T), \Ff)$, and suppose there are
functions $\EE, \Dd, \RR \in C([0,T),[0,+\infty))$ such that (i)-(iii) hold: \parn
(i) $\fiap$ has the integral error estimate
\beq \| E(\fiap(t)) \| \leqs \EE(t) \qquad \mbox{for $t \in [0,T)$}~; \label{estimatm} \feq
(ii) one has
\beq \| \fiap(t) \| \leqs \Dd(t) \qquad \mbox{for $t \in [0,T)$}~; \label{bdm} \feq
(iii) with $K$ as in \rref{defk} and $\MM$ as in \rref{defmm}, $\RR$ solves the
control inequality
\beq \EE(t)  + K \int_{0}^t \! d s \, \um(t-s) (2 \Dd(s) \RR(s) + \RR^2(s)) \leqs \RR(t)
\qquad \mbox{for $t \in [0,T)$}~. \label{monm} \feq
Then, (a) and (b) hold: \parn
(a) \vop has a solution $\varphi : [0, T) \vain \Ff$; \parn
(b) one has
\beq \| \varphi(t) - \fiap(t) \| \leqs \RR(t) \qquad \mbox{for $t \in [0, T)$}~. \label{axbbm} \feq
\end{prop}
\textbf{Proof.} We refer to the control inequality \rref{monod} in
Proposition \ref{main} (with $t_0 = 0$). Due to Corollary
\ref{corell}, the growth of $\PP$ from $\fiap$ has the quadratic
estimator $\ell(r,t) := 2 K \| \fiap(t) \| \, r + K r^2$; binding
$\| \fiap(t) \|$ via \rref{bdm} we get another estimator, that we
call again $\ell$, of the form \beq \ell(r,t) = 2 K \Dd(t) r + K
r^2~. \feq With this choice of $\ell$, the control inequality
\rref{monod} takes the form \rref{monm} and (a) (b) follow from
Proposition \ref{main}. \fine \salto \textbf{How to handle the
control inequality \rref{monm}.} Rephrasing Remark \ref{runk} in
the present case, we repeat that $\RR$ is the only unkown in
\rref{monm}. In fact, the functions $\EE, \Dd$ appearing therein
can be determined when $\fiap$ is given, and $\um$, $K$ can be
obtained from $\AA$, $\PPP$ (as an example, the computation of
$\um$ and $K$ for the NS equations will be presented in Sections
\ref{inns}-\ref{basic}). Two basic strategies to find a function
$\RR$ solving \rref{monm} on an interval $[0,T)$, if it exists,
can be introduced: \parn (a) the analytical strategy: one makes
some ansatz for $\RR$, substitutes it into \rref{monm} and checks
whether the inequality is fulfilled; \parn (b) the numerical
strategy. \parn Let us mention that a numerical approach was
presented in \cite{uno}, for the simpler control inequality
considered therein; in that case it was possible to transform the
control equality (with $\leqs$ replaced by $=$) into an equivalent
Cauchy problem for $\RR$, and then solve it by a standard package
for ODEs. \parn The approach of \cite{uno} can not be used for
\rref{monm} due to the singularity of $\um(t)$ for $t \vain
0^{+}$; a different numerical attack could be used, but this is
not so simple and its features suggest to treat it extensively
elsewhere. For the above reasons, in the present work we only give
an introductory sketch of the numerical strategy: see Appendix
\ref{outli}. In the rest of the paper, starting from the next
paragraph, we will use the analytical strategy (a). \salto
\textbf{Some special results on \vop.} All these results will be
derived solving the control inequality \rref{monm} by analytic
means, in special cases.
\begin{prop}
\label{esglobbm}
\textbf{Proposition.} Let $f_0 \in \Ff$ and $T \in (0,+\infty]$. Let us consider for \vop an
approximate solution $\fiap \in C([0,T), \Ff)$, and suppose there are
functions $\EE, \Dd \in C([0,T],[0,+\infty))$ such that (i)-(iii) hold: \parn
(i) $\EE$ is nondecreasing and binds the integral error as in \rref{estimatm}; \parn
(ii) $\Dd$ is nondecreasing and binds $\fiap$ as in \rref{bdm}; \parn
(iii) with $K$ as in \rref{defk} and $\MM$ as in \rref{defmm}
{(\footnote{Of course, in the case $T=+\infty$ \rref{coquaam} implies $\MM(+\infty) < +
\infty$}})~,
\beq 2 \sqrt{K \MM(T) \EE(T)} + 2 K \MM(T) \Dd(T) \leqs 1~.  \label{coquaam} \feq
\parn
Then \vop has a solution $\varphi : [0, T) \vain \Ff$ and, for all $t \in [0, T)$,
\beq \| \varphi(t) - \fiap(t) \| \leqs \RR(t)~, \label{xquaam} \feq
$$ \hspace{-0.0cm} \, \RR(t):=
\left\{ \barray{ll}
\! \! \dd{1 \! \! - \! \! 2 K \MM(t) \Dd(t) - \sqrt{(1 \! \! - \! \! 2 K \MM(t) \Dd(t))^2 - 4 K \MM(t) \EE(t)}
\over 2 K \MM(t)}
& \! \!  \mbox{if $t \in (0,T)$,} \\
\! \! \EE(0) & \! \! \mbox{if\, $t = 0$;} \farray \right.
$$
the above prescription gives a well defined, nondecreasing function $\RR$ $\in$ $C([0,T)$, $[0,+\infty))$.
\end{prop}
\textbf{Proof.} We refer to Proposition \ref{geprop}, and try to fulfill
the control inequality \rref{monm} with a \textsl{nondecreasing} $\RR \in C([0,T),[0,+\infty))$.
Noting that $\RR(s) \leqs \RR(t)$, $\Dd(s) \leqs \Dd(t)$ for $s \in [0,t]$, we have
\beq \EE(t)  + K \int_{0}^t \! d s \, \um(t-s) \, (2 \Dd(s) \RR(s) + \RR^2(s))  \label{dician} \feq
$$ \leqs \EE(t) + K (2 \Dd(t) \RR(t) + \RR^2(t)) \int_{0}^t d s \, \um(t-s) $$
$$ \leqs \EE(t) + K (2 \Dd(t) \RR(t) + \RR^2(t)) \, \MM(t)~; $$
the last inequality follows from $\int_{0}^t \! d s \, \um(t-s) = \int_{0}^t \, d s \,
\um(s) \leqs \MM(t)$. Due to \rref{dician}, \rref{monm} holds if
$\EE(t) + K (2 \MM(t) \Dd(t) \RR(t) + \MM(t) \RR^2(t)) \leqs \RR(t)$, i.e.,
\beq K \MM(t) \RR(t)^2 - (1 - 2 K \MM(t) \Dd(t)) \RR(t) + \EE(t) \leqs 0~. \label{nonnegm} \feq
This inequality is fulfilled as an equality if we define
$\RR$ as in \rref{xquaam}, provided that $2 \sqrt{K \MM(t) \EE(t)} + 2 K \MM(t) \Dd(t) \leqs 1$;
this happens for each $t \in [0,T)$ due to the assumption \rref{coquaam} ({\footnote{Obviously
enough, we take for $\RR(t)$ the definition \rref{xquaam} since this gives the smallest nonnegative solution of
Eq. \rref{nonnegm}.}). \parn
The function $\RR$ on $[0,T)$ defined by \rref{xquaam} is continuous and nonnegative; to conclude the
proof, we must check it to be nondecreasing. To this purpose, we note that
\beq \RR(t) = \Rr(K \MM(t), \Dd(t), \EE(t)),~~~\feq
$$ \Rr(\mu, \delta, \ep) :=
\left\{ \barray{ll}
\dd{{1 - 2 \mu \delta -
\sqrt{(1 - 2 \mu \delta)^2 - 4 \mu \epsilon} \over 2 \mu}}
& \mbox{if\, $\mu > 0$}, \\
\epsilon & \mbox{if\, $\mu = 0$.} \farray \right.$$
The above function $\Rr$ has domain
\beq Dom \Rr := \{ (\mu, \delta, \ep)~|~\mu, \delta, \ep \geqs 0,~
2 \mu \delta + 2 \sqrt{\mu \ep}  \leqs 1~\}~, \feq
containing all triples $(K \MM(t), \Dd(t), \EE(t))$
due to \rref{coquaam}; one checks by elementary means (e.g., computing derivatives) that $\Rr$ is a nondecreasing function of
each one of the variables $\mu, \delta, \ep$, when the other two are fixed. \fine
\begin{rema}
\textbf{Remark.} The inequality
\rref{coquaam} is certainly fulfilled if $(\EE(T), \Dd(T))$ or $T$ are sufficiently small (recall that
$\MM(T)$ vanishes for $T \vain 0^{+}$).
\end{rema}
\salto
\textbf{An example: the zero approximate solution.} For simplicity, we suppose
\beq u(t) \leqs 1 \qquad \mbox{for all $t \in [0,+\infty)$}~. \label{withu} \feq
Let $f_0 \in \Ff$, $T \in (0,+\infty]$; we apply Proposition \ref{esglobbm}, choosing for
\vop the trivial approximate solution
\beq \fiap(t) := 0 \qquad \mbox{for $t \in [0,T)$}~. \feq
\begin{prop}
\textbf{Lemma.} $\fiap := 0$ has the integral error
\beq E(\fiap)(t) = - e^{t \AA} f_0 - \int_{0}^t d s \, e^{-(t-s) \AA} \xi(s)~; \label{erzero} \feq
with $\Xim$ as in Definition \ref{dexi}, $\| E(\fiap)(t) \|$ has the estimate
\beq \| E(\fiap)(t) \| \leqs \FFf(t)~,
\qquad \FFf(t) := \de + \Xim(t) \, \MM(t)~. \label{proem} \feq
\end{prop}
\textbf{Proof.} Eq. \rref{erzero} follows from the general definition
\rref{ier} of integral error, and from the fact that $\PP(\fiap(s), s) = \xi(s)$. \parn
Having Eq. \rref{erzero}, we derive \rref{proem} in the following way. First of all,
\beq \| E(\fiap)(t) \| \leqs u(t) \de + \int_{0}^t d s~ \um(t-s) \nom{\xi(s)}~; \feq
but $u(t) \leqs 1$, $\nom{\xi(s)} \leqs \Xim(t)$ for $s \in [0,t]$, so
\beq \| E(\fiap)(t) \| \leqs \de + \Xim(t)\, \int_{0}^t d s~ \um(t-s)
\leqs \de + \Xim(t) \, \MM(t)~. \qquad \square \feq
From the previous Lemma and Proposition \ref{esglobbm}, we infer the following
result.
\begin{prop}
\label{xvxm}
\textbf{Proposition.} With $u$, $\FFf$ as in \rref{withu} \rref{proem}, assume
\beq  4 K \MM(T) \FFf(T) \leqs 1~. \label{coquaazm} \feq
Then \vop has a solution $\varphi : [0, T) \vain \Ff$ and,
for all $t$ in this interval,
\beq \| \varphi(t) \| \leqs \FFf(t) \, \XX(4 K \MM(t) \FFf(t)) ~. \label{xbbzm} \feq
Here $\XX \in C([0,1], [1,2])$ is the increasing function defined by
\beq \XX(z) := \left\{ \barray{ll} \dd{1 - \sqrt{1 - z} \over (z/2)} & \mbox{for $z \in (0,1]$}~, \\
1 & \mbox{for $z = 0$}~. \farray \right.
\label{xquaazm} \feq
\end{prop}
\textbf{Proof.} According to the previous Lemma, we can apply Proposition \ref{esglobbm} with
$\fiap = 0$, $\EE = \FFf$; obviously, we have for
$\|\fiap(t)\|$ the estimator $\Dd(t) := 0$. Eqs. \rref{coquaam},
\rref{xquaam}  yield
the present relations \rref{coquaazm}, \rref{xbbzm}, \rref{xquaazm}; in particular,
the function $\RR$ in \rref{xquaam} is given by $\RR(t) = (1 - \sqrt{1 - 4 K \MM(t) \FFf(t)})$ $/ 2 K
\MM(t) $ $=\FFf(t) \, \XX(4 K \MM(t) \FFf(t))$. \fine
Of course, the basic inequality \rref{coquaazm} is fulfilled if $(f_0, \Xim(T))$ or $T$ are sufficiently small.
\salto
\textbf{Further results (global in time) for \vop, under special assumptions.}
\salto
We keep the assumptions (\PU)-(\PC) of Section \ref{prelim}} and (\Q)(\X) at the beginning of this section, and put
more specific requirements on the semigroup estimators $u, \um$. More precisely, we add
to (\PQ) (\PC) the following conditions, involving two constants
\beq B \in [0, +\infty)~, N \in (0,+\infty)~. \feq
(\PQ') The semigroup estimator $u$ has the form
\beq u(t) = e^{-B t} ~ \mbox{for $t \in [0,+\infty)$}~. \label{espu} \feq
(\PC') The semigroup estimator $\um$ has the form
\beq \um(t) = \uv(t)\,  e^{-B t} \quad \mbox{for $t \in (0,+\infty)$}~, \label{espum} \feq
\beq \uv \in C((0,+\infty), (0,+\infty)), \quad \uv(t) = O({1 \over t^{1-\esp}})
~~\mbox{for $t \vain 0^{+}$} \quad (\esp \in [0,1))~, \feq
\beq \int_{0}^t d s \, \uv(t-s) \, e^{-B s} \leqs N \qquad \mbox{for $t \in [0,+\infty)$}~.\label{defnn} \feq
\begin{prop}
\label{esglobb}
\textbf{Proposition.} Given $f_0 \in \Ff$, let us consider for \vop an
approximate solution $\fiap \in C([0,+\infty), \Ff)$. Suppose there are
constants $E, D  \in [0,+\infty)$ such that:\parn
(i) $\fiap$ admits the integral error estimate
\beq \| E(\fiap)(t) \| \leqs E \, e^{-B t} \qquad \mbox{for $t \in [0, +\infty)$}~; \label{estimat} \feq
(ii) for all $t$ as above,
\beq \| \fiap(t) \| \leqs D \, e^{-B t}~; \label{bd} \feq
(iii) with $N$ as in \rref{defnn} and $K$ as in \rref{defk}, one has
\beq 2 \sqrt{K N E} + 2 K N D  \leqs 1~. \label{coquaa} \feq
Then \vop has a global solution $\varphi : [0, +\infty) \vain \Ff$ and, for all $t \in [0, +\infty)$,
\beq \| \varphi(t) - \fiap(t) \| \leqs \R e^{-B t}~, \label{xquaa} \feq
$$ \R := {1 - 2 K N D - \sqrt{(1 - 2 K N D)^2 - 4 K N E} \over 2 K N}~. $$
\end{prop}
\textbf{Proof.} \textsl{Step 1: the control inequality.} We use again
Proposition \ref{geprop} and the control inequality
\rref{monm}. With the notations of the cited proposition, we have
$\EE(t) := E e^{-B t}$, $\Dd(t) := D \, e^{-B t}$ and $\um$ has the expression
\rref{espum}. So, \rref{monm} takes the form
\beq E e^{-B t} + K \int_{0}^t d s \, \uv(t-s) e^{-B(t-s)} (2 D e^{-B s} \RR(s)
+ \RR^2(s)) \leqs \RR(t)~,  \label{mon} \feq
that we regard as an inequality for an unknown nonnegative function $\RR$.
\parn
\textsl{Step 2: searching for a global solution $\RR$ of \rref{mon}.} We try to fulfill
\rref{mon} with
\beq \RR(t) := \R e^{-B t} \qquad \mbox{for all $t \in [0,+\infty)$}, \feq
with $\R \geqs 0$ an unknown constant. Then, the left hand side of \rref{mon} is
\beq e^{-B t} \left[E + K (2 D \R + \R^2) \int_{0}^t ds \, \uv(t-s) e^{- B s} \right]\feq
$$ \leqs
e^{-B t} \left(E + K (2 D \R + \R^2) N \right)~, $$ where the inequality depends on
\rref{defnn}. The last expression is bounded by $\RR(t)$ if $\R$ fulfills the inequality
$E + K (2 D \R + \R^2) N \leqs \R$, i.e.,
\beq K N \R^2 - (1 - 2 K N D) \R + E \leqs 0~. \label{nonneg} \feq
This condition is fulfilled as an equality if we define $\R$ as in \rref{xquaa};
this $\R$ is well defined and nonnegative due to the assumption \rref{coquaa}
({\footnote{and, in fact, is the smallest
nonnegative solution of \rref{nonneg}.}). Due to the above considerations, the thesis is proved. \fine
\salto
\textbf{Applications to cases with exponentially decaying forcing.}
\salto
From here to the end of the section, we add to (\PU)-(\PC')(\Q)(\X) the following condition: \parn
(\Xp) There is a constant $J \in [0, +\infty)$ such that (with $B$ as in (\PQ')(\PC'))
\beq \nom{\xi(t)} \leqs J e^{-2 B t} \qquad \mbox{for all $t \in [0,+\infty)$}~. \label{eqxp} \feq
Two cases where (\Xp)~ holds are: (i) the trivial case $\xi(t) = 0$ for all $t$; (ii)
situations where the external forcing is switched off in the future. \parn
Hereafter we present two applications of Proposition \ref{esglobb}, corresponding
to different choices for $\fiap$. Both of them give global existence for the exact
solution $\varphi$ of \vop when the datum $f_0$ is sufficiently small, with suitable estimates of the
form \rref{xquaa}.
\salto
\textbf{Example: the zero approximate solution.} Let $f_0 \in \Ff$; we reconsider,
from the viewpoint of Proposition \ref{esglobb}, the \vop approximate solution
\beq \fiap(t) := 0 \qquad \mbox{for $t \in [0,+\infty)$}~. \feq
\begin{prop}
\textbf{Lemma.} $\fiap := 0$ has the integral error estimator
\beq \| E(\fiap)(t) \| \leqs \A e^{-B t}~, \qquad \A := \de + N J~. \label{proe} \feq
\end{prop}
\textbf{Proof.} The integral error $E(\fiap)$ was already computed, see Eq. \rref{erzero}.
From this equation and the assumptions \rref{espu} \rref{espum} on $u$, $\um$ we infer
\beq \| E(\fiap)(t) \| \leqs e^{-B t} \de + \int_{0}^t d s~ e^{-B (t-s)} \uv(t-s) \nom{\xi(s)}~; \feq
inserting here the bound \rref{eqxp} for $\nom{\xi(s)}$, and using Eq. \rref{defnn} for $\uv$ we get
\beq \| E(\fiap)(t) \| \leqs e^{-B t} \de + J e^{-B t} \int_{0}^t d s~ \uv(t-s)  e^{-B s} \leqs
e^{-B t} \de + J e^{-B t} N~. \feq
\fine
From the previous Lemma and Proposition \ref{esglobb}, we infer the following.
\begin{prop}
\label{xvx}
\textbf{Proposition.} With $\A$ as in \rref{proe}, let
\beq 4 K N \A \leqs 1~. \label{coquaaz} \feq
Then \vop has a global solution $\varphi : [0, +\infty) \vain \Ff$ and,
for all $t \in [0, +\infty)$,
\beq \| \varphi(t) \| \leqs \A \, \XX(4 K N \A) \, e^{-B t}~, \label{xbbz} \feq
with $\XX$ as in \rref{xquaazm}.
\end{prop}
\textbf{Proof.} We apply Proposition \ref{esglobb} with $\fiap :=0$. The constants
of the cited proposition are $E=\A$ and $D=0$ (the first equality
follows from the previous lemma, the second one is obvious).
Eqs. \rref{coquaa} \rref{xquaa} yield
the present relations \rref{coquaaz} \rref{xbbz} \rref{xquaazm}; in particular, the
constant in \rref{xquaa} is given by $\R = (1 - \sqrt{1 - 4 K N \A})/(2 K N)$
$= \A \, \XX(4 K N \A)$. \fine
\textbf{Example: the "\boma{\AA}-flow" approximate solution.} Given $f_0 \in \Ff$, we
consider for \vop the approximate solution
\beq \fiap(t) := e^{t \AA} f_0 + \int_{0}^t d s~\UU{(t-s)} \xi(s)
\qquad \mbox{for $t \in [0,+\infty)$} \label{heatff} \feq
(i.e., we use the flow of the linear equation $\dot f = \AA f + \xi(t)$).
\begin{prop}
\textbf{Lemma.} Let $\A$ be as in \rref{proe}. The above $\fiap$ has the integral error estimator
\beq \| E(\fiap)(t) \| \leqs K N \A^2 e^{-B t} \label{eeh} \feq
and fulfills for all $t \in [0,+\infty)$ the norm estimate
\beq \| \fiap(t) \| \leqs \A e^{-B t}~. \label{nest} \feq
\end{prop}
\textbf{Proof.} We first derive Eq. \rref{nest}. To this purpose, we note that
the definition \rref{heatff} of $\fiap$, the assumptions \rref{espu} \rref{espum} on $u, \um$ and Eqs.
\rref{eqxp} for $\nom{\xi(\cdot)}$, \rref{defnn} for $\uv$ imply
\beq \| \fiap(t) \| \leqs e^{-B t} \de + \int_{0}^t d s~e^{-B (t-s)} \uv(t-s) \nom{\xi(s)} \feq
$$ \leqs e^{-B t} \de + J e^{- B t} \int_{0}^t d s~\uv(t-s) e^{- B s} \leqs
e^{- B t} \de + J e^{-B t} N~; $$
by comparison with the definition \rref{proe} of $\A$, we get the thesis \rref{nest}. \parn
Let us pass to the proof of \rref{eeh}. To this purpose we note that the definition
\rref{ier} of integral error gives, in the present case,
\beq E(\fiap)(t) = - \int_{0}^t d s \, e^{(t - s) \AA} \PPP(\fiap(s), \fiap(s))~. \feq
From here and from Eqs. \rref{espum} for $\um$, \rref{defk} for $\PPP$ we infer
\beq
\| E(\fiap)(t) \|  \leqs \int_{0}^t \! \! d s~\uv(t-s) e^{-B (t - s)} \nom{ \PPP(\fiap(s), \fiap(s)) } \feq
$$ \leqs K \int_{0}^t \! \! d s~\uv(t-s) e^{-B (t - s)} \| \fiap(s) \|^2~. $$
In the last inequality, we insert the bound \rref{nest} and then recall \rref{defnn}. This gives
\beq \| E(\fiap)(t) \|  \leqs K \A^2 e^{-B t} \int_{0}^t d s~ e^{- B s} \uv(t-s) \leqs
K N \A^2 e^{-B t}~. \hspace{1cm} \square \feq
Let us return to Proposition \ref{esglobb}; with the previous Lemma, this implies the following result.
\begin{prop}
\label{aflow}
\textbf{Proposition.} Let us keep the definition \rref{proe} for $\A$, and the assumption
\rref{coquaaz} $4 K N \A \leqs 1$.
The global solution $\varphi : [0, +\infty) \vain \Ff$ of \vop is such that,
for all $t \in [0, +\infty)$,
\beq \| \varphi(t) - e^{t \AA} f_0 \| \leqs K N \A^2 \, \XXX(4 K N \A) \, e^{-B t}~; \label{xbbh} \feq
here $\XXX \in C([0,1], [1,4])$ is the increasing function defined by
\beq \XXX(z) := \left\{ \barray{ll} \dd{1 - (z/2) - \sqrt{1 - z} \over (z^2/8)} & \mbox{for $z \in (0,1]$}~, \\
1 & \mbox{for $z = 0$}~. \farray \right. \label{xquaah} \feq
\end{prop}
\textbf{Proof.} According to the previous Lemma, we can apply Proposition \ref{esglobb} with
$E = K N \A^2$ and $D = \A$. Eqs. \rref{coquaa} \rref{xquaa} yield
the present relations \rref{coquaaz} \rref{xbbh} \rref{xquaah}; in particular, the
constant in  Eq. \rref{xquaa} is
$$ \R = (1 - 2 K N \A - \sqrt{1 - 4 K N \A})/(2 K N)
= K N \A^2 \, \XXX(4 K N \A)~. \hspace{0.2cm} \square $$
\begin{rema}
\textbf{Remark.} Most of the results presented in this section could be extended to
the case $\PP(f, t) = \PPP(f,....,f) + \xi(t)$, where $\PPP: \Ff^m \vain \Fm$
is~~continuous and $m$-linear for some integer $m \geqs 3$.
In this case, the growth of $\PP$ from any approximate solution admits
an estimator $\ell(r,t)$ more general than \rref{moregen}, which is polynomial of degree
$m$ in $r$.
One could  extend the analysis as well to the case $\PP(f, t) = \PPP(f,....,f, t) + \xi(t)$,
involving a time dependent multilinear map
$\PPP : \Ff^m \times [0,+\infty) \vain \Fm$, $(f_1,...,f_m,t) \mapsto \PPP(f_1,...,f_m,t)$.
These generalizations are not written only to save space.
\end{rema}
\section{The Na\-vier-Stokes (NS) equations on a torus.}
\label{inns}
From here to the end of the paper, we work in any space dimension
\beq d \geqs 2~. \feq
\textbf{Preliminaries: distributions on $\boma{\Td}$, Fourier series and
Sobolev spaces.}
Throughout this section, we use
$r,s$ as indices running from $1$ to $d$ and employ for them the Einstein summation
convention on repeated, upper and lower indices; $\delta_{r s}$ or $\delta^{r s}$ is the Kronecker symbol.
\parn
Elements $a, b$,.. of $\reali^d$ or $\complessi^d$ will be written with upper or
lower indices, according to convenience: $(a^s)$ or
$(a_s)$, $(b^s)$ or $(b_s)$. In any case, $a \sc b$ is the sum of product of the
components of $a$ and $b$, that will be written in different ways to accomplish with the Einstein
convention. Two examples, corresponding to different positions for the indices of $a$, are
\beq a \sc b = a_s b^s~, \qquad a \sc b = \delta_{r s} a^r b^s~. \feq
Let us consider the $d$-dimensional torus
\beq \Td := \To \times ... \times \To \qquad \mbox{($d$ times)}~, \qquad \To := \reali/(2 \pi \interi)~. \feq
A point of $\Td$ will be generally written $x = (x^r)_{r=1,...d}$.
We also consider the ``dual'' lattice $\Zd$ of elements $k =
(k_r)_{r=1,...,d}$ and the Fourier basis $(e_k)_{k \in \Zd}$, made of the functions
\beq e_k: \Td \vain \complessi~, \qquad e_k(x) := {1 \over (2 \pi)^{d/2}} \, e^{i k \sc \, x}~
\feq
($k \sc \, x := k_r x^r$ makes sense as an element of $\Td$).
We introduce the space of periodic distributions $\DD'(\Td, \complessi) \equiv \DD'_{\complessi}$, which is
the dual of $C^{\infty}(\Td, \complessi) \equiv C^{\infty}_{\complessi}$ (equipping the latter
with the topology of uniform convergence of all derivatives); we write
$\la v, f \ra$ for the action of a distribution $v \in \DD'_{\complessi}$
on a test function $f \in C^{\infty}_{\complessi}$. The weak topology on $\DD'_\complessi$ is the
one induced by the seminorms $p_f$ ($f \in C^{\infty}_{\complessi}$), where
$p_f(v) := | \la v, f \ra |$. \parn
Each $v \in \DD'_\complessi$
has a unique (weakly convergent) series expansion
\beq v = \sum_{k \in \Zd} v_k e_k~, \label{fs} \feq
with coefficients $v_k \in \complessi$ for all $k$, given by
\beq v_k = \la v, e_{-k} \ra~. \label{uk} \feq
The "Fourier series transformation" $v \mapsto (v_k)$ is one-to-one between $\DD'_{\complessi}$
and the space of sequences $\es'(\Zd, \complessi) \equiv \es'_\complessi$, where
\beq \es'_\complessi := \{ c = (c_k)_{k \in \Zd}~|~c_k \in \complessi,~
| c_k | = O( | k |^p) ~\mbox{as $k \vain \infty$,~ for some $p \in \reali$}\}~.
\feq
In the sequel, we often use the \textsl{mean} of a distribution $v \in \DD'_{\complessi}$, which is
\beq \la v \ra := {1 \over (2 \pi)^d} \la v, 1 \ra = {1 \over (2 \pi)^{d/2}} v_0 \label{mean} \feq
(in the first passage, $\la v, 1 \ra$ means the action of $v$ on the test function $1$;
the second relation follows from \rref{uk} with $k=0$, noting that $e_0 = 1/(2 \pi)^{d/2}$. Of course,
$\la v, 1 \ra = \int_{\Td} v(x) d x$ if $v$ is an ordinary, integrable function). \parn
The complex conjugate of a distribution $v \in \DD'_{\complessi}$ is the unique distribution $\overline{v}$ such that
$\overline{\la v, f \ra} = \la \overline{v}, \overline{f} \ra$ for each $f
\in C^{\infty}_{\complessi}$; one has $\overline{v}= \sum_{k \in \Zd} \overline{v_{k}} e_{-k}$. \parn
From now on we will be mainly interested in the space of
\textsl{real} distributions $\DD'(\Td, \reali) \equiv \DD'$, defined as follows:
\beq \DD' := \{ v \in \DD'_{\complessi}~|~\overline{v} = v \} =
\{ v \in \DD'_{\complessi}~|~\overline{v_k} = v_{-k}~\mbox{for all $k \in \Zd$} \}~; \feq
we note that $v \in \DD'$ implies $\la v \ra \in \reali$. The weak
topology on $\DD'$ is the one inherited from $\DD'_{\complessi}$. \parn
Let us write $\partial_{s}$ $(s=1,...,d)$ for the distributional derivative
with respect to the coordinate $x^s$; from these derivatives, we define the distributional
Laplacian
$\Delta := \delta^{r s} \partial_{r} \partial_s : \DD'_{\complessi} \vain \DD'_{\complessi}$. Of course
$\partial_s e_k = i k_s e_k$,  $\Delta e_k = - | k |^2 e_k$ for each $k$.
For any $v \in \DD'_{\complessi}$, this implies
\beq \partial_s v = i \sum_{k \in \Td} k_s v_k e_k~, \qquad
\Delta v = - \sum_{k \in \Td} | k |^2 v_k e_k~; \feq
\beq (1 - \Delta)^m  v = \sum_{k \in \Td} (1 + | k |^2)^m v_k e_k \label{reg} \feq
for $m \in {0,1,2,3,...}$. For any $m \in \reali$, we will regard \rref{reg} as the
definition of $(1 - \Delta)^m$ as a linear operator from $\DD'_{\complessi}$ into itself.
Comparing the previous Fourier expansions with \rref{mean}, we find
\beq \la \partial_s v \ra = 0~, \qquad \la \Delta v \ra = 0~. \label{ofc} \feq
All the above differential operators leave invariant the space of real
distributions, more interesting for us; in the sequel we will fix the attention
on the maps $\partial_s$, $\Delta$, $(1-\Delta)^m :$ $\DD' \vain \DD'$. \parn
Let us consider the real Hilbert space $L^2(\Td, \reali, d x) \equiv L^2$, i.e.,
\beq L^2 := \{ v : \Td \vain \reali~|~\int_{\Td} v^2(x) d x < + \infty \} =
\{ v \in \DD'~|~\sum_{k \in \Zd} | v_k |^2 < + \infty \}~; \feq
this has the inner product and the associated norm
\beq \la v | w \ra_{L^2} := \int_{\Td} v(x) w(x) d x = \sum_{k \in \Zd} \overline{v_k} w_k~, \label{inner} \feq
\beq \| v \|_{L^2} := \sqrt{ \int_{\Td} v^2(x) d x} = \sqrt{\sum_{k \in \Zd} | v_k |^2}~. \label{norms} \feq
To go on, we introduce the Sobolev spaces $H^n(\Td, \reali) \equiv H^n$. For each $n \in \reali$,
\beq H^n := \{ v \in \DD'~|~(1 - \Delta)^{n/2} v \in L^2 \} =
\{ v \in \DD'~|~\sum_{k \in \Td} (1 + | k |^2)^n | v_k |^2 < + \infty~\}~;
\label{hn} \feq
this is also a real Hilbert space with the inner product
\beq \la v | w \ra_{n} := \la (1 - \Delta)^{n/2} v \, | \, (1 - \Delta)^{n/2} w \ra_{L^2}
= \sum_{k \in \Td} (1 + | k |^2)^n \, \overline{v_k}  w_k~\label{refur} \feq
and the corresponding norm
\beq \| v \|_{n} := \| (1 - \Delta)^{n/2}~ v \|_{L^2} = \sqrt{\sum_{k \in \Td} (1 + | k |^2)^n | v_k
|^2}~. \label{repfur} \feq
One proves that
\beq n \geqs n' ~~~\Rightarrow~~~ H^{n} \hookrightarrow H^{n'}~, ~~ \|~\|_{n'} \leqs \|~\|_{n}~. \label{imb} \feq
In particular, $H^0$ is the space $L^2$ and contains $H^n$ for each $n \geqs 0$; for
any real $n$, $\Delta$ is a continuous map of $H^n$ into $H^{n-2}$. Finally, let us
recall that
\beq H^n \hookrightarrow \DD' \qquad \mbox{for each $n \in \reali$}~; \label{emb0} \feq
\beq H^{n} \hookrightarrow C^q \qquad \mbox{if $q \in \naturali$, ~ $n \in (q + d/2, + \infty)$}~. \label{emb} \feq
In the above $H^n$ carries its Hilbertian topology, and $\DD'$ the weak topology;
$C^q$ stands for the space $C^q(\Td,\reali)$, with the topology of uniform convergence
of all derivatives up to order $q$. \parn
Obviously enough, we could define as well the complex Hilbert spaces $L^2_{\complessi}$
and $H^n_{\complessi}$; however, these are never needed in the sequel.
\salto
\textbf{Spaces of vector valued functions on $\Td$.}
To deal with the NS equations, we need vector extensions of all the above spaces and
mappings. Let us stipulate the following: if $\Vi(\Td, \reali) \equiv \Vi$ is any vector space of \textsl{real}
functions or distributions on $\Td$, then
\beq \bVi := \Vi^d = \{ v = (v^1,...,v^d)~|~v^r \in \Vi \quad \mbox{for all $r$}\}~. \feq
This notation allows to define the spaces $\bb{\DD}'$, $\bb{L}^2$,
$\HO{n}$.
Any $v = (v^r) \in \bb{\DD}'$ will be referred to
as a \textsl{vector field} on $\Td$. We note that $v$ has a unique Fourier series expansion \rref{fs} with coefficients
\beq v_k = (v^r_k)_{r =1,...,d} \in \complessi^d~, \qquad v^r_k  := \la v^r, e_{-k} \ra~; \feq
again, the reality of $v$ ensures $\overline{v_k} = v_{-k}$.
We define componentwisely the mean $\la v \ra \in \reali^d$
of any $v \in \bb{\DD}'$ (see Eq. \rref{mean}), the derivative operators
$\partial_s : \bb{\DD}' \vain \bb{\DD}'$, their iterates and, consequently,
the Laplacian $\Delta$. The prescription \rref{reg} gives a map $(1-\Delta)^m :
\bb{\DD}' \vain \bb{\DD}'$ for all real $m$.
Whenever $\Vi$ is made of ordinary functions, a $d$-uple $v \in \bVi$ can be identified
with a function $v : \Td \vain \reali^d$, $x \mapsto v(x) = (v^r(x))_{r=1,...,d}$. \parn
$\bb{L}^2$ is a real Hilbert space. Its
inner product is as in \rref{inner}, with $v(x) w(x)$ and $\overline{v_k} w_k$
replaced by
\beq v(x) \sc \, w(x) = \delta_{r s} v^r(x) \, w^s(x)~,
\qquad \overline{v_k} \sc \, w_k = \delta_{r s} \overline{v^r_{k}} w^{s}_{k}~; \feq
the corresponding norm is as in \rref{norms}, replacing $v^2(x)$ with
$| v(x) |^2 = \sum_{r=1}^d v^r(x)^2$ and intending $| v_k |^2 = \sum_{r=1}^d | v^r_{k} |^2$. \parn
For any real $n$, the Sobolev space $\HO{n}$ is made
of all $d$-uples $v$ with components
$v^r \in H^n$; an equivalent definition can be given via Eq.\rref{hn},
replacing therein $L^2$ with $\bb{L}^2$.
$\HO{n}$ is a real Hilbert space with the inner product
\beq \la v | w \ra_{n} := \la (1 - \Delta)^{n/2} v \, | (1 - \Delta)^{n/2} w \ra_{L^2}
= \sum_{k \in \Td} (1 + | k |^2)^n \, \overline{v_k}  \sc \, w_k~.\label{refurd} \feq
The corresponding norm $\|~\|_{n}$ is given, verbatim, by Eq. \rref{repfur}; Eq.
\rref{imb} holds as well for $\HO{n}$, $\HO{n'}$ and their norms. Let us consider
the Laplacian operator
$\Delta : \bb{\DD} \vain \bb{\DD}\,$; for any real $n$
\beq \Delta \HO{n} \subset \HO{n-2}~, \feq
and $\Delta$ is continuous with respect to the norms $\|~\|_n$, $\|~\|_{n-2}$.
The embeddings \rref{imb} \rref{emb0} \rref{emb} have obvious vector analogues. \parn
\salto
\textbf{Zero mean vector fields.} The space of these vector fields is
\beq \Dz := \{ v \in \bb{\DD}'~| ~\la v \ra = 0 \} \feq
(of course, $\la v \ra = 0$ is equivalent to the vanishing of the Fourier coefficient $v_0$).
\salto
\textbf{Divergence free vector fields.} Let us consider the divergence operator (linear, weakly
continuous)
\beq \dive : \bb{\DD}' \vain \DD'~, \qquad v \mapsto \dive \,v := \partial_r v^r~; \feq
we put
\beq \Ds := \{ v \in \bb{\DD}'~|~\dive \,v = 0 \} \feq
and refer to this as to the space of divergence free (or solenoidal) vector fields.
The description of these objects in terms of Fourier transform is obvious, namely:
\beq \dive \, v = i \, \sum_{k \in \Zd}  (k \sc \, v_k) e_k \qquad \mbox{for all}~ v = \sum_{k \in \Zd} v_k e_k \in
\bb{\DD}'~; \feq
{\vbox{
\beq \Ds = \{ v \in \bb{\DD}'~|~k \sc \, v_k = 0~\mbox{for all}~ k \in \Zd~\}
\label{redv} \feq
$$ = \{ v \in \bb{\DD}'~|~v_k \in \prec k \succ^{\perp} ~\mbox{for all}~ k \in \Zd~\}~, $$ }}
where $\prec k \succ$ is the subspace of $\complessi^d$ spanned by $k$, and ${~}^\perp$ is the
orthogonal complement with respect to the inner product
$(b, c) \in \complessi^d \times \complessi^d
\mapsto \overline{b} \sc \, c$. \parn
In the sequel, we will consider as well the subspace
\beq \bb{D}'_{\so} := \Ds \cap \Dz~. \feq
\salto
\textbf{Gradient vector fields.} Let us consider the gradient operator
(again linear, weakly continuous)
\beq \partial : \DD' \vain \bb{\DD}'~;
\qquad p \mapsto \partial p := (\partial_s p)_{s=1,...,d}~; \feq
if $p = \sum_{k \in \Zd} p_k e_k$, then
\beq \partial p = i \, \sum_{k \in \Zd} k p_k e_k~. \label{gradp} \feq
The image
\beq \Dg = \{ \partial p~|~p \in \DD' \} \feq
is a linear subspace of $\bb{\DD}'$, hereafter referred to as the space of gradient
vector fields; for any vector field $w$, comparison with \rref{gradp} gives
\beq \Dg = \{ w \in \bb{\DD}'~|~w_k \in \prec k \succ \quad \mbox{for all $k \in \Zd$} \}~.
\label{regra} \feq
Of course, if $w$ is in this subspace, the distribution $p$ such that $w = \partial p$ is defined up to
an additive constant.
We put
\beq \partial^{-1} w := \mbox{unique $p \in \DD'$ such that $w = \partial p$ and $p_0 = 0$}~. \label{put} \feq
This gives a linear map
\beq \partial^{-1} : \Dg \vain \DD'~. \feq
\salto
\textbf{The Leray projection.} Using the Fourier representations \rref{redv} \rref{regra}, one
easily proves the following facts. \parn
(i) One has
\beq \bb{\DD}' = \Ds \oplus \Dg \label{dec} \feq
in algebraic sense, i.e., any $v \in \bb{\DD}'$ has a unique decomposition as the
sum of a divergence free and a gradient vector field. \parn
(ii) The projection
\beq \LP : \bb{\DD}' \vain \Ds~, \qquad v \mapsto \LP v \feq
corresponding to the decomposition \rref{dec} is given by \parn
{\vbox{
\beq \LP v = \sum_{k \in \Zd} (\LP_k v_k) e_k \qquad \mbox{for all}~ v = \sum_{k \in \Zd} v_k e_k \in
\bb{\DD}~, \feq
$$ \LP_k := \mbox{orthogonal projection of $\complessi^d$ onto $\prec k \succ^{\perp}$}~; $$ }}
more explicitly, for all $c \in \complessi^d$,
\beq \LP_0 c = c~, \qquad \LP_k c = c - {(k \sc \, c) k \over | k |^2} \quad \mbox{for $k \in \Zd, k \neq 0$}~. \feq
As usually, we refer to $\LP$ as to the \textsl{Leray projection}; this operator is
weakly continuous. From the Fourier representations of $\LP$, of the mean and of the
derivatives one easily infers, for all
$v \in \bb{\DD}'$,
\beq \la \LP v \ra = \la v \ra~,\qquad
\LP (\partial_s v) = \partial_s (\LP v)~, \qquad \LP (\Delta v) = \Delta (\LP v)~.\label{meanlp} \feq
\salto
\textbf{A Sobolev framework for the previous decomposition.} For $n \in \reali$, let us define
\beq \Hs{n} := \HO{n} \cap \Ds = \{ v \in \HO{n}~|~\dive \,v = 0 \}~; \feq
\beq \Hg{n} := \HO{n} \cap \Dg = \{ w \in \HO{n}~|~w = \partial p,~ p \in \DD'\}~; \feq
\beq \partial H^n := \{ \partial p~|~p \in H^n~\}~. \feq
Then the following holds for each $n$: \parn
(i) $\Hg{n}$ is a closed subspace of the Hilbert space $(\HO{n}, \la~|~\ra_n)$
(because $\mbox{div}$ is continuous between this Hilbert space and $\DD'$ with the weak topology). \parn
(ii) One has
\beq \Hg{n} = \partial H^{n+1} \feq
and $\Hg{n}$ is also a closed subspace of $\HO{n}$.
The map $\partial^{-1}$ of Eq. \rref{put} is continuous between $\Hg{n}$ and
$H^{n+1}$. \parn
(iii) Denoting with ${~}^{\perp_n}$ the orthogonal complement in $(\HO{n}, \la~|~\ra_n)$, we have
\beq {\Hs{n}}^{\perp_n} = \Hg{n} \feq
and $\LP \restriction \HO{n}$ is the orthogonal projection of $\HO{n}$ onto
$\Hs{n}$; so, as usual for Hilbertian orthogonal projections,
\beq \| \LP v \|_n \leqs \| v \|_n \qquad \mbox{for all $v \in \HO{n}$}~. \label{hilb} \feq
\salto
\textbf{Other Sobolev spaces of vector fields.} For $n \in \reali$, we put
\beq \Hz{n} := \HO{n} \cap \Dz:= \{ f \in \HO{n}~|~\la f \ra = 0 \}~; \feq
\beq \HM{n} := \Hs{n} \cap \Hz{n} := \{ f \in \HO{n}~|~\dive f = 0~, \la f \ra = 0 \}
\label{defhm} \feq
Then, $\Hz{n}$ is a closed subspace of $\HO{n}$ (by the continuyity
of $\la ~ \ra : \HO{n} \vain \complessi$); the same holds for $\HM{n}$,
since this is the intersection of two closed subspaces. \parn
The space \rref{defhm} plays an important role in the sequel; we will often use
the Fourier representation
\beq \HM{n} = \{ f \in \bb{\DD}~|~\sum_{k \in \Zd}
(1 + | k |^2)^n | f_k |^2 < + \infty~, \feq
$$ k \sc \, f_k = 0~\mbox{for all $k \in \Zd$}, f_0 = 0 \}~. $$
\salto
\textbf{Some inclusions.}
We note that the relations $\Delta \HO{n} \subset \HO{n-2}$ and
$\dive (\Delta v) = \Delta (\dive v)$, $\la \Delta v \ra = 0$  for all $v \in \bb{D}$ imply
the following, for each real $n$:
\beq \Delta \HO{n} \subset \Hz{n-2}~, \qquad \Delta \Hs{n} \subset \HM{n-2}~. \feq
\salto
\textbf{A digression: estimates on certain series.}
Let us define
\beq \Zd_0 := \Zd \setminus \{ 0 \}~; \feq
throughout the section, $n$ is a real number such that
\beq n > {d \over 2}~. \label{nd} \feq
The series considered hereafter are used shortly
afterwards to derive quantitative estimates on the fundamental bilinear map appearing in
the NS equations: by this we mean the map sending two vector fields $v, w$ on $\Td$
into the vector field $v \sc \, \partial w$ (see the next paragraph). The estimates we give
are also useful for the numerical computation of those series. In both Lemmas hereafter,
\beq \Ze := \Zd ~\mbox{or}~\Zd_0~. \label{ze} \feq
\begin{prop}
\label{lemsi}
\textbf{Lemma.} One has
\beq \Sigma_n := {1 \over (2 \pi)^d} \sum_{h \in \Ze} {1 \over (1 + |h|^2)^n} < + \infty~. \label{desig} \feq
For any real "cutoff" $\lambda \geqs 2 \sqrt{d}$, one has
\beq \Sigm_n(\lambda) < \Sigma_n \leqs \Sigm_n(\lambda) + \delta \Sigm_n(\lambda) \label{decosig} \feq
where
\beq \Sigm_n(\lambda) := {1 \over (2 \pi)^d}  \sum_{h \in \Ze, |h| < \lambda} {1 \over (1 + |h|^2)^n}~, \label{desgm} \feq
\beq \delta \Sigm_n(\lambda) := {(1 + d)^n \over 2^{d-1} \pi^{d/2} \Gamma(d/2) (2 n - d)~
(\lambda - \sqrt{d})^{2 n - d}}~. \label{reph} \feq
\end{prop}
\textbf{Proof.} See Appendix \ref{appemu}.
\fine
\begin{prop}
\label{lemci}
\textbf{Lemma.} For $k \in \Ze$, define
\beq \KK_n (k) \equiv \KK_{n d}(k) := {(1 + | k |^2)^{n-1} \over (2 \pi)^d}
\sum_{h \in \Ze} {| k - h |^2 \over (1 + |h|^2)^n (1 + | k - h |^2)^n}~;  \label{mndi2} \feq
then, (i) (ii) hold. \parn
(i) One has $\KK_n(k) < + \infty$ for all $k \in \Ze$;
furthermore, with $\Sigma_n$ as in \rref{desig},
\beq \KK_n(k) \vain \Sigma_n \qquad \mbox{for $k \vain \infty$}~. \label{likn} \feq
Thus,
\beq \sup_{k \in \Ze} \KK_n(k) < + \infty~. \feq
(ii) Let us choose any "cutoff function"
$\Lambda : \Ze \vain [2 \sqrt{d}, + \infty)$ and define
\beq \lambda : \Ze \vain (0,+\infty)~,
\qquad k \mapsto \lan(k) := \left\{ \barray{ll} 1 + | k |^2 & \mbox{if $\Lambda(k) < | k |$}~, \\
\dd{{1 + | k |^2 \over 1 + (\Lambda(k) - | k |)^2}} & \mbox{if $\Lambda(k) \geqs | k |$}~. \farray \right.
\label{defom} \feq
Then, for all $k \in \Ze$,
\beq \Ki_n(k) < \KK_n(k) \leqs \Ki_n(k) + \delta \Ki_n(k)~, \label{kkn} \feq
where
\beq \Ki_n(k) := {(1 + | k |^2)^{n-1} \over (2 \pi)^d}
\sum_{h \in \Ze, |h | < \Lambda(k)}
{| k - h |^2 \over (1 + |h|^2)^n (1 + | k - h |^2)^n}~, \label{mndig2} \feq\beq \delta \Ki_n(k) := {(1 + d)^n \lan(k)^{n-1} \over 2^{d-1} \pi^{d/2} \Gamma(d/2) (2 n - d)~
(\Lambda(k) - \sqrt{d})^{2 n - d}}~. \label{deci} \feq
Finally, suppose the cutoff $\Lambda$ has the property
\beq \alpha | k | \leqs \Lambda(k) \leqs
\beta | k | \qquad \mbox{for all $k \in \Ze$ such that $| k | \geqs \chi$} \label{cutoff} \feq
$$ (\chi \geqs 0,~~1 < \alpha \leqs \beta)~; $$
then
\beq \Ki_n(k) \vain \Sigma_n~, \qquad \delta \Ki_n(k) =
O({1 \over | k |^{2 n - d}}) \vain 0 \qquad \mbox{for $k \vain \infty$}~.
\label{nona} \feq
\end{prop}
\textbf{Proof.} See Appendix \ref{appemu}. \fine
\salto
\textbf{The fundamental bilinear map.}
Let $v, w \in \HO{n}$. For each $r, s \in \{1, 2, ..., d\}$, $v^r \in H^n$ and $\partial_r w^s \in H^{n-1}$ are ordinary
real functions: note that $v^r \in C^0$ by the
embedding \rref{emb}, and $\partial_r w^s \in L^2$ since $n - 1 > d/2 - 1 \geqs 0$. These
functions can be multiplied pointwisely, and this allows to define
\beq  v \sc \, \partial w := \left( v \sc \, \partial w^s \right)_{s=1,...,d}
\qquad v \sc \, \partial w^s := v^r \partial_r w^s : \Td \vain \complessi~. \feq
\begin{prop}
\label{multi}
\textbf{Proposition.} (i) Consider $v, w \in \HO{n}$. The vector field
$v \sc \, \partial w$ has Fourier coefficients
\beq (v \sc \, \partial w)_k = {i \over (2 \pi)^{d/2}}~\sum_{h \in \Zd} [(v_{h} \sc \, (k - h)] w_{k - h}~.
\label{fuco} \feq
Furthermore,  $v \sc \, \partial w \in \HO{n-1}$. \parn
(ii) The bilinear map
\beq \HO{n} \times \HO{n} \vain \HO{n-1}~, \qquad
(v,w) \mapsto v \sc \, \partial w \feq
admits an estimate (indicating continuity)
\beq \| v \sc \, \partial w \|_{n-1} \leqs K_{n} \,  \| v \|_n \| w \|_n \label{esto} \feq
for all $v, w$ as above, with a suitable constant $K_{n} \equiv K_{n d} \in (0,+\infty)$. For the latter
one can take any constant such that
\beq \sqrt{\sup_{k \in \Zd} \KK_{n}(k)} \leqs K_n~, \label{mndi} \feq
\beq \KK_n (k) := {(1 + | k |^2)^{n-1} \over (2 \pi)^d}
\sum_{h \in \Zd} {| k - h |^2 \over (1 + |h|^2)^n (1 + | k - h |^2)^n} \label{mmnnddii} \feq
(as in \rref{mndi2}, with $\Ze = \Zd$).
\end{prop}
\textbf{Proof.} See Appendix \ref{appemult}. \fine
\begin{prop}
\label{lemean}
\textbf{Lemma.} Let $v \in \Hs{n}$, $w \in \HO{n}$. Then
\beq \la v \sc \, \partial w \ra = 0~; \label{elem} \feq
combined with Proposition \ref{multi}, this gives $v \sc \, \partial w \in \Hz{n-1}$.
\end{prop}
\textbf{Proof.} For $s=1,...,d$, integration by parts and the assumption
$0 = \mbox{div} v = \partial_ r v^r$ give
$$\hspace{2cm}
 \la v \sc \, \partial w^s \ra
= {1 \over {(2 \pi)}^d} \int_{\Td} \! \! \! \! \! d x \, v^r \partial_r w^s
= - {1 \over {(2 \pi)}^d} \int_{\Td} \! \! \! \! \! d x \, (\partial_ r v^r) w^s = 0~.
\hspace{2cm} \square $$
\begin{prop}
\label{mult}
\textbf{Proposition.} The bilinear map
\beq \HM{n} \times \HM{n} \vain \Hz{n-1}~, \qquad
(f,g) \mapsto f \sc \, \partial g \feq
admits an estimate
\beq \| f \sc \, \partial g \|_{n-1} \leqs K_{n} \,  \| f \|_n  \, \| g \|_n \label{est} \feq
for all $f, g \in \HM{n}$, with a suitable constant $K_{n} \equiv K_{n d} \in (0,+\infty)$. One can take
for the latter any constant such that
\beq \sqrt{\sup_{k \in \Zd_0} \KK_{n}(k)} \leqs K_n~, \label{mnd} \feq
\beq \KK_n (k) := {(1 + | k |^2)^{n-1} \over (2 \pi)^d}
\sum_{h \in \Zd_0} {| k - h |^2 \over (1 + |h|^2)^n (1 + | k - h |^2)^n} \label{mnd2} \feq
(as in \rref{mndi2}, with $\Ze = \Zd_0$).
\end{prop}
\textbf{Proof.} See Appendix \ref{appemult}. \fine
\salto
To conclude, we report another integral identity frequently used in the sequel.
\begin{prop}
\label{td}
\textbf{Lemma.} For any $v \in \Hs{n}$, one has
\beq \la v | v \sc \, \partial v \ra_{L^2} = 0 \label{etd} \feq
\end{prop}
\textbf{Proof.} We have (with $q,s,r \in \{1,...,d\}$),
\beq \la v | v \sc \, \partial v \ra_{L^2} = \int_{\Td} d x \, \delta_{q s} v^q (v^r \partial_r v^s)~. \label{ehand} \feq
From here we infer, integrating by parts,
$$ \la v | v \sc \, \partial v \ra_{L^2} =
- \int_{\Td} d x \, \delta_{q s} \partial_r (v^q v^r) v^s  $$
\beq = - \int_{\Td} d x \, \delta_{q s} (\partial_r v^q) v^r v^s - \int_{\Td} d x \, \delta_{q s} v^q (\partial_r v^r) v^s
= - \int_{\Td} d x \, \delta_{q s} (\partial_r v^q) v^r v^s \feq
since $\dive v = 0$; by comparison with \rref{ehand} we obtain
\beq \la v | v \sc \, \partial v \ra_{L^2} = - \la v | v \sc \, \partial v \ra_{L^2}~, \feq
whence the thesis \rref{etd}. \fine
\begin{rema}
\textbf{Remark.} The inequality \rref{esto} (or \rref{est}) is known from the literature, in this
form or in some variant: see, for example, one of the standard references on NS equations cited in the
Introduction. To our knowledge, the novelty of Proposition \ref{multi} (or \ref{mult}) with respect
to the already published material is the rule \rref{mndi} (or \rref{mnd}) to
determine  $K_n$, that can be used with Lemmas \ref{lemsi}, \ref{lemci} to provide
a numerical value for this constant. Examples of this computation appear in Section
\ref{nume} and Appendix \ref{appekdue}.
\parn
The method employed in Appendix \ref{appemult} to prove Propositions \ref{multi}, \ref{mult}
is very similar to one employed in \cite{MP} to estimate the product
of two scalar functions in $H^{n}(\reali^d)$. In that paper, already mentioned in the Introduction,
we have given a rule similar to \rref{mndi} (or \rref{mnd}) to find a constant $C_{n d} \equiv C_n$
such that $\| p  q \|_{H^{n}(\reali^d)} \leqs C_n \| p  \|_{H^n(\reali^d)} \| q \|_{H^n(\reali^d)}$;
furthermore, using convenient trial functions $p, q$ we have shown this constant to be
very close to the smallest one fulfilling the inequality. Due to
the similarities with \cite{MP}, the constant $K_n$
provided by \rref{mndi} (or \rref{mnd}) is hopingly close to the smallest one
for the inequality \rref{esto} (or \rref{est}).
\end{rema}
\salto
\textbf{Functions on $\Td$, depending on time.} Suppose we have a function
\beq \rho : [t_0, T) \subset \reali \vain \Vi~ \mbox{or}~\bVi~,
\qquad t \mapsto \rho(t) \feq
where $\Vi$, $\bVi$ stand for some spaces of $\reali$ or $\reali^d$ valued functions on $\Td$.
At each "time" $t \in [t_0, T)$, this gives a function $x \in \Td \mapsto \rho(t)(x)$; of course,
we will use the more common notation
\beq \rho(x, t) := \rho(t)(x)~. \feq
\salto
\textbf{The incompressible NS equations, in the Leray formulation.}
\salto
Let us recall that $C^{0,1}$ indicates the locally Lipschitz maps.
\begin{prop}
\label{defnsl}
\textbf{Definition.}  The \textsl{incompressible NS Cauchy problem with initial
datum} $v_0 \in \Hs{n+1}$ and forcing term
$\eta \in C^{0,1}([0,+\infty), \Hs{n-1})$, in the \textsl{Leray formulation}, is the following.
$$ \mbox{Find $\nu \in C([0,T),\Hs{n+1}) \cap C^1([0,T), \Hs{n-1})$,
such that} $$
\beq {\dot \nu}(t) = \Delta \nu(t)- \LP \big( \nu(t) \sc \, \partial \nu(t) \big) + \eta(t)
\quad \mbox{for $t \in [0,T)$}~, \qquad \nu(0) = v_0~ \label{nsl} \feq
(for some $T \in (0,+\infty]$).
\end{prop}
\begin{rema}
\textbf{Remarks. (i)} The requirement $\nu \in C([0,T),\Hs{n+1})$ ensures
by itself that the right hand side of the above differential equation
is in $C([0,T),\Hs{n-1})$. \parn
\textbf{(ii)} The differential equation in \rref{nsl} can be interpreted as
the usual NS equation for an incompressible fluid, in a convenient adimensional
formulation. For each $t \in [0,T)$, $\nu(t) : x \in \Td \mapsto \nu(x,t)$ is the
velocity field of the fluid at time $t$; $\eta(t) : x \mapsto \eta(x,t)$ is the Leray projection
of the density of external forces. Of course, $\nu(t)$ is taken in the divergence free space
$\Hs{n+1}$ to fulfill the condition of incompressibility; the pressure gradient does not
appear in \rref{nsl}, having been eliminated by
application of $\LP$. For completeness, all these facts are surveyed in Appendix \ref{appens}.
\ffine
\end{rema}
The statement that follows refers to the time evolution of the functions
$t \mapsto \la \nu(t) \ra = (2 \pi)^{-d} \int_{\Td} \nu(x, t) d x$ and
$t \mapsto$ $(1/2) \| {\nu}(t) \|^2_{L^2}$. Up to dimensional factors, the first one
gives the mean value of the velocity field or, equivalently, the total momentum;
the second one gives the total kinetic energy.
\begin{prop}
\label{meanco}
\textbf{Proposition.} Suppose $\nu$ fulfills \rref{nsl} on an interval $[0,T)$.
Then, for all $t$ in this interval, we have the following relations: \parn
(i) Balance of momentum:
\beq {d \over d t} \la \nu(t) \ra = \la \eta(t) \ra~. \label{consm} \feq
(ii) Balance of energy:
\beq {1 \over 2} {d \over d t} \| {\nu}(t) \|^2_{L^2}
=  \la \nu(t) | \Delta \nu(t) \ra_{L^2} + \la \nu(t) | \eta(t) \ra_{L^2}
\leqs \la \nu(t) | \eta(t) \ra_{L^2}~. \label{conse} \feq
\end{prop}
\textbf{Proof.} (i) Let us observe that
\beq {d \over d t} \la \nu \ra = \la \dot \nu \ra = \la \Delta \nu \ra -
\la \LP(\nu \sc \, \partial\nu) \ra + \la \eta \ra~; \feq the means of
$\Delta \nu$ and $\LP(\nu \sc \, \partial\nu)$ vanish due
to Eqs. \rref{ofc} \rref{meanlp} \rref{elem}, so we get the thesis \rref{consm}. \parn
(ii) Let us write
\beq {1 \over 2} {d \over d t} \| {\nu} \|^2_{L^2} = \la \nu | \dot \nu \ra_{L^2} =
\la \nu | \Delta \nu \ra_{L^2} - \la \nu | \LP(\nu \sc \, \partial \nu) \ra_{L^2} + \la \nu | \eta \ra_{L^2}~;
\label{eu} \feq
on the other hand by the symmetry of $\LP$, the equality $\LP \nu = \nu$ and Lemma \ref{td},
\beq \la \nu | \LP(\nu \sc \, \partial \nu) \ra_{L^2} = \la \LP \nu | \nu \sc \, \partial \nu \ra_{L^2} =
\la \nu | \nu \sc \, \partial \nu \ra_{L^2} = 0~; \label{ed} \feq
these facts yield the equality in \rref{conse}. The subsequent inequality in \rref{conse}
follows via elementary integration by parts:
\beq \la \nu | \Delta \nu \ra_{L^2} = \delta_{r s} \int_{\Td} v^r \Delta v^s =
- \delta_{r s} \int_{\Td} \partial \nu^r \sc \, \partial \nu^s \leqs 0~. \qquad
\square \label{et} \feq
\salto
\textbf{Reducing the NS equations to the case of a zero mean velocity field.}
Let us recall the notation
$\HM{n}$ for the space of divergence free, zero mean vector fields (see Eq.
\rref{defhm} and subsequent comments); we regard this as a Hilbert space, with the inner
product $\la~|~\ra_n$ and the norm $\|~\|_n$ inherited from $\HO{n}$. \parn
The purpose of this paragraph is to show that the general Cauchy problem \rref{nsl}
can be reduced to a Cauchy problem for zero mean vector fields; let us define the latter
precisely.
\begin{prop}
\textbf{Definition.} Let $f_0 \in \HM{n+1}$ and $\xi \in C^{0,1}([t_0, +\infty), \Hs{n-1})$.
The \textsl{incompressible, zero mean NS Cauchy problem with initial
datum $v_0$ and forcing term $\xi$} is the following.
$$ \mbox{Find $\varphi \in C([0,T),\HM{n+1}) \cap C^1([0,T), \HM{n-1})$
such that} $$
\beq {\dot \varphi}(t) = \Delta \varphi(t)- \LP \big( \varphi(t) \sc \, \partial \varphi(t) \big) + \xi(t)
\quad \mbox{for $t \in [0,T)$}~, \qquad \varphi(0) = f_0~ \label{nsm} \feq
(for some $T \in (0,+\infty]$).
\end{prop}
Let us connect this problem with the previous one \rref{nsl}, for given $v_0$ and
$\eta$. To this purpose, we need a bit more regularity on the forcing $\eta$; to be
precise, we assume
\beq v_0 \in \Hs{n+1}~, \qquad \eta \in C^{0,1}([0,+\infty), \Hs{n-1}) \cap C([0,+\infty),
\Hs{n})~. \label{assumeta} \feq
Let us define from $v_0$ and $\eta$ the following objects:
\beq m_0 := \la v_0 \ra \in \reali^d~; \quad
m \in C([0,+\infty), \reali^d)~,~t \mapsto m(t) := m_0 + \int_{0}^t d s \, \la \eta(s) \ra~;
\label{defem} \feq
\beq h \in C^1([0,+\infty), \reali^d)~, \qquad t \mapsto h(t) := \int_{0}^t d s \, m(s)~; \feq
\beq f_0 := v_0 - \m_0 \in \HM{n+1}~; \label{efz} \feq
\beq \xi : [0,T) \vain \HM{n}~, \qquad t \mapsto \xi(t)~
\mbox{such that}~\xi(x, t) := \eta(x + h(t), t) - \la \eta(t) \ra~. \label{defxi} \feq
The statement $\xi(t) \in \HM{n}$ for each $t$ is evident from the definition. In
Appendix \ref{appexi} we will prove that
\beq \xi \in C^{0,1}([0,+\infty), \HM{n-1})~. \label{xilip} \feq
\begin{prop}
\textbf{Proposition.} Let us consider: \parn
(i) the Cauchy prolem \rref{nsl}, with any datum $v_0$ and forcing $\eta$ as
in \rref{assumeta}; \parn
(ii) the above definitions of $m_0, m, h, f_0, \xi$ and the Cauchy prolem \rref{nsm}. \parn
A function $\nu$ of domain $[0,T)$ fulfills \rref{nsl} if and only if there is a function
$\varphi$ on $[0, T)$ fulfilling \rref{nsm} such that, for all $x \in \Td$ and $t \in [0,T)$,
\beq \nu(x, t) = m(t) + \varphi(x - h(t), t)~. \label{nufi} \feq
\end{prop}
\textbf{Proof.} \textsl{Step 1.} We suppose \rref{nsm} with the above datum $f_0$
to have a solution $\varphi$ on $[0, T)$; we define $\nu$ as in \rref{nufi} and
prove that it solves problem \rref{nsl}. It is clear that $\nu$ is in the functional
space prescribed by \rref{nsl}, and that the following holds:
\beq \nu(x,0) = m_0 +\varphi(x,0) = m_0 + f_0(x) = v_0(x)~; \label{eq1} \feq
\beq {\dot \nu}(x,t) = \la \eta(t)\ra + {\dot \varphi}(x - h(t),t) - {\dot h}(t) \sc \, \partial \varphi(x - h(t),t)
\label{eq2} \feq
$$ = \la \eta(t)\ra + {\dot \varphi}(x - h(t),t) - m(t) \sc \, \partial \varphi(x - h(t)t,t)~; $$
\beq  \Delta \nu(x,t) = \Delta \varphi(x - h(t),t)~; \label{eq5} \feq
\beq (\nu \sc \, \partial \nu)(x,t) = (m(t) \sc \, \partial \varphi)(x - h(t),t) +
(\varphi \sc \, \partial  \varphi)(x - h(t),t)~. \feq
The Leray projection $\LP$ commutes with space translations, so the last equation
implies
\beq \LP(\nu \sc \, \partial \nu)(x,t) = \LP(m(t) \sc \, \partial  \varphi)(x - h(t),t) +
\LP(\varphi \sc \, \partial \varphi)(x - h(t),t)~. \label{eq3} \feq
To conclude, we note that
\beq \LP(m(t) \sc \, \partial  \varphi) = m(t) \sc \, \partial \varphi~; \label{eq4} \feq
to prove this, it suffices to check that $m(t) \sc \, \partial \varphi$ is divergence free. In fact,
\beq \mbox{div}(m(t) \sc \, \partial \varphi) = \partial_s (m^r(t) \partial_r \varphi^s) =
m^r(t) \partial_r (\partial_s \varphi^s) = m^r(t) \partial_r (\dive \varphi) = 0~,
\feq since $\varphi$ is divergence free at all times. From Eq. \rref{eq1} we see
that $\nu$ fulfills the initial condition in \rref{nsl}; from Eqs. \rref{eq2}
\rref{eq5} \rref{eq3} \rref{eq4} and \rref{nsm}, we see that $\nu$ fulfills the evolution equation in
\rref{nsl}. \parn \textsl{Step 2.} Let us consider function $\nu$ on $[0,T)$,
fulfilling \rref{nsl}; \ we will prove the existence of a
function $\varphi$ on $[0,T)$ fulfilling \rref{nsm},
such that $\nu$ and $\varphi$ are related by \rref{nufi}. To this purpose, let us
define a function $\varphi$ by
\beq \varphi(x, t) := \nu (x + h(t), t) - m(t)~; \feq
then, at each time $t \in [0, T)$,
\beq \la \varphi(t) \ra =  \la \nu(t) \ra -  m(t)  = 0 \feq
on account of Eq. \rref{consm} for $\la \nu \ra$ and of the definition \rref{defem} of $m$.
Besides having zero mean, $\varphi$ belongs to the function
spaces prescribed by \rref{nsl} due to the properties of $\nu$. Now, computations very similar to
the ones of Step 1 prove that $\varphi$ fulfills the initial condition and the evolution equation in
\rref{nsm}. \fine
\begin{rema}
\textbf{Remark.} Let us regard Eqs. \rref{efz} \rref{defxi} as defining a transformation
\beq \TT : \Hs{n+1} \times C^{0,1}([0,+\infty), \Hs{n-1}) \vain \HM{n+1}
\times C^{0,1}([0,+\infty), \HM{n-1})~, \feq
$$ (v_0, \eta) \mapsto (f_0, \xi) = \TT(v_0, \eta)~. $$
The map $\TT$ is onto, due to the trivial equality $(f_0, \xi) = \TT(f_0, \xi)$.
\end{rema}
\section{The NS equations in the general framework for evolution equations with quadratic
nonlinearity.}
\label{basic}
\textbf{Basic notations.}
The zero mean version \rref{nsm} is the final form for the NS Cauchy problem, to which we stick from now
on. Let us recall that $d \in \{2,3,...\}$; throughout the section, we fix
\beq n > {d \over 2}~, \qquad \xi \in C^{0,1}([0,+\infty), \HM{n-1}) \label{xiasin} \feq
(the function $\xi$ is regarded to be given by itself, independently of any function $\eta$ as in
the previous section).
Our aim is to apply the formalism of Section \ref{quadr} to the Cauchy problem \rref{nsm} (and
to the equivalent Volterra problem); in this case
\beq \Ff_{\pm} \equiv (\Ff_{\pm}, \|~\|_{\pm}) := (\HM{n \pm 1}, \|~\|_{n \pm 1})~,~~
\Ff \equiv (\Ff, \|~\|) := (\HM{n}, \|~\|_n)~, \feq
\beq \AA := \Delta : \HM{n+1} \vain \HM{n-1}~, \qquad f \mapsto \Delta f~; \feq
\beq \PPP : \HM{n} \times \HM{n} \vain \HM{n-1}~,~~~\PPP(f, g) := - \LP(f \sc \, \partial g)~;
\qquad \mbox{$\xi$ as in \rref{xiasin}}~. \label{depp} \feq
From $\PPP$ and $\xi$, we construct the function
\beq \PP : \HM{n} \times [0,+\infty) \vain \HM{n-1}~, \qquad \PP(f,t) :=
\PPP(f,f) + \xi(t) = - \LP(f \sc \, \partial f) + \xi(t)~, \label{dep} \feq
which appears in \rref{nsm}. \parn
For subsequent reference, we record the Fourier representations
$$ \| f \|_{n \pm 1} = \sqrt{\sum_{k \in \Zd_0} (1 + | k |^2)^{n \pm 1} | f_k |^2}, \quad
~\| f \|_n = \sqrt{\sum_{k \in \Zd_0} (1 + | k |^2)^{n} | f_k |^2},~$$
\beq (\Delta f)_k = - | k |^2 f_k~. \label{recall} \feq
\textbf{Verification of properties (\PU)-(\PC').} The above set $(\Fm, \Ff, \Fp, \AA)$, completed
with suitable semigroup estimators $u$, $\um$, has the properties (\PU)-(\PC) prescribed in
Section \ref{prelim},
and (\PQ')(\PC') of Section \ref{quadr}.
We will indicate which parts of the proof are obvious, and give details on the nontrivial parts.
\begin{prop}
\textbf{Proposition.} $\Ff = \HM{n}$, $\Ff_{\mp} = \HM{n \mp 1}$ and $\AA = \Delta$ fulfill conditions (\PU)-(\PD).
\end{prop}
\textbf{Proof.} Everything follows easily from the Fourier representations. \fine
\begin{prop}
\label{edelta}
\textbf{Proposition.} (i) $\Delta$ generates a strongly continuous semigroup on $\HM{n-1}$, given by
\beq e^{t \Delta} f = \sum_{k \in \Zd_0} e^{-t | k |^2} f_k e_k~, \qquad
\mbox{for $f \in \HM{n-1}$, $t \in [0,+\infty)$}~. \feq
So, (\PT) holds. \parn
(ii) For $f \in \HM{n}$ and $t \in [0,+\infty)$ one has
\beq e^{t \Delta} f \in \HM{n}~, \qquad \no{e^{t \Delta} f}_n \leqs u(t) \no{f}_n~,~~
u(t) := e^{-t}~; \label{reguff} \feq
the function $(f, t) \mapsto e^{t \Delta} f$ gives a strongly continuous semigroup on $\HM{n}$. \parn
(iii) For $f \in \HM{n-1}$ and $t \in (0,+\infty)$ one has
\beq e^{t \Delta} f \in \HM{n}~, \qquad \no{e^{t \Delta} f}_n \leqs \um(t) \no{f}_{n-1} ~, \label{eqvv} \feq
\beq \um(t) := e^{-t} \uv(t)~, \qquad
\uv(t) := \left\{ \barray{ll} \dd{e^{2 t} \over \sqrt{2 e t}} & \mbox{for $0 < t \leqs \dd{1 \over 4}$}~,
\\ \sqrt{2} & \mbox{for $t > \dd{1 \over 4}$}~; \farray \right. \label{eqww} \feq
note that $\um(t), \uv(t) = O(1/\sqrt{t})$ for $t \mapsto 0^+$.
The function $(f, t) \mapsto e^{t \Delta} f$
is continuous from $\HM{n-1} \times (0,+\infty)$ to $\HM{n}$. \parn
With $\um$ as in \rref{eqww}, the function $t \mapsto \MM(t) := \int_{0}^t d s \, \um(s)$ is given by
\beq \MM(t) := \left\{ \barray{ll} \dd{\gamma(t) \over \sqrt{2}} & \mbox{for $0 < t \leqs \dd{1 \over 4}$}~,
\\ \dd{ {\gamma(1/4) \over \sqrt{2}} + \sqrt{2} (e^{-1/4} - e^{-t})}
& \mbox{for $\dd{1 \over 4} < t \leqs +\infty$}~, \farray \right. \label{eqmm} \feq
\beq \gamma(t) := \int_{0}^t d s \, {e^{s} \over \sqrt{s}} \qquad \mbox{for $0 \leqs t \leqs \dd{1 \over 4}$}~. \feq
(In particular, $\MM(+\infty) = \gamma(1/4)/\sqrt{2} + \sqrt{2} \, e^{-1/4} \in (1.872, 1.873)$). \parn
(iv) With $\uv$ as in \rref{eqww}, one has
\beq \sup_{t \in [0,+\infty)}
\int_{0}^t d s \, \uv(t-s) e^{-s} = \sqrt{2}~. \label{defn} \feq
(v) In conclusion, (\PQ)(\PQ') and (\PC)(\PC') hold with
\beq B = 1~, \qquad \sigma = {1 \over 2}~, \qquad N = \sqrt{2}~. \label{holdw} \feq
\end{prop}
\textbf{Proof.} (i) This follows basically from Eq. \rref{recall} for $\Delta$. \parn
\parn (ii)(iii) We
only give details on the derivation of Eqs. \rref{reguff}-\rref{eqmm}.
Let $f \in \HM{n}$, $t \in [0,+\infty)$. Then, \parn
{\vbox{
\beq \sum_{k \in \Zd_0} (1 + | k |^2)^n |(e^{t \Delta} f)_k|^2 =
\sum_{k \in \Zd_0} (1 + | k |^2)^n e^{- 2 t | k |^2} |f_k|^2
\feq
$$
\leqs e^{-2 t} \sum_{k \in \Zd_0} (1 + | k |^2)^n |f_k|^2~,
$$}}
since $| k | \geqs 1$ for $k \in \Zd_0$;
this yields Eq. \rref{reguff}. Now, let $f \in \HM{n-1}$ and $t \in (0,+\infty)$; then,
\beq \sum_{k
\in \Zd_0} (1 + | k |^2)^n |(e^{t \Delta} f)_k|^2
= \sum_{k \in \Zd_0} (1 + | k |^2) e^{- 2 t | k |^2} (1 + | k |^2)^{n-1} |f_k|^2 \label{thenn} \feq
$$ \leqs \Big(\sup_{\te \in [1,+\infty)} U_t(\te) \Big) \Big(\sum_{k \in \Zd_0} (1 + | k
|^2)^{n-1} |f_k|^2 \Big)~, $$
$$ U_t(\te):= (1 + \te) e^{-2 t \te}~, $$
and an elementary computation gives
\beq \sup_{\te \in [1,+\infty)} U_t(\te) =
\left\{ \barray{ll} U_t\Big(\dd{1 \over 2 t} - 1\Big) = \dd{e^{2 t} \over 2 e t} & \mbox{for $0 < t \leqs \dd{1 \over 4}$}~,
\\ U_t(1) = 2 e^{-2 t} & \mbox{for $t > \dd{1 \over 4}$}~. \farray \right. \label{thee} \feq
From \rref{thee} we infer Eq. \rref{eqvv} with $\um(t) := \sqrt{\sup_{\te \in
[1,+\infty)} U_t(\te)}$, i.e.,
\beq \um(t) =
\left\{ \barray{ll} \dd{e^{t} \over \sqrt{2 e t}} & \mbox{for $0 < t \leqs \dd{1 \over 4}$}~,
\\ \sqrt{2} e^{- t} & \mbox{for $t > \dd{1 \over 4}$}~; \farray \right. \feq
this definition of $\um(t)$ agrees with Eq. \rref{eqww}, and \rref{eqmm} follows trivially.
\parn
(iv) See Appendix \ref{aplemen}. \parn
(v) Obvious consequence of items (i-iv). \fine
\salto
\textbf{Analysis of $\boma{\PPP}$ and $\boma{\PP}$.} We turn the attention to the functions in Eqs. \rref{depp} \rref{dep}.
\begin{prop}
\label{isbil}
\textbf{Proposition.}
(i) $\PPN$ is a bilinear map and admits the estimate (implying continuity)
\beq \| \PPN(f,g) \|_{n-1} \leqs K_n \| f \|_n \| g \|_n \feq
for all $f, g \in \HM{n}$, with $K_n \equiv K_{n d}$ any constant fulfilling \rref{est} (so,
condition (\Q)~holds for this map). \parn
(ii) As a consequence of (i), $\PN$ fulfills the Lipschitz condition (\PS). Furthermore,
for each function $\fiapp \in C([0,T), \HM{n})$, the growth of $\PN$ from $\fiapp$
admits the estimate
\beq \| \PN(f, t) - \PN(\fiapp(t), t) \|_{n-1} \leqs \ell_n(t, \| f - \fiapp(t) \|_n) \feq
for $t \in [0,T)$ and $f \in \HM{n}$, where
\beq \ell_n : [0,+\infty) \times [0, T) \vain [0,+\infty),~~
(r,t) \mapsto \ell_n(r,t) := 2 K_n \| \fiapp(t) \| \, r + K_n \, r^2~. \feq
\end{prop}
\textbf{Proof.} (i) The bilinearity is obvious, the estimate follows from Eqs.
\rref{hilb} for $\LP$ and \rref{est} for the map $(f, g) \mapsto f \sc \, \partial g$.
\parn
(ii) Use Corollaries \ref{corlip}, \ref{corell} on quadratic maps. \fine
\salto
\textbf{The function $\boma{\Ximn}$.} From here to the end of the paper, we denote in this
way any function
in $C([0,+\infty), [0,+\infty))$ such that
\beq \Ximn \in C([0,+\infty), [0,+\infty))~\mbox{nondecreasing}~, \label{dexin} \feq
$$ \| \xi(t) \|_{n-1} \leqs \Ximn(t) \qquad \mbox{for $t \in [0,+\infty)$} $$
(e.g., $\Ximn(t) := \sup_{s \in [0,t]} \| \xi(s) \|_{n-1}$).
\salto
\textbf{Cauchy and Volterra problems.}
\begin{prop}
\textbf{Definition.}
For any $f_0 \in \HM{n+1}$, \copn is the Cauchy problem \rref{nsm}, i.e.,
$$ \mbox{\textsl{Find}}~
\varphi \in C([0, T), \HM{n+1}) \cap C^1([0, T), \HM{n-1}) ~\mbox{\textsl{such that}} $$
\beq {\dot \varphi(t)} = \Delta \varphi(t) + \PN(\varphi(t), t)
\quad \mbox{\textsl{for all} $t \in [0, T)$}~, \qquad
\varphi(0) = f_0 \label{caupn} ~.\feq
For any $f_0 \in \HM{n}$, \vopn is the Volterra problem
$$ \mbox{\textsl{Find}}~
\varphi \in C([0, T), \HM{n}) \quad \mbox{\textsl{such that}}$$
\beq \varphi(t) = \UD{t}
f_0 + \int_{0}^t~ d s~ \UD{(t - s)} \PN(\varphi(s),s) \quad
\mbox{\textsl{for all} $t \in [0, T)$}~. \label{inttn} \feq
\end{prop}
\begin{rema}
\textbf{Remarks.}
\textbf{(i)} If $f_0 \in \HM{n+1}$, \vopn is equivalent to \copn by the reflexivity of
the Hilbert space $\HM{n-1}$ (see once more Proposition \ref{integral}). \parn
\textbf{(ii)} For any $f_0 \in \HM{n}$, uniqueness and local existence are granted for \vopn (Propo\-sitions
\ref{unique} and \ref{esloc}).
\end{rema}
\section{Results for the NS equations arising from the previous framework.}
\label{resns}
We keep the assumption \rref{xiasin} and all notations of the previous section; furthermore, we fix an initial datum
\beq f_0 \in \HM{n}~. \feq
The analysis of the previous section allows us to
identify \vopn with a Volterra problem of the general type discussed in Section
\ref{quadr}, with semigroup estimators $u, \um$ of the form considered therein and a
quadratic nonlinearity $\PN$. Due to Propositions \ref{edelta}, \ref{isbil}, the
constant $K$ and the functions $u$, $\um$, $\uv$, $\MM$, $\Xim$ of Section \ref{quadr} can be taken as
follows:
\beq  \mbox{$K =$ a constant $K_n$ fulfilling \rref{est}}~, \quad u(t) := e^{-t}~,~\um(t) = \uv(t) e^{-t}~,
\feq
$$ \uv~\mbox{as in \rref{eqww}}~, \quad \MM~\mbox{as in \rref{eqmm}}~,
\quad \Xim = \Ximn~\mbox{as in \rref{dexin}}~. $$
Hereafter, we rephrase Propositions \ref{geprop} and \ref{esglobbm} with the above specifications.
\begin{prop}
\label{gepropn}
\textbf{Proposition.}
Let us consider for \vopn an
approximate solution $\fiap \in C([0,T), \HM{n})$, where $T \in (0,+\infty]$.
Suppose there are
functions $\EE_n, \Dd_n, \RR_n$ $\in$ $C([0,T),$ $[0,+\infty))$ such that (i)-(iii) hold: \parn
(i) $\fiap$ has the integral error estimate
\beq \| E(\fiap(t)) \|_n \leqs \EE_n(t) \qquad \mbox{for $t \in [0,T)$}~; \label{estimatmn} \feq
(ii) one has
\beq \| \fiap(t) \|_n \leqs \Dd_n(t) \qquad \mbox{for $t \in [0,T)$}~; \label{bdmn} \feq
(iii) $\RR_n$ solves the control inequality
\beq \EE_n(t)  + K_n \int_{0}^t \! d s \, \um(t-s) (2 \Dd_n(s) \RR_n(s) + \RR_n^2(s)) \leqs \RR_n(t)
~~\mbox{for $t \in [0,T)$}~. \label{monmn} \feq
Then, (a) and (b) hold: \parn
(a) \vop has a solution $\varphi : [0, T) \vain \HM{n}$; \parn
(b) one has
\beq \| \varphi(t) - \fiap(t) \|_n \leqs \RR_n(t) \qquad \mbox{for $t \in [0, T)$}~. \label{axbbmn} \feq
\end{prop}
\begin{prop}
\label{esglobbmn}
\textbf{Proposition.} Let us consider for \vopn an
approximate solution $\fiap \in C([0,T), \HM{n})$, where $T \in (0,+\infty]$. Suppose
there are
functions $\EE_n, \Dd_n$ $\in C([0,T],$ $[0,+\infty))$ such that (i)-(iii) hold: \parn
(i) $\EE_n$ is nondecreasing, and binds the integral error as in \rref{estimatmn}; \parn
(ii) $\Dd_n$ is nondecreasing and binds $\fiap$ as in \rref{bdmn}; \parn
(iii) one has
\beq 2 \sqrt{K_n \MM(T) \EE_n(T)} + 2 K_n \MM(T) \Dd_n(T) \leqs 1~. \label{coquaamn} \feq
Then \vopn has a solution $\varphi : [0, T) \vain \HM{n}$ and, for all $t \in [0, T)$,
\beq \| \varphi(t) - \fiap(t) \|_n \leqs \RR_n(t) ~, \label{xquaamn} \feq
$$ \hspace{-0.0cm} \RR_n(t):=
\left\{ \barray{ll}
\! \! \dd{1 \! \! - \! \! 2 K_n \MM(t) \Dd_n(t) \! \! - \! \!
\sqrt{(1 \! \! - \! \! 2 K_n \MM(t) \Dd_n(t))^2 \! \! - \! \! 4 K_n \MM(t) \EE_n(t)} \over 2 K_n \MM(t)}
& \! \! \mbox{if $t \in (0,T)$,} \\
\! \! \EE_n(0) & \! \! \mbox{if $t = 0$;} \farray \right.
$$
the above prescription gives a well defined, nondecreasing function $\RR_n$ $\in$ $C([0,T)$, $[0,+\infty))$.
\end{prop}
In Section \ref{quadr}, from Proposition \ref{esglobbm}  we have inferred Proposition \ref{xvxm},
corresponding to the approximate solution $\fiap := 0$; in the present situation, this reads as follows.
\begin{prop}
\label{xvxmn}
\textbf{Proposition.} Let
\beq \FFf_n(t) := \de_n + \Ximn(t) \, \MM(t)~. \label{simi} \feq
Suppose $T \in [0,+\infty]$, and
\beq 4 K_n \MM(T) \FFf_n(T) \leqs 1~. \label{coquaazmn} \feq
Then \vopn has a solution $\varphi : [0, T) \vain \HM{n}$ and,
for all $t \in [0, T)$,
\beq \| \varphi(t) \|_n \leqs \FFf_n(t) \, \XX(4 K_n \MM(t) \FFf_n(t))~ \label{xbbzmn} \feq
(where, as in \rref{xquaazm}:
$\XX(z) := \dd{1 - \sqrt{1 - z} \over (z/2)}$ for $z \in (0,1]$, $\XX(0) := 1$).
\end{prop}
The other results of Section \ref{quadr} were about global existence and
exponential decay, under specific assumption. In the present framework the
constants $B, \sigma, N$ of the cited section are given by
Eq. \rref{holdw};
this allows to rephrase Proposition \ref{esglobb} in this way.
\begin{prop}
\label{esglobbn}
\textbf{Proposition.} Let us consider for \vopn an
approximate solution $\fiap \in C([0,+\infty), \HM{n})$. Suppose there are
constants $E_n, D_n \in [0,+\infty)$ such that: \parn
(i) $\fiap$ admits the integral error estimate
\beq \| E(\fiap)(t) \|_n \leqs E_n e^{-t} \qquad \mbox{for $t \in [0, +\infty)$}~; \feq
(ii) for all $t$ as above,
\beq \| \fiap(t) \|_n \leqs D_n e^{-t}~; \label{bdn} \feq
(iii) one has
\beq 2 {~}^{4}{\! \! \! \sqrt{2}} \sqrt{K_n E_n} + 2 \sqrt{2} K_n D_n \leqs 1~. \label{coquaan} \feq
Then \vopn has a global solution $\varphi : [0, +\infty) \vain \HM{n}$ and, for all $t \in [0, +\infty)$,
\beq \| \varphi(t) - \fiap(t) \|_n \leqs \R_n e^{-t}~, \label{xquaan} \feq
$$ \R_n := {1 - 2 \sqrt{2} K_n D_n - \sqrt{(1 - 2 \sqrt{2} K_n D_n)^2 - 4 \sqrt{2} K_n E_n}
\over 2 \sqrt{2} K_n}~. $$
\end{prop}
The applications of Proposition \ref{esglobb} considered in Section \ref{quadr} were based
on the assumption of exponential decay for the external forcing, that in the present framework
must be formulated in this way:\parn
(\Xp)$_n$ There is a constant $J_{n-1} \in [0, +\infty)$ such that
\beq \| \xi(t) \|_{n-1} \leqs J_{n-1} e^{-2 t} \qquad \mbox{for all $t \in [0,+\infty)$}~. \label{eqxpn} \feq
The above mentioned applications in Section \ref{quadr} were
Proposition \ref{xvx} (corresponding to the choice $\fiap := 0$) and Proposition
\ref{aflow} (with $\fiap$ the $\AA$-flow approximate solution). These can be restated,
respectively, in the following way:
\begin{prop}
\label{xvxn}
\textbf{Proposition.} Assume (\Xp)$_n$, and define
\beq \A_n := \de_n + \sqrt{2} J_{n-1}~; \label{proen} \feq
furthermore, assume
\beq 4 \sqrt{2} K_n \A_n \leqs 1~. \label{coquaazn} \feq
Then \vopn has a global solution $\varphi : [0, +\infty) \vain \HM{n}$ and,
for all $t \in [0, +\infty)$,
\beq \| \varphi(t) \|_n \leqs \A_n \, \XX(4 \sqrt{2} K_n \A_n) \, e^{-t}~ \label{xbbzn} \feq
(with $\XX$ as in \rref{xquaazm}).
\end{prop}
\begin{prop}
\label{aflown}
\textbf{Proposition.} Define
\beq \fiap \in C([0,+\infty), \HM{n})~, \qquad \fiap(t) := \UD{t} f_0 +
\int_{0}^t d s~\UD{(t-s)} \xi(s)~. \feq
Furthermore, let us keep the assumptions and definitions (\Xp)$_n$ \rref{proen}
\rref{coquaazn}. Then the global solution $\varphi : [0, +\infty) \vain \HM{n}$ of \vopn is such that,
for all $t \in [0, +\infty)$,
\beq \| \varphi(t) - \fiap(t) \|_n
\leqs \sqrt{2} K_n \A_n^2 \, \XXX(4 \sqrt{2} K_n \A_n) \, e^{-t}~ \label{xbbhn} \feq
(where, as in \rref{xquaah}:
$\XXX(z) := \dd{1 - (z/2) - \sqrt{1 - z} \over (z^2/8)}$ for $z \in (0,1]$, $\XXX(0) := 1$).
\end{prop}
\begin{rema} \label{remache}
\textbf{Remarks.} \textbf{(i)} Condition \rref{coquaazmn} can be fulfilled with either $f_0$,
$\xi$ small, $T$ large or $f_0$, $\xi$ large, $T$ small. Condition \rref{coquaazn} is fulfilled if $f_0$ and $\xi$
are sufficiently small. \parn
\textbf{(ii)} As a special case, suppose the external forcing
$\xi$ to be identically zero; then we can take $\Ximn =0$ and $J_{n-1} = 0$.
Eqs. \rref{coquaazmn} and \rref{coquaazn} ensure global existence if the datum fulfills
the conditions $4 K_n \MM(+\infty) \de_n \leqs 1$ and
$4 \sqrt{2} K_n \de_n \leqs 1$, respectively. The less restrictive condition on $f_0$ is the second one,
since $\sqrt{2} < \MM(+\infty)$. \parn
\textbf{(iii)} Let us return to Proposition \ref{gepropn}; this states, amongst else, that
the solution $\varphi$ of \vopn exists on $[0,T)$ if the inequality \rref{monmn}
has a solution $\RR : [0,T) \vain [0,+\infty)$. Let us compare this statement with a result presented
in the recent work \cite{Che}, that we rephrase here in our notations. \parn
Let us consider the (incompressible, zero mean)
NS equations with external forcing $\xi$ and initial datum $f_0$;
when we refer to \cite{Che} a solution of this Cauchy problem means a strong solution,
as defined therein. Now suppose  $\fiap$ to be an approximate solution
on an interval $[0,T)$, and
\beq \| \fiap(0) - f_0 \|_{n-1} \leqs \delta_{n-1}~, \qquad
\| \dot \fiap(t) - \Delta \fiap(t) - \PP(\fiap(s), s) \|_{n-1} \leqs \ep_{n-1}(t)~, \feq
$$ \| \fiap(t) \|_{n-1} \leqs \Dd_{n-1}(t)~, \qquad \| \fiap(t) \|_{n} \leqs \Dd_n(t) $$
for all $t \in [0,T)$, for suitable estimators $\delta_{n-1} \geqs 0$, $\ep_{n-1}, \Dd_{n-1}, \Dd_{n} : [0,T)
\vain [0,+\infty)$. According to \cite{Che} (page 065204-10), the NS Cauchy problem has a solution
$\varphi : [0,T) \vain \HM{n}$ if
\beq \delta_{n-1} + \int_{0}^T \! \! \! \! d s \, \ep_{n-1}(s)
< {1 \over C_n T} \, e^{\dd{-C_n \int_{0}^T \! \! \! \! d s \, (\Dd_{n-1}(s) +
\Dd_{n}(s))}} \label{ineche} \feq
(with $C_n > 0$ a constant not computed explicitly, whose role is analogous to the one of $K_n$);
\rref{ineche} is an inequality involving only the approximate solution, and plays a role similar
to our \rref{monmn} to grant the existence of an exact solution on $[0,T)$. \parn
Seemingly, Eq.\rref{ineche} is not suited to obtain
results of global existence for the exact solution. To explain this statement,
suppose $\fiap$ and its estimators to be defined on $[0,+\infty)$,
with $\delta_{n-1} \neq 0$ or $\ep_{n-1}$ non identically zero;
then \rref{ineche} surely fails for large $T$, even in the most favourable situation where all integrals
therein converge for $T \vain +\infty$. In fact,
\beq \mbox{l.h.s. of \rref{ineche}} \uper{$T \vain +\infty$} \, \delta_{n-1} + \int_{0}^{+\infty} d s \, \ep_{n-1}(s)
\in (0,+\infty], \feq
$$ 0 < \mbox{r.h.s. of \rref{ineche}} \leqs {1 \over C_n T} \uper{$T \vain +\infty$} 0~. $$
On the contrary, our control inequality \rref{monmn} can be used in certain cases
to derive the existence of $\varphi$ (and bind its distance from $\fiap$)
up to $T = +\infty$; some applications of this type have appeared in the present section,
further examples will be given in the next one on the Galerkin approximations
({\footnote{To conclude this remark we wish to point out that, under special assumptions, some global existence
results could perhaps be derived from the approach of \cite{Che}, with a different analysis
of the differential inequalities proposed by the authors to infer Eq. \rref{ineche}. A discussion
of this point, and of other interesting features of \cite{Che}, would occupy too much space here.}}).
\end{rema}
\salto
\section{Galerkin approximate solutions of the NS equations.}
\label{gale}
Throughout this section, we consider a set $G$ with the following features:
\beq \G \subset \Zd_0~, \qquad \G~\mbox{finite}~, \qquad k \in \G \Leftrightarrow - k \in \G~. \feq
Hereafter we write $\prec e_k \succ_{k \in G}$ for the linear subspace of $\bb{D}$ spanned by the
functions $e_k$ for $k \in G$. \salto
\textbf{Galerkin subspaces and projections.} We define them as follows.
\begin{prop}
\textbf{Definition.} The Galerkin
subspace and projection corresponding to $\G$ are
\beq \HG := \bb{D}'_{\so} \cap \prec e_k \succ_{k \in G} =
\{ \sum_{k \in \G} v_k e_k~|~v_k \in \complessi^d~, \overline{v_k} = v_{-k}~, k
\sc \, v_k = 0~\mbox{for all $k$} \}~. \label{hg} \feq
\beq \EG : \bb{\DD}'_{\so}~ \vain \HG~, \qquad v =
\sum_{k \in \Zd_0} v_k e_k \mapsto \EG v := \sum_{k \in \G} v_k e_k~. \label{pg} \feq
\end{prop}
It is clear that
\beq \HG \subset \bb{C}^{\infty} \cap \bb{\DD}'_{\so}~, \qquad
\Delta(\HG) \subset \HG~; \label{deg} \feq
\beq \HG \subset \HM{m}~,~~\EG(\HM{m}) = \HG \quad \mbox{for all $m \in \reali$}~. \feq
The following result will be useful in the sequel.
\begin{prop}
\textbf{Lemma.}
Let $n, p \in \reali$, $n \leqs p$ and $v \in \HM{p}$. Then,
\beq \| (1 - \EG) v \|_n \leqs {\| v \|_{p} \over |\G|^{p-n}}~, \qquad
|\G| := \inf_{k \in \Zd_0 \setminus \G} \sqrt{1 + |k|^2}~.\label{egmp} \feq
\end{prop}
\textbf{Proof.} We have $(1 - \EG) v = \sum_{k \in \Zd_0 \setminus G} v_k e_k$, implying
\beq \| (1 - \EG) v \|^2_n = \sum_{k \in \Zd_0 \setminus G} (1 + |k|^2)^n | v_k |^2 =
\sum_{k \in \Zd_0 \setminus G} {(1 + |k|^2)^{p} \over (1 + | k |^2)^{p-n}} | v_k |^2  \feq
$$ \leqs \Big(\sup_{k \in \Zd_0 \setminus G} {1 \over (1 + |k|^2)^{p-n}} \Big)
\Big(\sum_{k \in \Zd_0 \setminus G} (1 + |k|^2)^{p} | v_k |^2 \Big) \leqs
{1 \over |\G|^{2 (p-n)}} \| v \|^2_{p}~, $$
whence the thesis. \fine
\salto
\textbf{Galerkin approximate solutions.} Let
\beq \xi \in C([0,+\infty), \bb{\DD'}_{\so})~, \qquad f_0 \in \bb{\DD'}_{\!\!\so}~ \feq
(of course, in the sequel $\PN(f,t) := -\LP(f \sc \partial f) + \xi(t)$ whenever
this makes sense).
\begin{prop}
\textbf{Definition.} The \textsl{Galerkin
approximate solution} of NS corresponding to $\G$, with external forcing $\xi$ and
datum $f_0$, is the maximal solution
$\varphi^\G_{f_0} \equiv \varphi^\G$
of the following Cauchy problem, in the finite dimensional space $\HG$:
\beq \mbox{Find $\varphi^G \in C^1([0,T_\G), \HG)$ such that} \label{galer} \feq
$$ \dot \varphi^\G(t) = \Delta \varphi^\G(t) + \EG \PN(\varphi^\G(t)) \quad \mbox{for all $t$}~,
\qquad \varphi^G(0) = \EG f_0~. $$
\end{prop}
Of course, "maximal" means that $[0, T_\G)$ is
the largest interval of existence.
In certain cases, one can prove that $T_G = +\infty$ and derive estimates of $\varphi^\G$
(we return on this in the sequel).
In this section we use the functions
\beq \MM~\mbox{as in Eq. \rref{eqmm}}~, \qquad \XX~\mbox{as in Eq. \rref{xquaazm}}~. \feq
Let us consider any real number $m$, and assume $\xi \in C^{0,1}([0,+\infty), \HM{m-1})$.
We will use the notation $\Xi_{m-1}$ to indicate any
function in $C([0,+\infty), [0,+\infty))$ fulfilling Eq. \rref{dexin} with $n$ replaced by
$m$; in the sequel that equation will be referred to as \rref{dexin}$_m$.
When necessary we will make the assumption (\Xp)$_m$ of exponential decay for the
external forcing; this is like (\Xp)$_n$ with $n \vain m$, thus involving a constant
$J_{m-1} \in [0,+\infty)$. \parn
The forthcoming proposition gives estimates on the interval of existence of the
Galerkin solution and on its norm $\|~\|_{m}$,
which are in fact \textsl{independent of} $G$.
\begin{prop}
\label{duepro}
\textbf{Proposition.} Let $m > d/2$, $\xi \in C^{0,1}([0,+\infty), \HM{m-1})$ and $f_0 \in \HM{m}$;
then, (i)(ii) hold. \parn
(i) Define (similarly to \rref{simi})
\beq \FFf_{m}(t) := \de_{m} + \Xi_{m-1}(t) \, \MM(t)~; \label{assuem} \feq
furthermore, let $T \in (0,+\infty]$, and assume the inequality
\beq 4 K_{m} \MM(T) \FFf_m(T) \leqs 1~. \label{assufom} \feq
Then the Galerkin solution $\varphi^G$ with this datum exists on $[0,T)$ and fulfills
\beq \| \varphi^\G(t) \|_{m} \leqs \Dd_{m}(t)  \qquad \mbox{for $t \in [0,T)$}~, \label{esfigm} \feq
\beq \Dd_{m}(t) := \FFf_{m}(t) \XX(4 K_{m} \MM(t) \FFf_{m}(t))~. \label{esfig1m} \feq
(ii) Alternatively, assume (\Xp)$_m$; define (similarly to \rref{proen})
\beq \A_m := \de_{m} + \sqrt{2} J_{m-1}~, \label{assue} \feq
and suppose
\beq  4 \sqrt{2} K_{m} \A_{m} \leqs 1~. \label{assufo} \feq
Then the Galerkin solution $\varphi^G$ with this datum is global, and fulfills
\beq \| \varphi^\G(t) \|_{m} \leqs D_{m} \, e^{-t} \qquad \mbox{for $t \in [0,+\infty)$}~, \label{esfig} \feq
\beq D_{m} := \A_{m} \XX(4 \sqrt{2} K_{m} \A_{m})  \label{esfig1} \feq
(the above equations will be referred to in the sequel as \rref{assuem}$_m$, \rref{assufom}$_m$, etc.).
\end{prop}
\textbf{Proof.} We refer to the framework of Section \ref{quadr} on systems with quadratic
nonlinearities. In the present case $(\Ff, \|~\|) := (\HG, \|~\|_{m})$,
$(\Ff_{\mp}, \|~\|_{\mp} := (\HG, \|~\|_{m \mp 1})$ (we have three copies of the same finite dimensional space,
but equipped with different, though equivalent, norms);
the operator $\AA$ is $\Delta \restriction \HG$, and the bilinear map is
$\EG \PPN : \HG \times \HG \vain  \HG$, $(f,g) \mapsto \EG \PPN(f, g)$; the function
$\xi$ of Section \ref{quadr} is $\EG \xi \in C^{0,1}([0,+\infty), \HG)$; the initial datum is
$\EG f_0 \in \HG$. \parn
For the operator $\Delta \restriction \HG$ we use the same estimates given for
$\Delta$ in Proposition \ref{edelta}, with $n \vain m$; this justifies using the
scheme of Section \ref{quadr} with $\MM$ as in \rref{eqmm} and $B=1$, $N= \sqrt{2}$.
\parn
To estimate $\EG \PPN$, we use the inequalities on $\PPN$ in Proposition
\ref{samest} with $n \vain m$, and the obvious relation $\| \EG \cdot \|_{m-1}
\leqs \| \cdot \|_{m-1}$; this gives $\| \EG \PPN(f, g) \|_{m-1} \leqs K_{m}
\| f \|_{m} \| g \|_{m}$, and so the constant $K$ of Section \ref{quadr} is, in
this case, $K_{m}$.
\parn
For the initial datum $\EG f_0$ and for $\EG \xi$ we use the estimates
\beq \| \EG f_0 \|_{m} \leqs \| f_0 \|_{m}~, \feq
\beq \| \EG \xi(t) \|_{m-1} \leqs \| \xi(t) \|_{m-1} \leqs \Xi_{m-1}(t)~\mbox{or}~ J_{m-1}~e^{-2 t}~; \feq
of course, the bound via $\Xi_{m-1}$ refers to case (i) and the bound via $J_m$ is for case (ii).
Applying to this framework Propositions \ref{xvxm} and \ref{xvx} we get the statements in
(i) and (ii), respectively. \fine
\begin{rema}
\textbf{Remark.} Global existence of $\varphi^\G$ could be proved under much
weaker conditions than the ones in item (ii) of the above proposition. In fact, using
for $\varphi^\G$ an energy balance relation similar to \rref{conse},
one can derive
global existence and boundedness of $\| \varphi^\G(t) \|_{L^2}$
when $f_0$ is \textsl{arbitrary} and the external forcing makes finite both integrals
$\int_{0}^{+\infty} \! \!  d t~ \| \xi(t) \|_{L^2}$,
$\int_{0}^{+\infty} \! \!  d t~ \| \xi(t) \|^2_{L^2}$ : see, e.g., \cite{Tem}.
However, the energetic approach does not allow to derive
estimates of the specific type
appearing in Proposition \ref{duepro}.
\ffine
\end{rema}
\textbf{The distance between the exact NS solution and the Galerkin approximations.}
From here to the end of the paragraph, we fix two real numbers
\beq p \geqs n > {d \over 2}~; \feq
we also fix
\beq \xi \in C^{0,1}([0,+\infty), \HM{p-1}),~\qquad f_0 \in \HM{p} \feq
and denote with $\varphi^\G$ the Galerkin approximate solution with such forcing and
datum, for any $G$ as before. This will be compared with the solution $\varphi$ of the NS equations
with the same forcing and datum.
\begin{prop}
\textbf{Lemma.} Let us regard $\varphi^\G$ as an approximate solution of \vopn;
then the following holds. \parn
(0) The integral error of $\varphi^\G$ is
\beq E(\varphi^\G)(t) = -(1 - \EG)
\left[ \UD{t} f_0 + \int_{0}^t d s \, \UD{(t-s)} \PN(\varphi^\G(s)) \right]~. \label{eint} \feq
(i) Let us introduce the definitions or assumptions  \rref{assuem}$_{p}$
\rref{assufom}$_{p}$, for some $T \in (0,+\infty]$ (implying existence of $\varphi^\G$ on
$[0,T)$).
Then, for all $t \in [0,T)$ we have
\beq \| E(\varphi^\G)(t) \|_n \leqs {\YY_{p}(t) \over |\G|^{p-n}}~, \label{tesem} \feq
\beq \YY_{p}(t) := \FFf_{p}(t) \Big[1 + K_{p} \MM(t)
\FFf_{p}(t) \XX^2(4 K_{p} \MM(t) \FFf_{p}(t)) \Big]~. \label{tese1m} \feq
The function $\YY_p$ is nondecreasing. \parn
(ii) Alternatively,
introduce the definitions or assumptions (\Xp)$_{p}$  \rref{assue}$_{p}$
\rref{assufo}$_{p}$
(implying that $\varphi^\G$ is global). Then, for all $t \in [0,+\infty)$ we have
\beq \| E(\varphi^\G)(t) \|_n \leqs {Y_{p} \over |\G|^{p-n}} \, e^{-t}~, \label{tese} \feq
\beq Y_{p} := \A_{p}\Big[1 + \sqrt{2} K_{p}
\A_{p} \XX^2(4 \sqrt{2} K_{p} \A_{p}) \Big]~. \label{tese1} \feq
\end{prop}
\textbf{Proof.} \textsl{Derivation of \rref{eint}.} By definition
\beq E(\varphi^\G)(t) =
\varphi^\G(t) - \UD{t} f_0 - \int_{0}^t d s \, \UD{(t-s)} \PN(\varphi^\G(s))~; \label{bydef} \feq
on the other hand, the Cauchy problem \rref{galer} defining $\varphi^\G$ has the integral reformulation
\beq \varphi^\G(t) = \UD{t} \EG f_0 + \int_{0}^t d s \, \UD{(t-s)} \EG \PN(\varphi^\G(s))~; \feq
inserting this into \rref{bydef} we get
\beq E(\varphi^\G)(t) = - \UD{t} (1 - \EG)  f_0 - \int_{0}^t d s \, \UD{(t-s)} (1 - \EG)\PN(\varphi^\G(s))~. \feq
Finally, the operator $1 - \EG$ commutes with $\Delta$ and its semigroup (as made evident
by the Fourier representations); so, $1 - \EG$ can be factored out  and we obtain the thesis \rref{eint}. \parn
\textsl{Some preliminaries to the proof of (i) and (ii).} From \rref{eint}, the
estimates \rref{egmp} on $1 - \EG$ and \rref{reguff}--\rref{eqww} on the semigroup of $\Delta$ we get \parn
{\vbox{
\beq \| E(\varphi^\G)(t) \|_n \leqs {1 \over |\G|^{p-n}}
\left[ \| \UD{t} f_0 \|_{p} + \int_{0}^t d s \, \| \UD{(t-s)} \PN(\varphi^\G(s)) \|_{p} \right]
\label{intoe} \feq
$$ \leqs {1 \over |\G|^{p-n}}
\left[ e^{-t} \| f_0 \|_{p} + \int_{0}^t d s \, e^{-(t-s)} \uv(t-s) \| \PN(\varphi^\G(s)) \|_{p-1} \right]~.
$$}}
\textsl{Proof of (i).} From $\PN(\varphi^\G(s)) = \PPN(\varphi^\G(s), \varphi^\G(s)) +
\xi(s)$ we infer the following, for $s \in (0,t)$:
\beq \| \PN(\varphi^\G(s)) \|_{p-1} \leqs K_{p} \| \varphi^\G(s) \|^2_{p} + \| \xi(s) \|_{p - 1}
\label{weins} \feq
$$ \leqs
K_{p} \Dd^2_{p}(s) + \Xi_{p -1}(s) \leqs K_{p} \Dd^2_{p}(t) + \Xi_{p -1}(t) $$
(in the above, we have used  \rref{esfigm}$_p$ \rref{dexin}$_{p}$ and the relation
$\Dd_{p}(s) \leqs \Dd_{p}(t)$). \parn
We insert the result \rref{weins} into Eq. \rref{intoe}; in this way
we are left with an integral
$\int_{0}^t d s \, e^{-(t-s)} \uv(t-s)  = \int_{0}^t d s \, e^{-s} \uv(s)
\leqs \MM(t)$. From this bound and $e^{-t} \leqs 1$ we obtain
\beq \| E(\varphi^\G(t))\|_n \leqs { 1 \over |\G|^{p-n}}
\Big[ \| f_0 \|_{p} + \MM(t) \Big(K_{p} \Dd^2_{p}(t) + \Xi_{p-1}(t)\Big)~\Big] \feq
$$ = { 1 \over |\G|^{p-n}}
\Big[ \FFf_{p}(t) + \MM(t) K_{p} \Dd^2_{p}(t)~\Big]~, $$
where the last passage follows from definition \rref{assuem}$_{p}$;
now, explicitating $\Dd_{p}(t)$ we get the thesis \rref{tesem} \rref{tese1m}.
Finally, $\YY_p$ is nondecreasing because $\FF_p$, $\MM$ and $\XX$ are so. \parn
\textsl{Proof of (ii).} In this case, from the inequality
$ \| \PN(\varphi^\G(s)) \|_{p-1}$ $\leqs K_{p} \| \varphi^\G(s) \|^2_{p}$ $ + \| \xi(s) \|_{p - 1}$
we infer, by means of Eqs. \rref{esfig}$_{p}$ and (\Xp)$_p$~,
\beq \| \PN(\varphi^\G(s)) \|_{p-1} \leqs K_{p} D^2_{p} e^{-2 s} + J_{p - 1} e^{-2 s}~. \label{tns} \feq
Inserting this result into \rref{intoe} we are left with a term
$e^{-t} \int_{0}^t d s \, e^{-s} \uv(t-s)$, which is bounded by $\sqrt{2} \, e^{-t}$ due to \rref{defn};
the conclusion is
\beq \| E(\varphi^\G(t))\|_n \leqs {e^{-t} \over |\G|^{p-n}}
\Big[ \| f_0 \|_{p} + \sqrt{2} (K_{p} D^2_{p} + J_{p-1})~\Big] \feq
$$ = {e^{-t} \over |\G|^{p-n}}
\Big[ \A_{p} + \sqrt{2} K_{p} D^2_{p}~\Big] $$
(the last equality following from \rref{assue}$_{p}$).
Now, explicitating $D_{p}$ we get the thesis \rref{tese} \rref{tese1}. \fine
The following proposition contains the main result of the section.
\begin{prop}
\label{prgal}
\textbf{Proposition.}
(i) Let $T \in (0,+\infty]$; make the
assumptions and definitions \rref{assuem}$_{n}$ \rref{assufom}$_n$ and
\rref{assuem}$_{p}$ \rref{assufom}$_{p}$ (implying the existence of $\varphi^G$ on
$[0,T)$). Finally, with
$\Dd_n$ and $\YY_{p}$ defined by \rref{esfig1m}$_n$ \rref{tese1m},
assume
\beq 2 \sqrt{{K_n \MM(T) \YY_{p}(T) \over |G|^{p-n} }} + 2 K_n \MM(T) \Dd_n(T) \leqs 1~. \label{takesm} \feq
Then \vopn has a solution $\varphi$ of domain $[0,T)$ and, for all $t$ in this interval,
\beq \| \varphi(t) - \varphi^\G(t) \|_n \leqs {\WW_{n p \, |G|}(t) \over |G|^{p-n}}~, \label{fiig} \feq
\beq \WW_{n p \, \, |G|}(t) := {\YY_{p}(t) \over 1 - 2 K_n \MM(t) \Dd_n(t) } \,
\XX \left({4 K_n \MM(t) \YY_{p}(t)  \over (1 - 2 K_n \MM(t) \Dd_n(t))^2 |G|^{p-n}}\right)~. \label{wwnp} \feq
The function $t \mapsto \WW_{n p \, \, |G|}(t)$ is nondecreasing; a rough, $|G|$-independent bound for it is
\beq \WW_{n p \, \, |G|}(t) \leqs {2 \YY_{p}(T) \over 1 - 2 K_n \MM(T) \Dd_n(T) } \label{rough} \feq
for all $t \in [0,T)$. \parn
(ii) Alternatively, make the
assumptions and definitions (\Xp)$_n$  \rref{assue}$_{n}$ \rref{assufo}$_n$
and (\Xp)$_{p}$  \rref{assue}$_{p}$ \rref{assufo}$_{p}$ (implying global existence of
$\varphi^G$). Finally,
with $D_n$ and $Y_{p}$ defined by \rref{esfig1}$_n$ \rref{tese1}, assume
\beq 2 {~}^{4}{\! \! \! \sqrt{2}} \sqrt{{K_n Y_{p} \over |G|^{p-n} }} + 2 \sqrt{2} K_n D_n \leqs 1~. \label{takes} \feq
Then \vopn has a solution $\varphi$ of domain $[0,+\infty)$ and, for all $t$ in this interval,
\beq \| \varphi(t) - \varphi^\G(t) \|_n \leqs {W_{n p \, |G|} \over |G|^{p-n}} e^{-t}~, \label{ffiigg} \feq
\beq W_{n p \, \, |G|} := {Y_{p} \over 1 - 2 \sqrt{2} K_n D_n} \,
\XX \left({4 \sqrt{2} K_n Y_{p} \over (1 - 2 \sqrt{2} K_n D_n)^2 |G|^{p-n}}\right)~. \label{wwnnpp} \feq
The above constant has the rough, $|G|$-independent bound
\beq W_{n p \, \, |G|} \leqs {2 Y_{p} \over 1 - 2 \sqrt{2} K_n D_n}~. \label{rou} \feq
\end{prop}
\textbf{Proof.}
(i) A simple application of Proposition \ref{esglobbmn}, with
\beq \fiap = \varphi^{\G}~;~~\Dd_{n}~\mbox{as in \rref{esfig1m}}~;~~
\EE_n(t) = {\YY_{p}(t) \over |G|^{p-n}}~. \feq
The condition \rref{coquaamn} in the cited proposition takes the form \rref{takesm}. The proposition
ensures that $\varphi$ is defined on $[0,T)$, and gives the estimate
\beq \| \varphi(t) - \varphi^\G(t) \|_n \leqs \RR_n(t)~, \label{xybbm} \feq
involving the nondecreasing, continuous function
$$ \RR_n(t) := {1 - 2 K_n \MM(t) \Dd_n(t) -
\sqrt{(1 - 2 K_n \MM(t) \Dd_n(t))^2 - {4 K_n \MM(t) \YY_{p}(t)/|G|^{p-n}}} \over
2 K_{n} \MM(t)}~, $$
\beq ~~\RR_n(0) := {\| f_0 \|_{p} \over | \G |^{p-n}}~. \feq
We note that we can write
\beq \RR_n(t) = {\YY_{p}(t) \over (1 - 2 K_n \MM(t) \Dd_n(t)) | G |^{p-n}}
\XX \left({4 K_n \MM(t) \YY_{p}(t) \over (1 - 2 K_n \MM(t) \Dd_n(t))^2 |G|^{p-n}}\right)~; \feq
this yields the thesis \rref{fiig} \rref{wwnp}. \parn
The fact that $\WW_{n p |G|}$ is a nondecreasing function of time is apparent from its definition.
The bound \rref{rough} for it follows from the nondecreasing
nature of the function $t  \mapsto \YY_{p}(t)/(1 - 2 K_n \MM(t) \Dd_n(t))$ and
from the inequality $\XX(z) \leqs 2$ for all $z \in [0,1]$).
\parn
(ii) A simple application of Proposition \ref{esglobbn}, with
\beq \fiap = \varphi^{\G}~, \qquad D_n~\mbox{as in \rref{esfig1}}~, \qquad E_n = {Y_{p} \over |G|^{p-n}}~. \feq
The condition \rref{coquaan} in the cited proposition takes the form \rref{takes}. The proposition
ensures that $\varphi$ is defined on $[0,+\infty)$, and gives the estimate
\beq \| \varphi(t) - \varphi^\G(t) \|_n \leqs \R_n e^{-t}~, \label{xybb} \feq
\beq \R_n := {1 - 2 \sqrt{2} K_n D_n - \sqrt{(1 - 2 \sqrt{2} K_n D_n)^2 - {4 \sqrt{2} K_n Y_{p}/|G|^{p-n}}} \over
2 \sqrt{2} K_{n}}~. \label{xxquaa} \feq
We note that we can write
\beq \R_n = {Y_{p} \over (1 - 2 \sqrt{2} K_n D_n) | G |^{p-n}}
\XX \left({4 \sqrt{2} K_n Y_{p} \over (1 - 2 \sqrt{2} K_n D_n)^2 |G|^{p-n}}\right)~, \feq
yielding the thesis \rref{ffiigg} \rref{wwnnpp}; the rough bound \rref{rou} follows
again from the inequality $\XX(z) \leqs 2$. \fine
\begin{rema}
\textbf{Remark.} Of course, if $p > n$ the previous proposition implies convergence of the Galerkin
solution $\varphi^G$ to the exact solution $\varphi$ of \vopn. More precisely,
with the assumptions in (i) we infer from \rref{takesm} \rref{rough} that
\beq \sup_{t \in [0,T)} \| \varphi(t) - \varphi^\G(t) \|_n = O({1 \over |G|^{p-n}}) \vain 0 \qquad \mbox{for
$|G| \vain + \infty$}~; \feq
in case (ii), we infer from \rref{ffiigg} \rref{rou} that
\beq \sup_{t \in [0,+\infty)} e^t \, \| \varphi(t) - \varphi^\G(t) \|_n = O({1 \over |G|^{p-n}}) \vain 0 \qquad \mbox{for
$|G| \vain + \infty$}~. \feq
\end{rema}
\section{Numerical examples.}
\label{nume}
Given the necessary constants $K_n$, the datum norms and some bounds on the external forcing,
the framework of Sections \ref{resns} and \ref{gale} yields informations on
the time of existence of the solution $\varphi$ of \vopn, and on its $\HM{n}$ distance from an
approximate solution. In the sequel we exemplify such estimates referring to Section \ref{gale}, i.e., to
the Galerkin approximations. \parn
Throughout the section, we take
\beq d = 3~; \qquad n = 2~, \quad p = 4~. \feq
The constants $K_2$ and $K_4$ involved in calculations can be obtained from
Lemmas \ref{lemsi} \ref{lemci} and Proposition \ref{multi}; a MATLAB computation
illustrated in Appendix \ref{appekdue} yields the values
\beq K_2 = 0.20~, \qquad K_4 = 0.067~. \feq
The other calculations mentioned hereafter have been performed using MATHEMATICA.
\salto
\textbf{An application of Proposition \ref{prgal}, item (i).}
We suppose the external forcing has bounds \rref{dexin}$_2$
\rref{dexin}$_4$ with $\Xi_1(t) =$ const. $\equiv \Xi_1$ and
$\Xi_3(t) =$ const. $\equiv \Xi_3$ for all $t \in [0,+\infty)$. Conditions
\rref{assufom}$_2$ and \rref{assufom}$_4$ are satisfied with $T = +\infty$ if
\beq \| f_0 \|_2 + 1.88 \, \Xi_1 < 0.667~, \qquad  \| f_0 \|_4  + 1.88 \, \Xi_3 < 1.99~; \label{satisf} \feq
under the above inequalities for $f_0$ and the forcing, the Galerkin solution
$\varphi^\G$ exists on $[0,+\infty)$ for each $G$. As an example, conditions \rref{satisf} are satisfied
in the case
\beq \| f_0 \|_2 = 0.15~,~~\| f_0 \|_4 = 1.50~, \quad
\Xi_1 = 0.025~,~~\Xi_3 = 0.25~, \label{satisfe} \feq
to which we stick hereafter. In the above case, condition \rref{takesm} with $n=2$, $p=4$ and
$T = +\infty$ becomes $0.161 + 2.31/|G| \leqs 1$, which is fulfilled if
\beq |G| \geqs 2.76~. \label{gie} \feq
Assuming \rref{satisfe} \rref{gie} the solution $\varphi$ of
\vopdue is also global, and
\beq \| \varphi(t) - \varphi^\G(t) \|_2 \leqs {\WW_{2 4 \, |G|}(t) \over |G|^2} \leqs
{8.71 \over |G|^2} \qquad \mbox{for all $G$ as above,
$t \in [0,+\infty)$}~. \feq
The numerical value of $\WW_{2 4 \, |G|}(t)$ can be computed at will from
definition \rref{wwnp}; here we have used the rough bound
$\WW_{2 4 \, |G|}(t) \leqs 8.71$, coming from
\rref{rough}.
\salto
\textbf{Another application of Proposition \ref{prgal}, item (i).}
We maintain the assumptions $\Xi_1(t) =$ const. $\equiv \Xi_1$ and
$\Xi_3(t) =$ const. $\equiv \Xi_3$ for all $t \in [0,+\infty)$. We take
\beq \| f_0 \|_2 = 0.20~,~~\| f_0 \|_4 = 2.00~, \qquad \Xi_1, \Xi_3~~
\mbox{as in \rref{satisfe}}. \label{satisfef} \feq
Now conditions \rref{satisf} are not fulfilled, indicating that
\rref{assufom}$_2$ and \rref{assufom}$_4$ are not satisfied with $T = +\infty$. On the contrary,
\rref{assufom}$_2$ and \rref{assufom}$_4$ are found to hold with
\beq T = 1.51~, \feq
i.e., the Galerkin solution $\varphi^\G$ exists for any $G$ on the time interval $[0,1.51)$. To go on,
we note that condition \rref{takesm} with $n=2$, $p=4$ and
$T$ as above becomes $0.163 + 2.41/|G| \leqs 1$, which is fulfilled if
\beq |G| \geqs 2.88~. \label{giee} \feq
Under the assumption \rref{giee} the solution $\varphi$ of
\vopdue exists on the same interval, and
\beq \| \varphi(t) - \varphi^\G(t) \|_2 \leqs {\WW_{2 4 \, |G|}(t) \over |G|^2} \leqs
{11.1 \over |G|^2} \qquad \mbox{for all $G$ as above, $t \in [0,1.51)$}~. \feq
Again, we can compute the numerical value of $\WW_{2 4 \, |G|}(t)$ from
the definition \rref{wwnp}; here we have used
the bound $\WW_{2 4 \, |G|}(t) \leqs 11.1$, coming from \rref{rough}.
\salto
\textbf{An application of Proposition \ref{prgal}, item (ii).} Let us recall that
this case refers to exponentially decaying forcing. From the datum norms $\| f_0 \|_m$
and the constants $J_{m-1}$ in the forcing bounds, as in \rref{assue} we define the
coefficients $\A_m := \| f_0 \|_m + \sqrt{2} J_{m-1}$ for $m=2,4$. Conditions \rref{assufo}
for $m=2,4$ become, respectively, $1.14 \A_2 \leqs 1$ and $0.380 \A_4 \leqs 1$; these are fulfilled if
\beq \A_2 \leqs 0.877~, \qquad \A_4 \leqs 2.63~, \feq
and in this case the Galerkin solution $\varphi^\G$ is global for each $G$.
As an example, let us suppose
\beq \A_2 = 0.20~, \qquad \A_4 = 2.00~; \feq
then, condition \rref{takes} becomes $0.121 + 1.75/|G| \leqs 1$, which is fulfilled if
\beq | G | \geqs 2.00~. \feq
With these assumptions the exact solution $\varphi$ of \vopdue is global, and
\beq \| \varphi(t) - \varphi^\G(t) \|_2 \leqs {W_{2 4 \, |G|} \, e^{-t} \over |G|^2}
\leqs {6.10 \, e^{-t} \over |G|^2}
~~\mbox{for all $G$ as above, $t \in [0,+\infty)$}~. \feq
The expression of $W_{2 4 \, |G|}$ is provided by \rref{wwnnpp}; here we have used
the rough bound $W_{2 4 \, |G|} \leqs 6.10$, coming from \rref{rou}.
\vskip 1cm\noindent
\appendix
\section{Appendix. Proof of Lemma \ref{lemmaz}.}
\label{appez}
First of all, we put
\beq Z := \sup_{t \in [t_0, \tau]} z(t)~; \feq
we continue in two steps.
\salto
\textsl{Step 1. For all $k \in \naturali$, one has}
\beq z(t) \leqs Z~ {\Lambda^k \Gamma(\esp)^k (t - t_0)^{k \esp} \over \Gamma(k \esp + 1)}
\qquad \mbox{\textsl{for $t \in [t_0, \tau]$}}~.  \label{claim} \feq
To prove this, we write \rref{claim}$_k$ for the above equation at order $k$, and
proceed by recursion. Eq. \rref{claim}$_0$ is just the inequality $z(t) \leqs Z$.
Now, we suppose that \rref{claim}$_k$ holds and infer from it Eq. \rref{claim}$_{k+1}$. To this purpose,
we substitute \rref{claim}$_k$ into the basic inequality \rref{diseqz}, which gives
$$ z(t) \leqs Z \Lambda^{k+1} {\Gamma(\esp)^k \over \Gamma(k \esp + 1)}~
\int_{t_0}^t~d s~{(s - t_0)^{k \esp} \over (t - s)^{1 - \esp}} ~;
$$
expressing the integral via the known identity \rref{hand}, we get the thesis
\rref{claim}$_{k+1}$.
\parn
\textsl{Step 2. $z(t) =0$ for all $t \in [t_0, \tau]$.} In Eq. \rref{claim}, let us
fix $t$ and send $k$ to $\infty$; the right hand side of this inequality vanishes
in this limit, yielding the thesis. \fine
\section{Appendix. A scheme to solve numerically the control inequality \rref{monm}.}
\label{outli}
\textbf{Notations.}
In this Appendix we often write $\{0,...,\Mm \}$ where
$\Mm$ is an integer or $+\infty$. If $\Mm$ is a nonnegative integer this
will mean, as usually, the set of integers $0,1,2....,\Mm$. If $\Mm$ is a negative
integer, we will intend $\{0,...,\Mm \} := \emptyset$. If $\Mm = +\infty$,
$\{0,...,\Mm\}$ will mean the set $\naturali$ of all natural numbers. \parn
We often consider finite or infinite sequences of real numbers
of the form $(t_m)_{m \in \{0,...,\Mm\}}$; if $\Mm = +\infty$, we intend
$t_{\Mm} := \lim_{m \vain +\infty} t_m$ whenever the limit exists.
\vskip 0.1cm \noindent
\textbf{The numerical scheme.}
Let us be given an approximate solution $\fiap \in C([0,T'), \Ff)$ of
\vop, where $T' \in (0,+\infty]$; in the sequel we always intend
\beq \EE(t) := \| E(\fiap)(t) \|~, \qquad \Dd(t) := \| \fiap(t) \| \qquad \mbox{for $t \in [0, T')$}~. \feq
Hereafter we outline a numerically implementable algorithm to construct a
solution $\RR$ of the integral inequality \rref{monm} on some interval $[0,T)
\subset [0,T')$; this solution $\RR$ will be piecewise linear. \parn
In order to construct the algorithm,
we choose a sequence of instants $(t_m)_{m=0,...,\EMP}$, where
$\EMP$ is a positive integer or $+\infty$. We assume
\beq 0 = t_0 < t_1 < t_2 < ...< t_{\EMP} = T'~.  \feq
Furthermore, we denote with $\EE_m$, $\Dd_m$, $H_{m k}$, $I_{m k}$, $N_{m k}$
some constants such that
\beq \sup_{t \in [t_m, t_{m+1})}
\EE(t) \leqs \EE_m~, \sup_{t \in [t_m, t_{m+1})} \Dd(t) \leqs \Dd_m~; \label{edk} \feq
$$ \sup_{t \in [t_m, t_{m+1})} \int_{t_k}^{t_{k+1}} \! \! \! \! \! \! d s~\um(t-s)
\left({s - t_k \over t_{k+1} - t_k}\right)^2 \leqs H_{m k},
\sup_{t \in [t_m, t_{m+1})} \int_{t_k}^{t_{k+1}} \! \! \! \! \! \! d s~\um(t-s)
{s - t_k \over t_{k+1} - t_k} \leqs I_{m k}~, $$
$$ \sup_{t \in [t_m, t_{m+1})} \int_{t_k}^{t_{k+1}} \! \! \! \! \! \! d s~\um(t-s) \leqs N_{m k}~
\quad \mbox{for $m \in \{1,...,\EMP-1\}, k \in \{0,...,m-1\}$}~; $$
$$ \sup_{t \in [t_m, t_{m+1})} \int_{t_m}^{t} \! \! \! \! \! d s~\um(t-s)
\left({s - t_m \over t_{m+1} - t_m}\right)^2 \leqs H_{m m},
\sup_{t \in [t_m, t_{m+1})} \int_{t_m}^{t} \! \! \! \! \! d s~\um(t-s)
{s - t_m \over t_{m+1} - t_m} \leqs I_{m m}~, $$
\beq \sup_{t \in [t_m, t_{m+1})} \int_{t_m}^{t} \! \! \! \! d s~\um(t-s) \leqs N_{m m}~
\quad \mbox{for $m \in \{0,...,\EMP-1\}$}~. \label{b6} \feq
Finally, for $m, k$ as above and all $a, x \in \reali$, we define
\beq \Phi_{m k}(a, x) := (H_{m k}  + N_{m k} - 2 I_{m k}) a^2 + 2 (I_{m k} - H_{m k}) a x +
H_{m k} x^2  \label{b7} \feq
$$ + 2 (N_{m k} - I_{m k}) \Dd_{k} a + 2 I_{m k} \Dd_k x ~.  $$
\begin{prop}
\textbf{Proposition.} Suppose there is a finite or infinite sequence
of nonnegative
reals $(\RR_m)_{m \in \{0,..., \EM \}}$ (with $1 \leqs \EM \leqs \EMP$) such that
\beq \EE_{m} + K \sum_{k=0}^{m} \Phi_{m k}(\RR_k, \RR_{k+1}) \leqs \min(\RR_m, \RR_{m + 1})
\qquad \mbox{for $m \in \{0,...., \EM - 1\}$}~. \label{recur} \feq
Let  $\RR \in C([0,t_\EM), [0,+\infty))$ be the unique piecewise
linear map with values $\RR_m$ at the times $t_m$, i.e.,
\beq \RR(t) = \RR_m + (\RR_{m+1} - \RR_m) {t - t_m \over t_{m+1} - t_m}
~~\mbox{for $t \in [t_m, t_{m+1})$, $m \in \{0,..., \EM-1\}$}~. \feq
Then, $\RR$ solves the integral inequality \rref{monm} on $[0,t_\EM)$.
\end{prop}
\textbf{Proof.} Let $\RR$ be defined as above, and $t$ in some subinterval
$[t_m, t_{m+1})$ ($m \in \{0,..., \EM -1\}$). Then
\beq \mbox{l.h.s. of \rref{monm}} \leqs \EE_{m} + K  \int_{0}^{t}
\! d s \, \um(t-s) (2 \Dd(s) + \RR(s)) \RR(s) \feq
$$ \leqs \EE_{m} + K  \left( \sum_{k=0}^{m-1} \int_{t_k}^{t_{k+1}} + \int_{t_m}^{t} \right)
\! d s \, \um(t-s) (2 \Dd_k + \RR(s)) \RR(s) $$
$$ = \EE_{m} + K \left( \sum_{k=0}^{m-1} \int_{t_k}^{t_{k+1}} + \int_{t_m}^{t} \right) \! d s \, \um(t-s)
\left(2 \Dd_k + \RR_k + (\RR_{k+1} - \RR_k) {s - t_k \over t_{k+1} - t_k} \right) \times $$
$$ \times \left(\RR_k + (\RR_{k+1} - \RR_k) {s - t_k \over t_{k+1} - t_k}  \right) \leqs
\EE_{m} + K \sum_{k=0}^{m} \Phi_{m k}(\RR_k, \RR_{k+1})~, $$
the last passage following from the inequalities \rref{b6} and the definition \rref{b7}
of $\Phi_{m k}$. From here and from \rref{recur} we infer, for $t$ in the same interval,
\beq \mbox{l.h.s. of \rref{monm}} \leqs
\min(\RR_m, \RR_{m+1}) \leqs \RR(t) = \mbox{r.h.s. of \rref{monm}} \feq
In conclusion, \rref{recur} ensures $\RR$ to fulfill \rref{monm} on $[0,t_{\EM})$.
\fine
\begin{rema}
\label{re}
\textbf{Remarks. (i)} A sequence of constants $(\EE_k)$ fulfilling the first
inequality \rref{edk} is
easily obtained if $\fiap \in C([0,T'), \Fp) \cap C^1([0,T'), \Fm)$,
and there are suitable estimators for (the semigroup and)
for the datum and differential errors $d(\fiap)$, $e(\fiap)$.
More precisely suppose that
\beq \| d(\fiap) \| \leqs \delta~, \label{supp2} \feq
and that, for $m \in \{0,...,\EMP-1\}$,
\beq \sup_{t \in [t_m, t_{m+1})} u(t) \leqs u_m~, \qquad
\sup_{t \in [t_m, t_{m+1})} \nom{ e(\fiap(t))}  \leqs \ep_m~, \label{supp3} \feq
$(u_m)_{m=0,...,\EMP-1}$ and $(\ep_m)_{m=0,...,\EMP-1}$ being sequences of nonnegative reals. \parn
From Lemma \ref{lap} on the integral error we obtain, for $t \in [t_m, t_{m+1})$,
\beq \EE(t) \leqs u(t) \delta +
\left(\sum_{k=0}^{m-1} \int_{t_k}^{t_{k+1}} + \int_{t_m}^t \right)
\! d s \, \um(t-s) \nom{ e(\fiap(s))}~; \feq
now, from \rref{supp3} and \rref{b6} we infer, for
$m \in \{0,1,...,\EMP-1\}$,
\beq \EE(t) \leqs \EE_m~~\mbox{if $t \in [t_m, t_{m+1})$}~, \qquad
\EE_m := u_m \delta + \sum_{k=0}^{m} N_{m k} \ep_k~. \feq
(In fact, one could extend this result to the case where $\fiap$ is
continuous from $[t_0, T')$ to $\Fp$ and \textsl{piecewise} $C^1$ from
$[t_0, T')$ to $\Fm$: this typically occurs for the approximate solutions
defined by finite difference schemes in time). \parn
\textbf{(ii)} For any $m$, Eq. \rref{recur} holds if and only if
\beq \mbox{either}~\RR_{m + 1} \in [0, \RR_{m})~,
~~\EE_{m} + K \sum_{k=0}^{m} \Phi_{m k}(\RR_k, \RR_{k+1}) \leqs \RR_{m +1}~, \label{one} \feq
\beq \mbox{or}~~ \RR_{m + 1} \in [\RR_{m}, + \infty)~,
~~\EE_{m} + K \sum_{k=0}^{m} \Phi_{m k}(\RR_k, \RR_{k+1}) \leqs \RR_{m}~; \label{two} \feq
note that each $\Phi_{m k}$ in these inequalities is a quadratic polynomial.
For $m=0$, Eq.s \rref{one} \rref{two} define a problem for
two unknowns $(\RR_0, \RR_1)$; for $m > 0$, we can see \rref{one} \rref{two} as a problem
to determine recursively $\RR_{m+1}$
from $\RR_0, ..., \RR_{m}$. \parn
For each $m$, if the problem has solutions it seems convenient to choose for $\RR_{m+1}$
the smallest admissible value. This criterion could be applied for $m=0$ as well,
choosing among all solutions $(\RR_0, \RR_1)$ the one with the smallest $\RR_1$. \parn
\textbf{(iii)} In practical computations, the determination of $\RR_{m+1}$ from $\RR_0, ..., \RR_m$
goes on until problem \rref{one} or \rref{two} have solutions.
The iteration ends if, for some finite $\EM$, both the above inequalities for $\RR_{\EM + 1}$ have no solutions.
In this case, we have a function $\RR$ solving \rref{monm} on
the interval $[0, t_{\EM})$. Alternatively, the iteration might go on indefinitely. \parn
\textbf{(iv)} The recursive scheme \rref{recur} has a typical feature of the iterative methods
to solve integral equations or inequalities of the Volterra type: to find $\RR_{m+1}$
one must compute a ''memory term'' involving $\RR_0,...\RR_{m}$. The memory term depends nontrivially
on $m$ (through the coefficients $I_{m k}$, etc.), so it must be fully redetermined at each step;
this makes the computation more and more expensive while $m$ grows. \parn
An exception to this framework occurs if the semigroup estimator $t \mapsto \um(t)$ is (a constant $\times$)
an exponential, at least for $t$ greater than some fixed time $\te$; this is just the case of the NS equations,
see the forthcoming Remark \ref{remns} (ii). In this special situation, Eq. \rref{recur} can be rephrased
as a pair of recursion relations  for two real sequences $(\RR_m)$, $(\SS_m)$; at each step,
computation of $\RR_{m+1}$ and $\SS_{m+1}$ does not involve the whole previous history, but only
the values of $\SS_m$ and of $\RR_k$ for $t_m - \te < t_{k+1} \leqs t_{m+1}$.
The forthcoming Proposition explains all this.
\end{rema}
\begin{prop}
\label{special}
\textbf{Proposition.} Let us suppose there are $\te \geqs 0$ and $A, B >  0$ such that
\beq \um(t) = A e^{- B t} \qquad \mbox{for $t \in (\te, + \infty)$}~; \feq
furthermore, let us intend that $k$ always means an integer in $\{0,....,\EMP - 1\}$.
Then, i) ii) hold. \parn
i) Let
\beq m \in \{0,...,\EMP\}~, \quad t_{k+1} \leqs t_{m} - \te~; \feq
then, conditions \rref{b6} are fulfilled with
$$ H_{m k} := A e^{- B t_{m}} H_{k}, \quad H_k :=
{2(e^{B t_{k+1}} - e^{B t_k}) - 2 B e^{B t_{k+1}} (t_{k+1} - t_k) + B^2 e^{B t_{k+1}} (t_{k+1} - t_k)^2
\over B^3 (t_{k+1} - t_k)^2}~; $$
$$ I_{m k} := A e^{- B t_{m}} I_{k}~, \qquad I_k := { e^{B t_k} - e^{B t_{k+1}} + B
e^{B t_{k+1}} (t_{k+1} - t_k) \over B^2 (t_{k+1} - t_k)}~; $$
\beq N_{m k} := A e^{- B t_{m}} N_{k}~, \qquad N_k := {e^{B t_{k+1}} - e^{B t_k} \over B}~.
\label{bb6} \feq
Consequently, for all real $a, x$ one has
\beq \Phi_{m k}(a, x) = A e^{- B t_m} \Phi_k(a, x)~, \label{bb7} \feq
$$ \Phi_{k}(a, x) := (H_{k}  + N_{k} - 2 I_{k}) a^2 + 2 (I_{k} - H_{k}) a x +
H_{k} x^2
+ 2 (N_{k} - I_{k}) \Dd_{k} a + 2 I_{k} \Dd_k x ~.  $$
ii) Consider a sequence $(\RR_{m})_{m \in \{0, ...,\EM\}}$ of nonnegative reals. Then,
$(\RR_{m})$ fulfills Eq. \rref{recur} if and only if there is a sequence
of reals $(\SS_{m})_{m \in\{0,...,\EM-1\}}$ such that
\beq \SS_{m} +
\ski \sum_{\{k | t_{m} - \te < t_{k+1} \leqs t_{m+1} - \te \}} \ski \Phi_{k}(\RR_k, \RR_{k+1})
\leqs \SS_{m+ 1} \quad \mbox{for $m \in\{0,...,\EM-2\}$}~;  \label{recur2} \feq
\beq \EE_{m} +
K \ski \sum_{\{ k | t_{m} - \te < t_{k+1} \leqs t_{m+1} \}} \ski \Phi_{m k}(\RR_k, \RR_{k+1})
+ K A e^{- B t_m}\SS_{m} \leqs \min(\RR_m, \RR_{m + 1}) \label{recur1} \feq
$$ \mbox{for $m \in \{0,...,\EM-1\}$}~. $$
iii) In particular, suppose the instants $t_{m}$ and $\te$ to be integer multiples of
a basic spacing $\tau > 0$:
\beq t_{m} = m \tau \quad \mbox{for $m \in \{0,...,\EMP\}$}~;
\qquad \te = \h \tau \quad \mbox{for some $\h \in \naturali$}~. \feq
Then, for $m \in \{0,...,\EM - 1\}$,
\beq \{ k | t_{m} - \te < t_{k+1} \leqs t_{m+1} \} = \left\{
\barray{ll} \{ m - \h, ..., m \} & \mbox{if $m \geqs \h$}~, \\
\{ 0,...,m \} & \mbox{if $m < \h$}~, \farray \right. \feq
\beq \{ k | t_{m} - \te < t_{k+1} \leqs t_{m+1} - \te \} = \left\{
\barray{ll} \{ m - \h \} & \mbox{if $m \geqs \h$}~, \\
\emptyset & \mbox{if $m < \h$}~. \farray \right. \feq
\end{prop}
\textbf{Proof.} i) Let $t \in [t_{m},t_{m+1})$.
For $s \in [t_k, t_{k+1})$ one has $t - s > t_{m} - t_{k+1} \geqs \te$, implying
$\um(t-s) = A e^{-B (t-s)}$; so,
\beq \int_{t_k}^{t_{k+1}} \! \! \! \! \! \! d s~\um(t-s)
\left({s - t_k \over t_{k+1} - t_k}\right)^2 =
A e^{- B t} \int_{t_k}^{t_{k+1}} \! \! \! \! \! \! d s~e^{B s}
\left({s - t_k \over t_{k+1} - t_k}\right)^2 = A e^{- B t} H_k \feq
$$ \leqs A e^{- B t_m} H_k = H_{m k}~\mbox{as in \rref{bb6}}~. $$
In conclusion, defining $H_{m k}$ as in \rref{bb6}
we fulfill the first inequality in \rref{b6} (note that $k \leqs m-1$ due to
$t_{k+1} \leqs t_m$). Similarly, the other inequalities
\rref{b6} are fulfilled with $I_{m k}$, $N_{m k}$ as in \rref{bb6}. \parn
Finally, inserting Eq. \rref{bb6} into Eq. \rref{b7} for $\Phi_{m k}$ we obtain
the thesis \rref{bb7}. \parn
ii) Let us rephrase Eq. \rref{recur} in the case under examination. To this purpose,
we reexpress the sum therein writing \beq
\sum_{k=0}^m = \sum_{\{ k | t_{m} - \te < t_{k+1} \leqs t_{m+1} \}} +
\sum_{ \{ k | t_{k+1} \leqs t_{m} - \te \} }~, \feq
and then use Eq. \rref{bb7} for the summands with $t_{k+1} \leqs t_m - \vartheta$; in this way, Eq.
\rref{recur} becomes
$$ \EE_{m} +
K \ski \sum_{\{ k | t_{m} - \te < t_{k+1} \leqs t_{m+1} \}} \ski \Phi_{m k}(\RR_k, \RR_{k+1})
+ K A e^{- B t_m}~~\ski \sum_{ \{ k | t_{k+1} \leqs t_{m} - \te \} } \ski~~ \Phi_{k}(\RR_k, \RR_{k+1}) $$
\beq \leqs \min(\RR_m, \RR_{m + 1}) \qquad \mbox{for $m \in \{0,...,\EM-1\}$}~. \label{recurr} \feq
Let us consider any sequence $(\RR_{m})_{m \in \{0,...,\EM\}}$ of nonnegative reals. If
$(\RR_{m})$ fulfills \rref{recurr}, define
\beq \SS_{m} := \sum_{ \{ k | t_{k+1} \leqs t_{m} - \te \} } \ski~~ \Phi_{k}(\RR_k, \RR_{k+1})
\qquad \mbox{for $m \in \{0,...,\EM - 1\}$}~; \label{defesel} \feq
then \rref{recur2} \rref{recur1} follow immediately (in fact,
with $\leqs$ replaced by $=$ in \rref{recur2}). \parn
Conversely, suppose there is a
sequence of reals $(\SS_{m})_{m \in \{0,...,\EM-1\}}$ fulfilling  Eq.s \rref{recur2} and
\rref{recur1} with $(\RR_{m})$; then
$\sum_{ \{ k | t_{k+1} \leqs t_{m} - \te \} } \Phi_{k}(\RR_k, \RR_{k+1}) \leqs
\SS_m$, and it is easy to infer Eq. \rref{recurr} for $(\RR_{m})$. \parn
iii) Obvious. \fine
\begin{rema}
\label{remns}
\textbf{Remark. (i)} Eq. \rref{recur2} does not prescribe $\SS_0$. It is convenient to choose
$\SS_0 := 0$; with this position Eq. \rref{recur1} with $m= 0$ is a problem for two unknowns $(\RR_0, \RR_1)$,
for which we could repeat the comments of Remark \ref{re} (ii). After these initial steps, we can use Eq.s
\rref{recur2} \rref{recur1} as recursion relations to obtain $\SS_1, \RR_2, \SS_2, \RR_3$, and so on.
\parn
\textbf{(ii)}
As anticipated, Proposition \ref{special} can be
applied to the NS equations, in the framework
of Section \ref{basic}. Eq. \rref{eqww} of the cited section gives the semigroup estimator
$\um(t) := e^{t}/\sqrt{2 e t}$ for $t \leqs 1/4$, and
$\um(t) := \sqrt{2} e^{-t}$ for all $t > 1/4$. So, the conditions of the previous
proposition are fulfilled with $\te = 1/4$, $A = \sqrt{2}$, $B = 1$. These values
must be sustituted into Eq. \rref{bb6} for $H_{m k}$, $I_{m k}$ and $N_{m k}$ giving,
for example,
\beq N_{m k} := \sqrt{2} e^{- t_{m}} N_{k}~,~~ N_k := e^{t_{k+1}} - e^{t_k}
\qquad \mbox{for $t_{k+1} \leqs t_m - 1/4$}~. \feq
Explicit expressions could be derived as well for $H_{m k}$, $I_{m k}$ and $N_{m k}$ when $t_m - 1/4 < t_{k+1}
\leqs t_m$, using elementary bounds on $\um$ derived from the expression
\rref{eqww}. However, this requires a tedious analysis of a number of cases, since
the parameter $t-s$ in Eq. \rref{b6} can be smaller or greater that $1/4$;
details will be given elsewhere, when we will treat systematically the approach outlined
in this Appendix.
\end{rema}
\section{Appendix. Proof of Lemmas \ref{lemsi} and \ref{lemci}.}
\label{appemu}
We begin with two auxiliary Lemmas.
\begin{prop}
\label{somme}
\textbf{Lemma.} Let us consider two radii $\rho, \rho_1$ such that $2 \sqrt{d} \leqs \rho < \rho_1 \leqs + \infty$, and
a nonincreasing function $\chi \in C([\rho - 2 \sqrt{d}, \rho_1), [0,+\infty)$). Then,
\beq \sum_{h \in \Zd, \rho \leqs | h | < \rho_1} \chi(|h|) \leqs
{2 \pi^{d/2} \over \Gamma(d/2)} \int_{\rho - 2 \sqrt{d}}^{\rho_1} dt \, (t + \sqrt{d})^{d-1} \chi(t)~,
\label{asinthe} \feq
\end{prop}
\textbf{Proof.} Let us introduce the cubes
\beq h + [0,1]^d  \qquad (h \in \Zd) \feq
and the annulus
\beq A(\rho - \sqrt{d}, \rho_1 + \sqrt{d}) :=
\{ q \in \reali^d~|~\rho - \sqrt{d} \leqs |q | <  \rho_1 + \sqrt{d}~\}~; \feq
we claim that
\beq \cup_{h \in \Zd\!, \, \rho \leqs |h | < \rho_1} (h + [0,1]^d) \subset A(\rho - \sqrt{d}, \rho_1 + \sqrt{d})~.
\label{claimq} \feq
In fact, $q \in h + [0,1]^d$ implies $|h| - \sqrt{d} \leqs |q| \leqs |h | + \sqrt{d}$; now, if
$\rho \leqs |h | < \rho_1$ we conclude $\rho - \sqrt{d} \leqs | q | < \rho_1 + \sqrt{d}$. \parn
The inclusion \rref{claimq} implies
\beq \sum_{h \in \Zd, \rho \leqs |h | < \rho_1} \int_{h + [0,1]^d} \!\!\!\! dq \,\chi(|q| - \sqrt{d})
\leqs \int_{A(\rho- \sqrt{d}, \rho_1 + \sqrt{d})} \!\!\!\! dq \, \chi(|q| - \sqrt{d})~. \label{lhsq} \feq
On the other hand, for $h$ as in the above
sum and $q \in h + [0,1]^d$, we have $|q| - \sqrt{d} \leqs | h |$, whence
$\chi(|q| - \sqrt{d}) \geqs \chi(|h|)$; this implies
\beq \int_{h + [0,1]^d} \!\!\!\! dq \, \chi(|q| - \sqrt{d}) \geqs \chi(|h|)
\int_{h + [0,1]^d} \!\!\!\! dq = \chi(|h|)~. \feq
From here and from \rref{lhsq} we obtain
\beq \sum_{h \in \Zd, \rho \leqs |h | < \rho_1} \chi(|h|) \leqs
\int_{A(\rho- \sqrt{d}, \rho_1 + \sqrt{d})} \!\!\!\! dq \, \chi(|q| - \sqrt{d})~. \label{lhsqq} \feq
The right hand side of \rref{lhsqq} can be expressed in terms of the one-dimensional variable $r = |q|$;
as well known, $d q = (2 \pi^{d/2}/\Gamma(d/2)) r^{d-1} d r$, so
\beq \sum_{h \in \Zd, \rho \leqs | h | < \rho_1} \chi(|h|) \leqs {2 \pi^{d/2}
\over \Gamma(d/2)} \int_{\rho - \sqrt{d}}^{\rho_1 + \sqrt{d}}  \!\!\!\! dr \,  r^{d-1} \chi(r - \sqrt{d})~;
\label{asinthes} \feq
now, a change of variables $t  = r - \sqrt{d}$ gives the thesis \rref{asinthe}. \fine
To go on, let us recall the convention \rref{ze}
$$ \Ze := \Zd ~\mbox{or}~\Zd_0~. $$
\begin{prop}
\label{lemsiz}
\textbf{Lemma.} (i) Generalizing \rref{desgm}, let
\beq \Sim_\nu(\lambda) := {1 \over (2 \pi)^d} \sum_{h \in \Ze, | h | < \lambda}
{1 \over (1 + | h |^2)^\nu} \qquad
\mbox{for $\nu, \lambda \in (0,+\infty)$}~. \label{sim} \feq
Then,
\beq \Sim_{\nu}(\lambda) \leqs S_{\nu} + {(1 + d)^\nu \over 2^{d-1} \pi^{d/2} \Gamma(d/2)}~
F_{\nu}(\lambda) \quad \mbox{for $\nu > 0$, $\lambda > 2 \sqrt{d}$}~, \label{essim} \feq
$$ S_\nu := {1 \over (2 \pi)^d} \sum_{h \in \Ze, | h | < 2 \sqrt{d}}
{1 \over (1 + | h |^2)^\nu}, $$
$$ F_{\nu}(\lambda) :=
\left\{ \barray{ll} \dd{{1 \over 2 \nu - d} \Big({1 \over \sqrt{d}^{2 \nu - d}}
- {1 \over (\lambda + \sqrt{d})^{2 \nu - d}}\Big)} & \mbox{if $\nu \neq d/2$}~, \\
\dd{\log({\lambda + \sqrt{d} \over \sqrt{d}})}  & \mbox{if $\nu = d/2$}~. \farray \right. $$
For fixed $\nu$ and $\lambda \vain + \infty$, this implies
\beq \Sim_{\nu}(\lambda) = \left\{ \barray{lll} O(1) & \mbox{if $\nu > d/2$}~, \\
O(\log \lambda)  & \mbox{if $\nu = d/2$}~, \\
O(\lambda^{d - 2 \nu}) & \mbox{if $0 < \nu < d/2$}~. \farray \right. \label{esim} \feq
(ii) Let
\beq \Sip_\nu(\lambda) := {1 \over (2 \pi)^d} \sum_{h \in \Ze, | h | \geqs \lambda}
{1 \over (1 + | h |^2)^\nu} \qquad \mbox{for $\nu > d/2$, $\lambda > 0$}~. \label{sip} \feq
Then
\beq \Sip_{\nu}(\lambda) \leqs \delta \Sigm_{\nu}(\lambda)
\quad \mbox{for $\nu > d/2$, $\lambda > 2 \sqrt{d}$}~, \label{essip} \feq
where (generalizing \rref{reph})
\beq \delta \Sigm_{\nu}(\lambda) :=
{(1 + d)^\nu \over 2^{d-1} \pi^{d/2} \Gamma(d/2) (2 \nu - d)~
(\lambda - \sqrt{d})^{2 \nu - d}}~. \label{essipp} \feq
\end{prop}
\textbf{Proof.}
(i) Let $\nu > 0$, $\lambda > 2 \sqrt{d}$. Dividing in two parts the sum defining $\Sim_{\nu}(\lambda)$, we get
\beq  \Sim_\nu(\lambda) = S_\nu + {1 \over (2 \pi)^d} \sum_{h \in \Ze, 2 \sqrt{d} \leqs | h | < \lambda}
{1 \over  (1 + | h |^2)^\nu}~. \feq
Now, we bind the sum using \rref{asinthe} with $\chi(t) := 1/(1 + t^2)^{\nu}$, $\rho := 2 \sqrt{d}$
and $\rho_1 := \lambda$, yielding
\beq \Sim_\nu(\lambda) \leqs S_\nu + {1 \over 2^{d-1} \pi^{d/2} \Gamma(d/2)}
\int_{0}^{\lambda} dt \, {(t + \sqrt{d})^{d-1} \over (1 + t^2)^{\nu}}~. \label{thisgives} \feq
On the other hand, one establishes by elementary means the inequality
\beq {1 \over 1 + t^2} \leqs {1 + d \over (t + \sqrt{d})^2} \qquad \mbox{for $t \in [0,+\infty)$}~ \label{elemen} \feq
(holding as an equality when $t = 1/\sqrt{d}$). Inserting this into \rref{thisgives}, we get
\beq \Sim_\nu(\lambda) \leqs S_\nu + {(1 + d)^\nu \over 2^{d-1} \pi^{d/2} \Gamma(d/2)}
\int_{0}^{\lambda} dt \, {1 \over (t + \sqrt{d})^{2 \nu - d + 1} }~. \feq
The last integral equals $F_{\nu}(\lambda)$, so
\rref{essim} is proved. Having this result, the statement \rref{esim} on the limit $\lambda \vain + \infty$
is obvious. \parn
(ii) Let $\nu > d/2$, $\lambda > 2 \sqrt{d}$. To bind
$\Sip_n(\lambda)$ we use
\rref{asinthe} with $\chi(t) := 1/(1 + t^2)^{\nu}$, $\rho := \lambda$
and $\rho_1 := +\infty$, and subsequently employ the inequality \rref{elemen}. This gives \parn
{\vbox{
\beq  \Sip_\nu(\lambda) \leqs {1 \over 2^{d-1} \pi^{d/2} \Gamma(d/2)}
\int_{\lambda - 2 \sqrt{d}}^{+\infty} dt \, {(t + \sqrt{d})^{d-1} \over (1 + t^2)^{\nu}} \feq
$$ \leqs {(1 + d)^\nu \over 2^{d-1} \pi^{d/2} \Gamma(d/2)}
\int_{\lambda - 2 \sqrt{d}}^{+\infty} dt \, {1 \over (t + \sqrt{d})^{2 \nu - d + 1} }~, $$}}
and computing the last integral we justify Eqs. \rref{essip} \rref{essipp}. \fine
From here to the end of the Appendix we fix a real number $n$, fulfilling the relation \rref{nd}
$$ n > {d \over 2}~; $$
here are the proofs of Lemmas \ref{lemsi}, \ref{lemci}. \salto
\textbf{Proof of Lemma \ref{lemsi}.} We take any  $\lambda \geqs 2 \sqrt{d}$; with the notations
of the previous Lemma, we have
\beq \Sigma_n = \Sim_n(\lambda) + \Sip_n(\lambda)~. \feq
Both terms in the right hand side are finite; the
term $\Sip_{n}(\lambda)$ has the upper bounds \rref{essip} and the obvious
lower bound $\Sip_{n}(\lambda) > 0$. From these facts we infer the finiteness of $\Sigma_n$, and the bounds
\rref{decosig} for it. \fine
\salto
\textbf{Proof of Lemma \ref{lemci}.} We consider any cutoff function $\Lambda : \Ze
\vain [2 \sqrt{d}, + \infty)$. For $k \in \Ze$, we introduce the decomposition
\beq \KK_n(k) = \Ki_n(k) + \Delta \Ki_n(k)~; \label{dexx} \feq
here $\Ki_n (k)$ is defined by \rref{mndig2}, and
\beq \Delta \Ki_n (k) = {(1 + | k |^2)^{n-1} \over (2 \pi)^d}
\sum_{h \in \Ze, |h | > \Lambda(k)}
{| k - h |^2 \over (1 + |h|^2)^n (1 + | k - h |^2)^n}~. \feq
For the term $\Ki_n(k)$, we furtherly introduce a decomposition
\beq \Ki_n(k) = \Ki'_n(k) + \Ki''_n(k)~, \feq
\beq \Ki'_n (k) := {(1 + | k |^2)^{n-1} \over (2 \pi)^d}
\sum_{h \in \Ze, |h | < | k |/2}
{| k - h |^2 \over (1 + |h|^2)^n (1 + | k - h |^2)^n}~, \feq
\beq \Ki''_n (k) := {(1 + | k |^2)^{n-1} \over (2 \pi)^d}
\sum_{h \in \Ze, | k |/2 \leqs | h | < \Lambda(k)}
{| k - h |^2 \over (1 + |h|^2)^n (1 + | k - h |^2)^n} \feq
(the sum defining $\Ki''(k)$ is meant to be zero if $\Lambda(k) \leqs |k|/2$). In the sequel
we analyse separately $\Delta \Ki_n$, $\Ki'_n$ and $\Ki''_n$; we will frequently use
the inequalities
\beq {| k - h |^2 \over (1 + | k - h |^2)^n} \leqs {1 \over (1 + | k - h |^2)^{n-1}}
\qquad \mbox{for all $k, h \in \Ze$}~; \label{freq} \feq
\beq {1 + z \over 1 + \eta z} \leqs \max(1, {1 \over \eta})
\qquad \mbox{for all $\eta > 0$, $z \geqs 0$}~. \label{ux} \feq
\textsl{Step 1. For all $k \in \Ze$, one has}
\beq 0 < \Delta \Ki_{n}(k) \leqs \delta \Ki_{n}(k)~; \label{thes} \feq
\textsl{here, as in \rref{deci}},
$$ \delta \Ki_{n}(k) :=
{(1 + d)^n \lan(|k|)^{n-1} \over 2^{d-1} \pi^{d/2} \Gamma(d/2) (2 n - d)~
(\Lambda(k) - \sqrt{d})^{2 n - d}}~, $$
\textsl{$\lan$ being defined by \rref{defom}. If $\Lambda$ fulfills \rref{cutoff}, the above upper bound
is such that}
\beq \delta \Ki_{n}(k) = O({1 \over | k |^{2 n - d}}) \vain 0
\quad \mbox{for $k \vain \infty$}~. \label{impc} \feq
The inequality $\Delta \Ki_n(k) > 0$ is obvious. To prove the rest we start from Eq.
\rref{freq}, implying
\beq \Delta \Ki_{n}(k) \leqs {(1 + | k |^2)^{n-1} \over (2 \pi)^d}
\sum_{h \in \Ze, | h | \geqs \Lambda(k)}
{1 \over (1 + |h|^2)^n (1 + | k - h |^2)^{n-1}}~; \label{muf} \feq
setting
\beq \mu(k) := \inf_{h \in \Ze, |h | \geqs \Lambda(k)} |k - h |~, \feq
we infer from \rref{muf} that
\beq \Delta \Ki_{n}(k)
\leqs {1 \over (2 \pi)^d} \left({1 + | k |^2 \over 1 + \mu^2(k)} \right)^{n-1} \!\!\!\!\!
\sum_{h \in \Ze, | h | \geqs \Lambda(k)} {1 \over (1 + |h|^2)^n}  \label{sint} \feq
$$ = \left({1 + | k |^2 \over 1 + \mu(k)^2}  \right)^{n-1} \Sip_{n}(\Lambda(k))~,$$
where $\Sip_n$ is defined following Eq. \rref{sip}. To go on we claim that, for all $k \in \Ze$,
\beq \mu(k) \geqs
\left\{ \barray{ll} 0 & \mbox{if $\Lambda(k) < | k |$}~, \\
\dd{\Lambda(k) - | k |} & \mbox{if $\Lambda(k) \geqs | k |$}~. \farray \right. \label{cl} \feq
The above inequality is trivial if $\Lambda(k) < | k |$; if $\Lambda(k) \geqs | k |$, it follows
noting that $| h | \geqs \Lambda(k)$ implies $| k - h | \geqs | h | - | k | \geqs \Lambda(k) - | k |$. \parn
Having proved \rref{cl}, we insert it into \rref{sint}; this gives the inequality
\beq \Delta \Ki_{n}(k) \leqs \lan(k)^{n-1} \Sip_{n}(\Lambda(k)) \feq
and substituting therein the upper bound \rref{essip} for $\Sip_n$ we get the upper bound in \rref{thes}. \parn
To go on, suppose $\Lambda$ fulfills \rref{cutoff}; then
\beq \Lambda(k) \geqs \alpha | k | \geqs |k|,~\lan(k) \leqs
{1 + | k |^2 \over 1 + (\alpha - 1)^2 |k|^2} \leqs \max(1, {1 \over (\alpha-1)^2})
~\mbox{for $|k | \geqs \chi$} \feq
(the last inequality follows from \rref{ux}, with $z = | k |^2$ and $\eta = (\alpha - 1)^2$). So,
\beq {1 \over \Lambda(k) - \sqrt{d}} = O({1 \over | k |})~,~~\lan(k) = O(1)
\qquad \mbox{for $k \vain \infty$}~; \feq
inserting this in the definition of $\delta \Ki_n$, we infer Eq. \rref{impc}. \parn
\textsl{Step 2. With $\Sigma_n$ as in \rref{desig}, one has}
\beq \Ki'_n (k) \vain \Sigma_n \qquad \mbox{for $k \vain \infty$}~. \feq
To prove this, we write
\beq \Ki'_n (k) = \sum_{h \in \Ze} c_{n k}(h)~, \feq
$$ c_{n k}(h) := {\theta(|k|/2 - | h |)
| k - h |^2 (1 + | k |^2)^{n-1} \over (2 \pi)^d (1 + |h|^2)^n (1 + | k - h |^2)^n}~,
\qquad \theta(z) := \left\{ \barray{ll} \dd{1} & \mbox{if $z \in (0, + \infty)$}~, \\
\dd{0} & \mbox{if $z \in (-\infty, 0]$}~. \farray \right. $$
For any fixed $h \in \Ze$, we have
\beq c_{n k}(h) \vain {1 \over (2 \pi)^d}  {1 \over (1 + | h |^2)^n} \qquad \mbox{for $k \vain + \infty$}~; \feq
this implies
\beq \Ki'_n(k) \vain {1 \over (2 \pi)^d} \sum_{h \in \Ze} {1 \over (1 + | h |^2)^n} = \Sigma_n~, \feq
if the limit $k \vain \infty$ can be exchanged with the sum over $h$. By the Lebesgue theorem
of dominated convergence, the exchange is possible if $c_{n k}(h)$ is bounded from above
by a summable function of $h$, uniformly in $k$; indeed this occurs, since
\beq c_{n k}(h) \leqs {4^{n-1} \over (2 \pi)^d} \,
{1 \over (1 + | h |^2)^n} \qquad \mbox{for all $h, k \in \Ze$}~. \label{since} \feq
Let us prove \rref{since}. The thesis is obvious if $|h| \geqs |k|/2$, since in this case $c_{n d}(k) =0$;
hereafter we assume $|h| < |k|/2$. First of all, we note that \parn
{\vbox{
\beq c_{n k}(h) = {| k - h |^2 (1 + | k |^2)^{n-1} \over (2 \pi)^d (1 + |h|^2)^n (1 + | k - h |^2)^n} \feq
$$ \leqs {(1 + | k |^2)^{n-1} \over (2 \pi)^d (1 + |h|^2)^n (1 + | k - h |^2)^{n-1}}~; $$}}
secondly, from $|h| < |k|/2$ we infer $| k - h | \geqs | k | - | h | \geqs |k|/2$, whence
\beq c_{n k}(h) \leqs
{1 \over (2 \pi)^d} \left({1 + | k |^2 \over 1 + |k|^2/4}\right)^{n-1} {1 \over (1 + |h|^2)^n}
\leqs  {4^{n-1} \over (2 \pi)^d} \,
{1 \over (1 + | h |^2)^n}~ \label{apto} \feq
(the last inequality follows from \rref{ux}, with $\eta = 1/4$ and $z = | k |^2$).  \parn
\textsl{Step 3. Suppose $\Lambda$ fulfills \rref{cutoff}; then, for $k \vain \infty$},
\beq \Ki''_{n}(k) = \left\{ \barray{lll} O(\dd{1 \over |k|^2}) & \mbox{if $n > d/2 + 1$}~, \\
O(\dd{\log | k | \over |k|^2})  & \mbox{if $n = d/2 + 1$}~, \\
O(\dd{1 \over | k|^{2 n - d}}) & \mbox{if $d/2 < n < d/2 +1$}~ \farray \right.
\quad \vain 0~. \label{esimm} \feq
To prove this, let us take any $k \in \Ze$ such that $|k | \geqs \chi$. First of all, from
\rref{freq} and $\Lambda(k) \leqs \beta | k |$} we infer
\beq \Ki''_n (k) \leqs {(1 + | k |^2)^{n-1} \over (2 \pi)^d}
\sum_{h \in \Ze, | k |/2 \leqs | h | < \Lambda(k)}
{1 \over (1 + |h|^2)^n (1 + | k - h |^2)^{n-1}} \feq
$$ \leqs {(1 + | k |^2)^{n-1} \over (2 \pi)^d} \sum_{h \in \Ze, | k |/2 \leqs | h | < \beta | k | }
{1 \over (1 + |h|^2)^n (1 + | k - h |^2)^{n-1}}~. $$
For $| h | \geqs | k |/2$ we have $1/(1 + | h |^2) \leqs 1/(1 + |k|^2/4)$, whence
\beq \Ki''_n (k) \leqs {1 \over (2 \pi)^d} {(1 + | k |^2)^{n-1} \over (1 + | k |^2/4)^n}
\sum_{h \in \Ze, | k |/2 \leqs | h | < \beta | k | }
{1 \over (1 + | k - h |^2)^{n-1}}~\feq
$$ = {1 \over (2 \pi)^d} \left({1 + | k |^2 \over 1 + | k |^2/4} \right)^{n-1}
{1 \over 1 + | k |^2/4} \sum_{h \in \Ze, | k |/2 \leqs | h | < \beta | k | }
{1 \over (1 + | k - h |^2)^{n-1}} $$
$$ \leqs {4^{n-1} \over (2 \pi)^d} \, {1 \over 1 + | k |^2/4} \sum_{h \in \Ze, | k |/2 \leqs | h | < \beta | k | }
{1 \over (1 + | k - h |^2)^{n-1}} $$
(the last passage uses again \rref{ux}, with $\eta = 1/4$ and $z = | k |^2$).
Now, a change of variable $h = q - k$ in the last sum gives
\beq \Ki''_n (k) \leqs {4^{n-1} \over (2 \pi)^d} {1 \over 1 + | k |^2/4}  \sum_{q \in \Ze, | k |/2 \leqs | q - k | < \beta | k | }
{1 \over (1 + | q |^2)^{n-1}}~. \feq
To go on we note that, for all $q \in \Ze$,
\beq | q - k | < \beta | k |~~ \Rightarrow~~| q | \leqs |q - k | + | k | < (\beta + 1) | k |~. \feq
Thus,
\beq \Ki''_n (k) \leqs {4^{n-1} \over (2 \pi)^d} \, {1 \over 1 + | k |^2/4}  \sum_{q \in \Ze, |q | < (\beta + 1) | k | }
{1 \over (1 + | q |^2)^{n-1}} \feq
$$ = {4^{n-1} \over 1 + | k |^2/4}~\Sim_{n-1}((\beta + 1) | k |)~, $$
where $\Sim_{n-1}$ is defined following Eq. \rref{sim}. The $k \vain \infty$ behaviour
of $\Sim_{n-1}((\beta + 1) | k |)$ is inferred from Eq. \rref{esim}; this function is multiplied
by $4^{n-1}/(1 + | k |^2/4) = O(1/|k|^2)$, and the conclusion is Eq. \rref{esimm}. \parn
\textsl{Step 4. Conclusion of the proof.} For any cutoff $\Lambda$,
the decomposition \rref{dexx} and the bounds \rref{thes} give the inequalities
\rref{kkn} for $\KK_n(k)$, also implying the finiteness of this quantity for
arbitrary $k$. \parn
From now on $\Lambda$ is supposed to fulfill \rref{cutoff}, and
we consider the limit $k \vain \infty$. Then,
the decomposition $\Ki_n = \Ki'_n + \Ki''_n$ and Steps 2, 3 give
$$ \Ki_n(k) \vain \Sigma_n~, $$
as claimed in \rref{nona}; the other statement in \rref{nona}, namely,
$\delta \Ki_n(k) = O(1/|k|^{2 n - d}) \vain 0$, is known from \rref{impc}.
The decomposition $\KK_n = \Ki_n + \Delta \Ki_n$, and
Eqs. \rref{thes} \rref{nona} also give
\beq \KK_n(k) \vain \Sigma_n~. \feq
Summing up, all statements in items (i)(ii) of this lemma are proved. \fine
\section{Appendix. Proof of Propositions \ref{multi} and \ref{mult}.}
\label{appemult}
\textbf{Proof of Proposition \ref{multi}.}  We consider two vector fields
$v, w \in \HO{n}$ with $n > d/2$, and proceed in two steps. \parn
\textsl{Step 1. Proof of Eq. \rref{fuco}.} We have
\beq v^r = \sum_{h \in \Zd} v^r_h e_h~, \qquad \partial_r w^s =
i \sum_{\ell \in \Zd} \ell_r w^s_\ell e_\ell~;\feq
this implies $v \sc \, \partial w^s = v^r \partial_r w^s = i \sum_{h, \ell \in \Zd}
(v^r_h \ell_r w^s_\ell) (e_h e_\ell)$ or, in vector form,
\beq v \sc \, \partial w = i \sum_{h, \ell \in \Zd}
(v_h \sc \, \ell) w_\ell (e_h e_\ell)~. \feq
On the other hand $e_h e_{\ell} = e_{h + \ell}/(2 \pi)^{d/2}$, so
\beq v \sc \, \partial w = {i \over (2 \pi)^{d/2}} \sum_{k \in \Zd}
\Big(\sum_{h, \ell \in \Zd\!\!,\,
h + \ell = k} (v_h \sc \, \ell) w_\ell \Big) e_k  \feq
$$ = {i \over (2 \pi)^{d/2}} \sum_{k \in \Zd}
\Big(\sum_{h \in \Zd} [v_{h} \sc \, (k - h)] w_{k - h} \Big) e_k~. $$
The term multiplying $e_k$ in the last expression is the
Fourier coefficient $(v \sc \, \partial w)_k$; so, Eq. \rref{fuco} is proved. \parn
\textsl{Step 2. Proof that $v \sc \partial w$ is in $\HO{n-1}$ and
fulfills \rref{esto}.} For any $k \in \Zd$, Eq. \rref{fuco}
implies
\beq |(v \sc \, \partial w)_k| \leqs {1 \over (2 \pi)^{d/2}} \sum_{h \in \Zd} |v_{h}| | k - h |  | w_{k - h} |  \feq
$$ = {1 \over (2 \pi)^{d/2}} \sum_{h \in \Zd}  {| k - h | \over \sqrt{1 + |h|^2}^{n} \sqrt{1 + | k - h |^2}^n}~
\sqrt{1 + |h|^2}^n | v_{h} | \sqrt{1 + | k - h |^2}^n~| w_{k - h} |~. $$
Now, H\"older's inequality $| \sum_h~ a_h b_h |^2 \leqs \Big(\sum_h | a_h |^2\Big)
\Big(\sum_h~| b_h |^2 \Big)$ gives
\beq | (v \sc \, \partial w)_k |^2 \leqs c_k p_k~, \label{dains} \feq
$$ c_k := {1 \over (2 \pi)^d} \sum_{h \in \Zd} {| k - h |^2 \over (1 + | h |^2)^n (1 + | k - h |^2)^n }~, $$
$$ p_k := \sum_{h \in \Zd} (1 + |h|^2)^n |  v_{h} |^2
(1 + | k - h |^2)^n | w_{k - h} |^2~. $$
The last inequality implies
\beq \sum_{k \in \Zd} (1 + | k |^2)^{n-1} | (v \sc \, \partial w)_k |^2 \leqs
\sum_{k \in \Zd} (1 + | k |^2)^{n-1} c_k p_k  \label{dains0} \feq
$$ \leqs \left( \sup_{k \in \Zd} (1 + | k |^2)^{n-1} c_k \right) \sum_{k \in \Zd} p_k ~. $$
On the other hand,
\beq (1 + | k |^2)^{n-1} c_k = \KK_{n}(k)~, \qquad \sup_{k \in \Zd} (1 + | k |^2)^{n-1} c_k \leqs K^2_{n}
\label{dains1} \feq
with $\KK_{n}$, $K_n$ as in Eqs. \rref{mmnnddii} \rref{mndi};
the finiteness of $\KK_n(k)$ for any $k$, and of its sup over $k$ are known from
Lemma \ref{lemci}. To go on, we note that
\beq \sum_{k \in \Zd} p_k = \Big(\sum_{h \in \Zd} (1 + |h |^2)^n | v_{h} |^2 \Big)
\Big(\sum_{h \in \Zd} (1 + |h |^2)^n | w_h |^2 \Big)= \| v \|^2_n \, \| w \|^2_n~.
\label{dains2} \feq
Inserting Eqs. \rref{dains1}, \rref{dains2} into \rref{dains0},
we see that
$\sum_{k \in \Zd} (1 + | k |^2)^{n-1} | (v \sc \, \partial w)_k |^2 < + \infty$, implying $v \sc \, \partial w
\in \HO{n-1}$. Furthermore,
\beq \| v \sc \, \partial w \|^2_{n-1} = \sum_{k \in \Zd} (1 + | k |^2)^{n-1} | (v \sc \, \partial w)_k |^2
\leqs \, K^2_{n} \, \| v \|^2_{n} \, \| w \|^2_n~, \feq yielding Eq. \rref{esto}.
\fine
\textbf{Proof of Proposition \ref{mult}.} This is a
simple variant of the proof given for Proposition
\ref{multi}. One takes into account the following facts: if $f, g \in \HM{n}$, then their zero order
Fourier coefficients are $f_0 = 0$, $g_0 = 0$; furthermore, $(f \sc \, \partial g)_0 = 0$, since this function
has mean zero by Lemma \ref{lemean}. \parn
\fine
\section{Appendix. Derivation of the NS equations \rref{nsl}.}
\label{appens}
\textbf{The NS equations in physical and adimensional units.}
\salto
In the space $\len$ of (oriented) lenghts we fix some positive lenght $\LS$, determining
the size of the system in consideration. The "space domain" of the system is modelled
as a torus $\len^d/(2 \pi \LS)^d$; we write $\XS$ for any point in this domain, and $\TS$
for the "physical" time. The NS equations are
\beq \DS \big({\boma{\dot \VS}}(\TS) + \VS(\TS) \sc \, \partial \VS(\TS)\big) = - \partial \PR(\TS) +
\Lam \Delta \VS(\TS) + \KS(\XS,\TS)
\label{eqphy} \feq
where
$\VS(\TS) : \XS \mapsto \VS(\XS,\TS)$, $\PR(\TS) : \XS \mapsto \PR(\XS, \TS)$ are the velocity
and the pressure fields, while $\KS(\TS) : \XS \mapsto \KS(\XS,\TS)$ is the density of
external forces; in the above, $\boma{\dot{~}}$ indicates the \textsl{partial} derivative
with respect to time $\TS$. The coefficients $\DS > 0$, $\Lam > 0$ are the density and viscosity of
the fluid, respectively. The velocity field $\VS(\TS)$ is required to be divergence free, to
fulfill the condition of incompressibility. \parn
One passes to the adimensional form introducing three functions
$\nu(x,t)$, $\pi(x,t)$, $\kappa(x,t)$ ($x \in \Td$, $t \in [0,T) \subset \reali)$ via the equations
\beq \VS(\XS,\TS) = {\Lam \over \DS \LS} \nu(x,t)~, \quad \PR(\XS,\TS) = {\Lam^2 \over \DS \LS^2} \pi(x,t)~,
\qquad \KS(\XS,\TS)= {\Lam^2 \over \DS \LS^3} \kappa(x,t) \feq
$$ \mbox{for}~~ x = {\XS \over \LS},~~t = {\Lam \TS \over \DS \LS^2}~. $$
With these positions, Eq. \rref{eqphy} is equivalent to
$$ {\dot \nu}(t) + \nu(t) \sc \, \partial \nu(t) = - \partial \pi(t) + \Delta \nu(t) + \kappa(t)~. $$
Hereafter we formalise this adimensional version, specifying the necessary functional spaces.
First of all, we suppose
\beq \kappa \in C^{0,1}([0,+\infty), \HO{n-1})~; \feq
secondly, we stipulate the following.
\begin{prop}
\textbf{Definition.}  The \textsl{incompressible NS Cauchy problem with initial
datum} $v_0 \in \Hs{n+1}$, in the \textsl{pressure formulation}, is the following.
$$ \mbox{Find $\nu \in C([0,T),\Hs{n+1}) \cap C^1([0,T), \Hs{n-1})$, $\pi \in C([0,T), H^{n})$
such that} $$
\beq {\dot \nu}(t) + \nu(t) \sc \, \partial \nu(t) = - \partial \pi(t) + \Delta \nu(t) + \kappa(t)
\quad \mbox{for $t \in [0,T)$}~, \qquad \nu(0) = v_0~ \label{nsp} \feq
(for some $T \in (0,+\infty]$).
\end{prop}
We note the following: the requirements $\nu \in C([0,T),\Hs{n+1})$ and $\pi \in C([0,T), H^{n})$ are
sufficient for the right hand side of the above differential equation to be in
$C([0,T), \Hs{n-1})$. \parn
\salto
\textbf{The equivalence between the pressure formulation \rref{nsp} and the
Leray formulation \rref{nsl} of the Cauchy problem.}
\salto
Let us introduce the function
\beq \eta \in C^{0,1}([0,+\infty), \Hs{n-1})~, \qquad t \mapsto \eta(t) := \LP \kappa(t)~. \feq
For convenience, we report here the Cauchy problem in the form \rref{nsl}:
$$ \mbox{Find $\nu \in C([0,T),\Hs{n+1}) \cap C^1([0,T), \Hs{n-1})$,
such that} $$
$$ {\dot \nu}(t) = \Delta \nu(t)- \LP \big( \nu(t) \sc \, \partial \nu(t) \big) + \eta(t)
\quad \mbox{for $t \in [0,T)$}~, \qquad \nu(0) = v_0~ $$
(for some $T \in (0,+\infty]$).
The above mentioned equivalence can be stated as follows.
\begin{prop}
\textbf{Proposition.} A function $\nu$ of domain $[0,T)$ fulfills \rref{nsl} if and
only if there is a function $\pi$ with the same domain, such that $(\nu,\pi)$ fulfills \rref{nsp}.
\end{prop}
\textbf{Proof.} Suppose a pair $(\nu, \pi)$ fulfills \rref{nsp}, and apply the projector $\LP$
to both sides of the differential equation. We have $\LP \nu(t) = \nu(t)$, and $\LP$ commutes
with both the time derivative $\dot{~}$ and the Laplacian $\Delta$; finally, $\LP \partial \pi(t) = 0$.
These facts yield Eq. \rref{nsl}. Conversely, suppose a function $\nu$ fulfills \rref{nsl} and define
\beq \gamma \in C([0,T), \HO{n-1})~, \qquad t \mapsto \gamma(t) :=
\Delta \nu(t) -  \nu(t) \sc \, \partial \nu(t) + \kappa(t) - \dot{\nu}(t) ~;\feq
then from
\rref{nsl} one infers $\LP \gamma(t) = 0$, i.e., $\gamma(t) \in \Hg{n-1}$ for all
$t$. Let
\beq \pi: t \in [0, T) \mapsto \pi(t) := \partial^{-1} \gamma(t) \feq
with $\partial^{-1}$ as in \rref{put}; then $\pi \in C([0,T) H^{n})$
(since $\partial^{-1}$ maps continuously $\Hg{n-1}$ into $H^n$).
One easily checks that $(\nu, \pi)$ fulfills \rref{nsp}. \fine
\section{Appendix. Proof of Eq. \rref{xilip}.}
\label{appexi}
Let us rephrase the definition \rref{defxi} of $\xi$ as
\beq \xi(t) = \zeta(t) - \la \eta(t) \ra~, \feq
having put
\beq \zeta : [0,T) \vain \Hs{n}~, \qquad t \mapsto \zeta(t)~
\mbox{such that}~\zeta(x, t) := \eta(x + h(t), t)~. \label{defze} \feq
Clearly, the thesis \rref{xilip} follows if we prove that
\beq \zeta \in C^{0,1}([0,+\infty), \HM{n-1})~. \label{zelip} \feq
To this purpose, we consider the Fourier coefficients $\eta_{k}(t)$,
$\zeta_k(t)$ of $\eta(t)$, $\zeta(t)$ and note that \rref{defze} implies
\beq \zeta_k(t) = \eta_{k}(t) \, e^{i k \sc h(t)} \qquad (k \in \Zd, t \in [0,T))~. \feq
Let $t, t' \in [0,T)$. We have
\beq \zeta_k(t) - \zeta_k(t') = \alpha_k(t,t') + \beta_k(t,t') \label{ab} \feq
$$ \alpha_k(t,t') := \eta_k(t) e^{i k \sc h(t')} \left(e^{i k \sc (h(t) - h(t'))} - 1 \right)~,~~
\beta_k(t, t') := e^{i k \sc h(t')} \left(\eta_k(t) - \eta_k(t')\right)~. $$
Therefore $\| \zeta(t) - \zeta(t') \|_{n-1} =
\sqrt{ \sum_{k \in \Zd} (1 + |k|^2)^{n-1} |\zeta_k(t) - \zeta_k(t')|^2}$ has the bound
\beq \| \zeta(t) - \zeta(t') \|_{n-1} \leqs A(t, t') + B(t,t')~, \label{AB} \feq
$$ A(t,t') := \sqrt{\sum_{k \in \Zd} (1 + |k|^2)^{n-1} |\alpha_k(t,t')|^2}~,
B(t,t') := \sqrt{\sum_{k \in \Zd} (1 + |k|^2)^{n-1} |\beta_k(t,t')|^2}~. $$
On the other hand, Eq. \rref{ab} and the elementary inequality $|e^{i k \sc y}  - 1|
\leqs | k | | y |$ (for all $y \in \reali^d$) give
\beq | \alpha_k(t,t') | \leqs | k | \, | \eta_k(t) | \, | h(t) - h(t')|~,
\qquad | \beta_k(t,t') | = | \eta_k(t) - \eta_k(t') |~. \feq
Inserting these bounds into the espressions \rref{AB} of $A(t,t')$, $B(t,t')$
(and using $(1 + | k |^2)^{n-1} | k |^2 \leqs (1 + | k |^2)^{n}$) we get
\beq \| \zeta(t) - \zeta(t') \|_{n-1} \leqs \| \eta(t) \|_{n} \, | h(t) - h(t')| +
\| \eta(t) - \eta(t') \|_{n-1}~. \feq
Now, let us consider any compact subset $I$ of $[0,+\infty)$. Then, the assumptions
\rref{assumeta} on $\eta$ and the $C^1$ nature of $h$ ensure the existence of
constants $Q, M_1, M_2$ such that $\| \eta(t) \|_{n} \leqs Q$, $| h(t) - h(t')|
\leqs M_1 | t - t'|$ and $\| \eta(t) - \eta(t') \|_{n-1} \leqs M_2 |t - t'|$ for
$t, t' \in I$. This implies
\beq \| \zeta(t) - \zeta(t') \|_{n-1} \leqs (Q M_1 + M_2) | t - t' |~, \feq
and \rref{zelip} is proved.
\salto
\section{Appendix. Proof of Proposition \ref{edelta}, item (iv).}
\label{aplemen}
Our aim is to derive Eq. \rref{defn}; we write this as
\beq \sup_{t \in [0,+\infty)} \I(t) = \sqrt{2}~, \label{deffn} \feq
\beq \I: [0,+\infty) \vain [0,+\infty)~, \qquad t \mapsto \I(t) := \int_{0}^t d s \,  \uv(t-s) e^{-s}~. \feq
Eq. \rref{deffn} will follow from Steps 1 and 2, giving separate estimates on $\I(t)$ for
$0 \leqs t \leqs 1/4$ and $t > 1/4$.
\parn
\textsl{Step 1. For $0 \leqs t \leqs 1/4$ one has $\I(t) < \sqrt{2}$.}
In fact, with this range for $t$ Eq. \rref{eqww} for $\uv$ implies
$$ \I(t) = {e^{2 t} \over \sqrt{2 \, e}} \int_{0}^{t} d s {e^{-3 s} \over \sqrt{t-s}} \leqs {e^{2 t}
\over \sqrt{2 \, e}} \int_{0}^{t} {d s \over
\sqrt{t-s}} = {e^{2 t} \over \sqrt{2 \, e}} \, 2 \sqrt{t}\leqs {1 \over \sqrt{2}} < \sqrt{2}~. $$
\textsl{Step 2. One has $\sup_{t > 1/4} \I(t) = \sqrt{2}$.} In fact, for all
$t > 1/4$, using again Eq. \rref{eqww} for $\uv$ we get
$$ \I(t) = (\int_{0}^{t-1/4} \! \! \! \! \! \! \! \! d s + \int_{t-1/4}^{t} \! \! \! \! \! d s) \uv(t-s) e^{-s} =
\sqrt{2} \, (1 - e^{1/4 - t}) + C e^{-t} = \sqrt{2} - (\sqrt{2} \, e^{1/4} - C) e^{-t}~, $$
$$ C := \int_{0}^{1/4} d s' \, {e^{3 s'} \over \sqrt{2 \, e s'}}~. $$
One has the estimate $C \leqs 0.6 < \sqrt{2} \, e^{1/4}$, implying $\I(t) <
\sqrt{2}$. From the above expression for
$\I(t)$, we also get $\lim_{t \vain +\infty} \I(t) = \sqrt{2}$; these facts yield the thesis. \fine
\salto
\section{Appendix. The constants $\boma{K_2}$ and $\boma{K_4}$ in dimension $\boma{d=3}$.}
\label{appekdue}
The above constants are needed for the numerical examples in Section
\ref{nume}; the route to compute them is outlined in Lemmas \ref{lemsi}, \ref{lemci} and Proposition
\ref{multi}.
\salto
\textbf{Computing $\boma{K_2}$.} We can take for it any constant such that
\beq \sqrt{\sup_{k \in \Zd_0} \KK_2(k)} \leqs K_2~. \label{mnd22} \feq
We have the bounds \rref{kkn}: $\Ki_2(k) < \KK_2(k) \leqs \Ki_2(k) + \delta \Ki_2(k)$,
with the explicit expressions \rref{mndig2} for $\Ki_2(k)$ and \rref{deci} for
$\delta \Ki_2(k)$ (and $\Ze = \Zd_0$ therein). Both $\Ki_2$ and $\delta \Ki_2$ are defined in terms of some
cutoff function $\Lambda_2$, that we choose in this way: $\Lambda_2(k) := 24$ if $| k | < 4$, and
$\Lambda_2(k) := 6 |k |$ if $|k | \geqs 4$. \parn
The above setting can be employed to evaluate $\KK_2(k)$ for some set of values of $k$;
computations performed for all $k$'s with
$|k_i| \leqs 10$ ($i=1,2,3$) seem to indicate that
\beq \sup_{k \in \Zd_0} \KK_2(k) = \lim_{k \vain \infty} \KK_2(k)~. \label{numtes}  \feq
On the other hand, according to \rref{likn},
\beq \lim_{k \vain \infty} \KK_2(k) = \Sigma_2~, \feq
where $\Sigma_2$ is the series \rref{desig} with $n = 2$ and $\Ze = \Zd_0$. To estimate
this series, we use the bounds \rref{decosig}:
$\Sigm_2(\lambda) < \Sigma_2 \leqs \Sigm_2(\lambda) + \delta \Sigm_2(\lambda)$
with $\Sigm_2(\lambda)$ and $\delta \Sigm_2(\lambda)$ as in Eqs.
\rref{desgm} \rref{reph}; these depend on a cutoff $\lambda$, to be chosen as large as possible
to reach a good precision. Taking $\lambda = 250$, we get
\beq 0.03607 \leqs \Sigma_2 \leqs 0.03934~, \feq
which implies, taking square roots,
\beq 0.1899 \leqs \sqrt{\sup_{k \in \Zd_0} \KK_2(k)} \leqs 0.1984~. \feq
Retaining only two meaningful digits, we take as a final upper bound for $\sqrt{\sup \KK_2}$
the quantity
\beq K_2 := 0.20~. \feq
\textbf{Computing $\boma{K_4}$.} We can take for it any constant such that
$\sqrt{\sup_{k \in \Zd_0} \KK_4(k)} \leqs K_4$.
We use again the bounds \rref{kkn} $\KK_4(k) \leqs \Ki_4(k) + \delta \Ki_4(k)$;
the cutoff $\Lambda_4$ defining $\Ki_4$ and $\delta \Ki_4$ is chosen setting
$\Lambda_4(k) := 10$ if $| k | < 10/3$ and $\Lambda_4(k) := 3 | k |$ if $| k | \geqs 10/3$.
Computing $(\Ki_4 + \delta \Ki_4)(k)$ for $|k_i| \leqs 6$ ($i=1,2,3$)
we obtain numerical evidence that $\sup_{k \in \Zd_0} (\Ki_4 + \delta \Ki_4)(k)$
is attained at $k = (3,0,0)$. So,
\beq \sup_{k \in \Zd_0} \KK_4(k) \leqs (\Ki_4 + \delta \Ki_4)(3,0,0) \leqs 0.004383~. \feq
We also have
\beq \sup_{k \in \Zd_0} \KK_4(k) \geqs \Ki_4(3,0,0) \geqs 0.004382~, \feq
and in conclusion, taking square roots,
\beq 0.06619 \leqs \sqrt{\sup_{k \in \Zd_0} \KK_4(k)} \leqs 0.06621~. \feq
Retaining only two meaningful digits, we take as a final upper bound for $\sqrt{\sup \KK_4}$
the quantity
\beq K_4 := 0.067~. \feq
\salto
\textbf{Acknowledgments.} We are grateful to G. Furioli and E. Terraneo for useful
bibliographical indications.
This work was partly supported by INdAM and by MIUR, PRIN 2006
Research Project "Geometrical methods in the theory of nonlinear waves and applications".

\end{document}